%% file: main.tex
\documentclass[preprint,12pt]{elsarticle}

\usepackage{amssymb}
\usepackage{amsmath} 
\usepackage{algorithmic}
\usepackage{algorithm}
\newlength\myindent
\newcommand\bindent{%
  \begingroup
  \setlength{\itemindent}{\myindent}
  \addtolength{\algorithmicindent}{\myindent}
}
\newcommand\eindent{\endgroup}

\usepackage{graphicx}
\usepackage{cancel}
\usepackage{subcaption}
\usepackage{placeins}
\usepackage{booktabs}
\usepackage{color}
\usepackage[normalem]{ulem}
\usepackage{url}

\journal{Computers \& Fluids}

\begin{document}

\begin{frontmatter}

\title{A Goal-Oriented Adaptive Sampling Procedure for Projection-Based Reduced-Order Models with Hyperreduction}

\author[inst1]{Calista Biondic}

\affiliation[inst1]{organization={Department of Mechanical Engineering, McGill University},
            addressline={817 Rue Sherbrooke Ouest \#270}, 
            city={Montreal},
            postcode={H3A 0C3}, 
            state={QC},
            country={Canada}}

\author[inst1]{Siva Nadarajah}

\begin{abstract}
Projection-based reduced-order models (PROMs) have demonstrated accuracy, reliability, and robustness in approximating high-dimensional, differential equation-based computational models across many applications. For this reason, it has been proposed as a tool for high-querying parametric design problems like those arising in modern aircraft design. Since aerodynamic simulations can be computationally expensive, PROMs offer the potential for more rapid estimations of high-fidelity solutions. However, the efficiency can still be tied to the dimension of the full-order model (FOM), particularly when projected quantities must be frequently recomputed due to non-linearities or parameter dependence. In the case of Petrov-Galerkin models, the projected residual and Jacobian are re-evaluated at every Newton iteration, thereby limiting the anticipated cost improvements. Hyperreduction is one of the tools available to approximate these quantities and address this issue. This work tests the energy-conserving sampling and weighting (ECSW) method as a potential approach for hyperreduction. It will be incorporated into the work in a previous article \cite{donJournal} which had developed an adaptive sampling procedure for building a reduced-order model (ROM) with a controlled functional error. The impacts of hyperreduction on computational cost and accuracy will be studied using the NACA0012 airfoil.
\end{abstract}


\end{frontmatter}



\input{01_introduction}
\input{03_modelRed}
\input{05_adaptSamp}
\input{06_results.tex}
\input{07_conclusion.tex}

\section*{Acknowledgments}
We would like to gratefully acknowledge the financial support of the Natural Sciences and Engineering Research Council of Canada (NSERC) Discovery Grant Program RGPIN-2019-04791.

 \bibliographystyle{elsarticle-num} 
 \bibliography{cas-refs}

\appendix
\input{appendix}

\end{document}

%% file: 01_introduction.tex
\section{Introduction}\label{sec:intro}

Many fields of design and engineering are increasingly reliant on computational models. In aircraft design, for instance, computational fluid dynamics (CFD) has become an essential tool for accelerating analysis and design processes and helps inform experimental testing decisions. Advancements in higher-order methods, faster solution algorithms, and parallel high-performance computing have resulted in access to high-fidelity solutions. These innovations have been instrumental in advancing the design of more fuel-efficient and environmentally friendly aircraft. In certain applications, the high computational costs and storage requirements of these models render them computationally intractable. For example, the time needed to produce high-fidelity CFD solutions poses a significant challenge in high-querying problems, such as aerodynamic shape optimization (ASO), where the model must be evaluated at numerous design parameter combinations, and in online real-time scenarios, like active control, where solutions must be computed within seconds. This limitation has driven substantial interest in methods to reduce the complexity and cost of large-scale computational models. Reduced-order modeling has been proposed as a possible solution to accelerate the process of computing new CFD solutions while preserving the fidelity achieved by higher-order methods.

Reduced-order models can broadly be classified under two types: non-intrusive and intrusive models. Non-intrusive models treat the problem as a black box and are independent of the full-order equations. There are a multitude of challenges that arise when using these models, particularly in large-scale dynamic systems. These include violations of physical constraints, high costs associated with generating large training data sets, and lack of confidence indicators \cite{swischuk2019}. Examples of non-intrusive models include proper orthogonal decomposition coupled with interpolation \cite{Bui-Thanh2003}, Isomap \cite{Franz2014, Franz2016}, Kriging \cite{martin2005, laurenceau2008}, and neural networks \cite{swischuk2019, renganathan2020, lee2020, maulik2021}. Intrusive models operate under the assumption that the solution of a large-scale system of dimension $N$ can be represented on a much lower dimensional subspace of dimension $n$ induced by parameter variation. The low-dimensional representation is found by projecting the governing equations onto the reduced-order subspace, which is why these models are often referred to as projection-based reduced-order models (PROMs). The two most commonly used approaches are the Galerkin projection \cite{holmes2012, sirovich1987, rowley2004model} and the least-squares Petrov-Galerkin projection (LSPG) \cite{carlberg2011, carlberg2017, GRIM2021}. These models are often preferable to non-intrusive models as they are physics-informed, making them more robust, and they offer accessible error metrics \cite{blaisTHESIS}, which help predict the expected accuracy of the model online. Noting that this work focuses on solving fluid flows problems, least-squares Petrov-Galerkin (LSPG) PROMS have been demonstrated to be both numerically stable and accurate in the context of convection-dominated laminar flows and turbulent flow problems, where Galerkin PROMs often exhibit instability \cite{GRIMBERG2020}. Therefore, an LSPG framework will be implemented in this work.

The reduced-order subspace used in the projection is obtained by collecting solution \textit{snapshots} and compressing this information in a reduced-order basis that spans the sampled solution space. Efficient construction of this basis can be a challenge, in particular, appropriate snapshot selection can directly impact the quality of the approximation. Some of the available sampling approaches include uniform sampling, Halton sequence method, Latin hypercube sampling, and greedy sampling \cite{mckay2000, Romero2005, Franz2016, veroy2003, veroy2005, grepl2005}. The shortcomings of these approaches are detailed in \cite{donJournal}. They are compared with a new goal-oriented adaptive sampling approach developed in \cite{donJournal}, \cite{blaisTHESIS}, \cite{AIAA_June2022}, and \cite{BLAIS2023}. The adaptive sampling procedure was designed to reduce the output error of an LSPG projection-based reduced-order model (ROM) to within a prescribed tolerance. The sampling procedure took advantage of dual-weighted residual (DWR) error indicators to estimate the error between the ROM and FOM at specific design parameter locations, which were then used to estimate the output error across the entire parameter space and determine new snapshot locations. It was shown to be more effective at building the reduced-basis than the other sampling approaches previously mentioned. 

The remaining major challenge not addressed in \cite{donJournal} is that the computational cost of the online stage of PROMs still scales with the dimension of the FOM, $N$. In parametric and/or non-linear problems, the parametric vectors and matrices such as the residual and Jacobian must be re-evaluated at each non-linear iteration and then projected onto the reduced-order basis. This can be computationally expensive, and the cost would scale with the dimension of the FOM. If not addressed, the expected time-saving benefits of the PROM cannot be fully realized. There are generally two approaches for addressing this computational bottleneck and approximating these quantities: exact and inexact \cite{FarhatBOOK}. Both share the fundamental strategy that underlies most PROM computations, of breaking down the evaluation of these quantities into offline and online components. Exact methods are applicable to specific classes of problems, two examples being parametric, linear FOMs admitting an efficient parameter-affine representation, and non-linear FOMs with low-order polynomial dependence of the internal force vector on the solution and a time-independent external force vector. These methods are often referred to as exact precomputation-based methods, as they compute the reduced matrices and vectors in two parts: the first part addresses the computational bottleneck and can be precomputed offline, while the second part, whose computational complexity scales with integer powers of the reduced-order subspace $n$, can be processed online and in real-time \cite{FarhatBOOK}.

The problems addressed in this work, however, do not fall into the classes suitable for precomputation. Therefore, inexact or approximate reconstruction methodologies must be used. Hyperreduction is an inexact method that has been developed for addressing more general cases which can be both linear or non-linear, parametric or non-parametric \cite{FarhatBOOK}. These approaches introduce an additional layer of approximation to evaluate the high-dimensional quantities in a manner that is independent of the dimension $N$ of the FOM. As a result, some accuracy in the ROM is traded for computational efficiency. Many methods have been developed in this field, most of which are related in some way, differing primarily by the level of theoretical support~\cite{FARHAT2015}.

Current hyperreduction techniques can be classified into two groups: \textit{approximate-then-project} methods and \textit{project-then-approximate} methods. The former was developed first and, as the name suggests, begins by generating an approximation of the high-dimensional quantity and then computes the exact projection of the approximation onto the left reduced-order basis (ROB) or test basis $\textbf{W}$ associated with the PROM. \cite{FarhatBOOK} notes that the underlying concept of this approach finds its origin in the gappy POD method \cite{everson1995}, which was originally developed for image reconstruction. Some of the notable approximate-then-project hyperreduction methods include the empirical interpolation method (EIM) \cite{BARRAULT2004, GREPL2007} and its discrete version discrete EIM (DEIM) \cite{DEIM, ERROR_DEIM}, the missing point approach \cite{astrid2008}, and the Gauss-Newton with approximate tensors (GNAT) method \cite{carlberg2011, carlberg2013}. These approaches differ in terms of the theoretical support they have, the types of projection frameworks to which they have been applied, and their level of success when applied to CFD problems. EIM was originally derived at the continuous level for PDEs and has some theoretical support for elliptic problems. The discrete version, DEIM, has become one of the most popular approaches. However, both methods have been applied only in Galerkin frameworks and have seen limited use in aerodynamics. When used, such as in Carlberg et al. \cite{carlberg2013, carlberg2017}, they have been reported to exhibit temporal instabilities for turbulent unsteady flows \cite{yanoBook}. The GNAT method, on the other hand, was formulated at the discrete level for Petrov-Galerkin frameworks and designed around approximations that satisfy consistency and discrete-optimality conditions \cite{carlberg2013}. It has proven highly effective for nonlinear structural dynamics problems \cite{carlberg2011}, demonstrating robustness, accuracy, and excellent computational cost savings, as well as for benchmark turbulent viscous flow problems \cite{carlberg2013}. However, Washabaugh \cite{washabaughThesis} notes that the GNAT method can be expensive for parameterized steady aerodynamic problems. Given the limitations of the approximate-then-project method, these approaches are not considered in this work.

The second group, project-then-approximate methods, was proposed more recently as a response to concerns about the stability of hyperreduction techniques. These methods have demonstrated greater robustness in certain cases \cite{FarhatBOOK}. Project-then-approximate methods first project the full-order model quantities onto the left ROB $\mathbf{W}$, then approximate the reduced-order vectors and matrices. These methods typically decompose the reduced-order quantities into a summation over the elements in the computational mesh and then seek to approximate the quantities by including only the contributions from a subset of the mesh elements, referred to as the reduced mesh. The way this set is determined varies, but all of the methods can be interpreted as generalized quadrature rules, where an empirical set of training data is collected, and a set of quadrature ``points" and associated weights are learned in a supervised procedure \cite{GRIM2021}. All project-then-approximate methods conduct some form of mesh sampling in which the quadrature points are mesh elements. Some examples of this class of approaches include the empirical quadrature procedure (EQP) \cite{PATERA2017, patera2019}, the energy-conserving sampling and weighting (ECSW) method \cite{FARHAT2014, FARHAT2015}, and the empirical cubature method (ECM) \cite{HERNANDEZ2017}. The EQP method is unique in that it can control the error introduced by hyperreduction in the norm of the solution \cite{patera2019}, and was later extended to control the error in the quantity of interest \cite{yano2020}. This differs from other hyperreduction methods that control the error introduced into the quantity being approximated, such as the residual, which must then be ``tuned” to achieve the desired functional error tolerance \cite{SLEEMAN2022}. However, both the EQP and ECM methods have only been applied in Galerkin frameworks. While the ECSW was originally developed for Galerkin frameworks, it was later extended to use in Petrov-Galerkin frameworks in \cite{GRIM2021}. For this reason, the ECSW hyperreduction method emerges as a suitable choice for the goals of this work. Moreover, it has both a substantial amount of theoretical and experimental support. The ECSW method is provably unconditionally stable for second-order hyperbolic models and has been shown to be both numerically stable and accurate in structural dynamics problems, where approximate-then-project methods have been known to fail \cite{FARHAT2015}. Additionally, when applied to first-order hyperbolic problems, it was shown to be reliable, accurate, and computationally efficient \cite{GRIM2021}.

The goal of this work is to advance the goal-oriented adaptive sampling procedure developed in \cite{donJournal}. We build a hyperreduced reduced-order model (HROM) which still achieves a prescribed error tolerance across the parameter space while accounting for the additional error introduced onto the output functional of interest by hyperreduction. The hyperreduction is implemented using the ECSW approach described in the previous section and is employed to approximate both the residual and the Jacobian in the HROM. Additionally, the second DWR error indicator introduced in \cite{donJournal}, which quantifies the error between a coarse and fine ROM (referring to the number of columns in the reduced order basis), is modified to capture the error introduced by the hyperreduction. The current work is only interested in steady state solutions. Unsteady CFD solutions are not addressed and are left as a subject for future work.

%% file: 03_modelRed.tex
\section{Model Reduction}\label{sec:model_red}

This section discusses the various components required to construct a hyperrreduced reduced-order model (HROM). These include the creation of a trial basis , or reduced-order basis (ROB) $\mathbf{V}$, via proper orthogonal decomposition, onto which the full-order model (FOM) will be projected. The basis is then used to derive an approximate representation $\tilde{\mathbf{w}}$ for a full-order CFD solution $\mathbf{w}$. By projecting the system of nonlinear equations onto the test basis $\mathbf{W}$, Newton's method can be employed to solve for the unknown reduced-order solution representation $\hat{\mathbf{w}}$. The least-squares Petrov-Galerkin (LSPG) framework is introduced. Finally, a detailed description of the chosen hyperreduction approach, the energy-conserving sampling and weighting (ECSW) approach, is given.

\subsection{Proper Orthogonal Decomposition}

As mentioned in the introduction, projection-based reduced order models take advantage of the lower dimensionality of the FOM solution manifold. They captures the physics of the problem by projecting the FOM onto a subspace of a smaller dimension, typically using a properly trained reduced-order basis $\mathbf{V} \in \mathbb{R}^{N \times n}$ \cite{TEZAUR2022}. This basis is a matrix whose columns, often referred to as basis vectors or modes, span the lower dimensional subspace \cite{blaisTHESIS}.  In fluid dynamics, the expansion procedure conducted to build the basis is referred to as proper orthogonal decomposition (POD). POD can extract dominant features in data or decompose a function into its underlying modes, which can then be used to build a reduced-order basis (ROB). Further details on the POD procedure and its derivation can be found in \cite{blaisTHESIS}.

The most common approach for identifying POD modes is singular value decomposition (SVD), a concept from linear algebra. To construct the ROB, a set of solution samples $\mathbf{w}^s \in \mathbb{R}^N$ at various parameter combinations must be collected. The way in which these samples are distributed in the parameter domain will be discussed further in Section \ref{sec:adapt}. Using a reference solution $\mathbf{w}_{\text{ref}} \in \mathbb{R}^N$, which in the previous work \cite{blaisTHESIS} was chosen to be the mean of the solution samples, the snapshots used to determine the basis vectors are of the form \( \mathbf{s} = \{\mathbf{w^s} - \mathbf{w}_{\text{ref}}\}_{s=1}^{n}\) to ensure consistency. These can be assembled into a matrix $\mathbf{S}$:
\begin{equation} \label{eq:snap}
    \mathbf{S} = \begin{bmatrix} \mathbf{s}_1 & \mathbf{s}_2 & \hdots & \mathbf{s}_n  \end{bmatrix} \in \mathbb{R}^{N \times n}.
\end{equation}
The SVD of the above matrix would then be given by:
\begin{equation}
\mathbf{S} = \mathbf{U}\mathbf{\Sigma}\mathbf{Z}^T ,
\end{equation}
where $\mathbf{U} \in \mathbb{R}^{N \times N}$, $\mathbf{Z} \in \mathbb{R}^{n \times n}$ are orthogonal matrices and \\$\mathbf{\Sigma} = \text{diag}(\sigma_1, \dots, \sigma_k,0,\dots, 0)$ is an $N \times n$ diagonal matrix containing singular values $\sigma$, where $k = \text{min}(N,n)$. The columns of $\mathbf{U}$ are equivalent to the POD modes identified using an eigenvalue problem. In this work, the number of snapshots $n$ is assumed to be much smaller than the dimension of the FOM $N$. As a result, the last $N-n$ columns of $\mathbf{U}$ are arbitrary, as they correspond to zero singular values. Therefore, by employing a ``thin" SVD without truncation, the first $n$ columns of $\mathbf{U}$ are selected as the basis vectors for the ROB $\mathbf{V}$ which spans the reduced-order subspace, commonly referred to as the trial space.

\subsection{Projection-Based Reduced-Order Model}

The current work is concerned with steady-state solutions to PDEs. The semi-discretization is accomplished using a discontinuous Galerkin method \cite{philip_2022}. The steady-state problem is formulated as:
\begin{equation}
\begin{split} \label{eqn:FOM}
& \mathbf{R}(\mathbf{w}_h, \boldsymbol{\mu}) = 0, \\
& \mathcal{J} = \mathcal{J} (\mathbf{w}_h, \boldsymbol{\mu}).
\end{split}
\end{equation}
where $\mathbf{w}_h \in \mathbb{R}^N$  is the discrete numerical solution approximating the exact solution $\mathbf{w}$, $\mathbf{R} \in \mathbb{R}^N$ is the residual, $\boldsymbol{\mu}$ is a set of input parameters, and $\mathcal{J}$  is an output functional of interest.

Once a basis $\mathbf{V} \in \mathbb{R}^{N \times n}$ has been determined using POD, the next step is to find an approximation of the FOM CFD solution $\mathbf{w} \in \mathbb{R}^{N}$. Projection-based ROMs operate under the assumption that these solutions lie in the reduced-order subspace. At a given parameter location $\boldsymbol{\mu}$, we would like to then find an approximate solution $\tilde{\mathbf{w}} \in \mathbb{R}^N$ of the form:
\begin{equation}
\label{eqn:rom_soln}
    \tilde{\mathbf{w}} = \mathbf{w}_{\text{ref}} + \mathbf{V}\hat{\mathbf{w}},
\end{equation}
where $\mathbf{w}_{\text{ref}} \in \mathbb{R}^N$ is the reference state and $\hat{\mathbf{w}} \in \mathbb{R}^n$ is the unknown reduced-order solution. Substituting the approximate solution in Equation \ref{eqn:rom_soln} into Equation \ref{eqn:FOM} results in an overdetermined system as there are $N$ equations and $n$ unknowns in $\hat{\mathbf{w}}$. The system of equations is then projected onto a test subspace via the test basis $\mathbf{W} \in \mathbb{R}^{N \times n}$. This results in the following reduced-order system of nonlinear equations:
\begin{equation}
 \mathbf{W}^T \mathbf{R}(\mathbf{w}_{\text{ref}} + \mathbf{V}\hat{\mathbf{w}}, \boldsymbol{\mu}) = \hat{\mathbf{R}}(\mathbf{w}_{\text{ref}} + \mathbf{V}\hat{\mathbf{w}}, \boldsymbol{\mu}) = \mathbf{0}^T,
\end{equation}
where $\hat{\mathbf{R}} \in \mathbb{R}^n$ is the reduced-order residual. Using a Taylor series expansion to approximate the solution, we have:
\begin{equation} 
\hat{\mathbf{R}}(\mathbf{w}_{\text{ref}} + \mathbf{V}\hat{\mathbf{w}}^{(k+1)},\boldsymbol{\mu}) \approx \hat{\mathbf{R}}(\mathbf{w}_{\text{ref}} + \mathbf{V}\hat{\mathbf{w}}^{(k)}, \boldsymbol{\mu}) + \frac{\partial \hat{\mathbf{R}}}{\partial \hat{\mathbf{w}}}^{(k)} (\hat{\mathbf{w}}^{(k+1)} - \hat{\mathbf{w}}^{(k)}),
\end{equation}
and given that $\hat{\mathbf{R}}(\mathbf{w}_{\text{ref}} + \mathbf{V}\hat{\mathbf{w}}, \boldsymbol{\mu}) = \mathbf{0}^T$, then:
\begin{equation} 
\mathbf{0}^T = \mathbf{W}^T \mathbf{R}(\mathbf{w}_{\text{ref}} + \mathbf{V}\hat{\mathbf{w}}^{(k)}, \boldsymbol{\mu}) + \frac{\partial(\mathbf{W}^T \mathbf{R}(\mathbf{w}_{\text{ref}} + \mathbf{V}\hat{\mathbf{w}}^{(k)}, \boldsymbol{\mu})}{\partial \hat{\mathbf{w}}} (\hat{\mathbf{w}}^{(k+1)} - \hat{\mathbf{w}}^{(k)}).
\end{equation}
Ignoring the second-order sensitivities created by taking the derivative of the test basis and using the chain rule on the remaining derivative of the residual with respect to $\hat{\mathbf{w}}$, we find:
\begin{equation}
\mathbf{W}^T\frac{\partial \mathbf{R}}{\partial {\mathbf{w}}}^{(k)} \mathbf{V} (\hat{\mathbf{w}}^{(k+1)} - \hat{\mathbf{w}}^{(k)})  = -\mathbf{W}^T \mathbf{R}(\mathbf{w}_{\text{ref}} + \mathbf{V}\hat{\mathbf{w}}^{(k)}, \boldsymbol{\mu}) .
\end{equation}
Therefore, the Newton iteration for $k = 1, ..., K$ are:
\begin{equation}
\begin{split}
& \left[\mathbf{W}^T \frac{\partial \mathbf{R}}{\partial \mathbf{w}}^{(k)} \mathbf{V} \right] \mathbf{p}^{(k)} = - \mathbf{W}^T \mathbf{R}^{(k)}, \\ 
& \hat{\mathbf{w}}^{(k+1)} = \hat{\mathbf{w}}^{(k)} + \mathbf{r}^{(k)} \mathbf{p}^{(k)},
\end{split}
\label{eq:gen_newton}
\end{equation}
where $\mathbf{p}^{(k)}$ is the search direction and $\mathbf{r}^{(k)}$ the step length.

\subsection{Least-squares Petrov-Galerkin Projection}

The final component of the PROM framework that must be selected is the test basis $\mathbf{W}$. Two options are available, the Galerkin and Petrov-Galerkin methods. Petrov-Galerkin (PG) PROMs have been shown to be ideal for parametric, nonlinear high-dimensional models like CFD problems. While it was previously thought that using PROMs in fluid dynamics would lead to numerical instabilities due to their inability to resolve the dissipative regime of the turbulent energy cascade, \cite{GRIMBERG2020} demonstrated that in fact, the Galerkin framework was the true source of these instabilities. Specifically, for problems where the Jacobian is not symmetric positive definite, which is generally the case in CFD problems, PROMs constructed using the Galerkin framework can perform poorly or exhibit instabilities \cite{rowley2004model, GRIMBERG2020}. However, for the same convection-dominated laminar and turbulent flows, PROMs constructed using the Petrov-Galerkin (PG) framework were found to be numerically stable and accurate \cite{GRIMBERG2020}.

Given that this work focuses CFD problems, the Petrov-Galerkin framework is the better fit. Specifically, the least-squares Petrov-Galerkin (LSPG) projection offers several attractive properties, including optimality \cite{carlberg2011} and monotonicity \cite{zahr2016}, which are discussed in detail in \cite{blaisTHESIS} and motivate its selection for use in this approach.

In the case of the least-squares Petrov-Galerkin projection \cite{benner2020, washabaughThesis}, the test basis is constructed such that $\mathbf{W} \equiv \frac{\partial \mathbf{R}}{\partial \mathbf{w}} \mathbf{V}$. Recalling Equation \ref{eq:gen_newton}, the Newton iterations for $k = 1, ..., K$ in the LSPG framework become:
\begin{equation}
\begin{split}
& \left[\mathbf{V}^T \frac{\partial \mathbf{R}}{\partial \mathbf{w}}^{(k)T} \frac{\partial \mathbf{R}}{\partial \mathbf{w}}^{(k)} \mathbf{V} \right] \mathbf{p}^{(k)} = - \mathbf{V}^T \frac{\partial \mathbf{R}}{\partial \mathbf{w}}^{(k) T} \mathbf{R}^{(k)}, \\ 
& \hat{\mathbf{w}}^{(k+1)} = \hat{\mathbf{w}}^{(k)} + \mathbf{r}^{(k)} \mathbf{p}^{(k)},
\end{split}
\label{eq:lspg_newton}
\end{equation}
where $\frac{\partial \mathbf{R}}{\partial \mathbf{w}}^{(k)} \in \mathbb{R}^{N \times N}$ is the Jacobian, $\mathbf{p}^{(k)}$ is the search direction and $\mathbf{r}^{(k)}$ the step length. It can be seen here that the test basis updates with each Newton iteration via the addition of the Jacobian, allowing it to better capture nonlinear effects \cite{carlberg2011}.

\subsection{Hyperreduction}

In Newton's method shown in Equation \ref{eq:lspg_newton} above, it can be seen that the residual and Jacobian must be re-computed at each iteration. This is one of the computational bottlenecks for a ROM of a non-linear, parametric model, as the cost of recomputing and projecting these quantities would scale with the dimension of the FOM $N$. Noting that the right-hand side of the first line in Equation \ref{eq:lspg_newton} is the reduced-order residual $\hat{\mathbf{R}}$, each iteration requires the evaluation of the high-dimensional residual $\mathbf{R}$, followed by multiplication with the test basis $\mathbf{W}^T$, which has asymptotic complexity $O(Nn)$ \cite{GRIMAIAA}. Therefore, while ROMs of the above form can provide accurate approximations of the FOM, they may not be significantly less computationally expensive. For the class of problems discussed in this work, an inexact methodology is required for addressing this computational bottleneck. Specifically, the ECSW hyperreduction method will be used, as discussed in the introduction.

\subsubsection{Energy-conserving Sampling and Weighting Method}

Farhat et al. \cite{FARHAT2014} note that previous hyperreduction approaches have primarily focused on the accuracy of the approximation of the FOM quantities, but have given little consideration to important properties of the resulting HROM, such as preservation of symmetry or numerical stability \cite{FARHAT2014}. The ECSW method was developed to preserve both symmetry and stability using the concepts of mesh sampling and the principle of virtual work \cite{FARHAT2014}. It has been shown that, in second-order hyperbolic problems, this approach preserves the Lagrangian structure associated with Hamilton's principle, which in turn enables the preservation of the numerical stability properties of the discrete system \cite{FARHAT2015}. When applied to realistic structural dynamics problems, the use of the ECSW method results in stable and accurate HROMs, whereas HROMs built with DEIM fail due to numerical instability \cite{FARHAT2015}. \cite{GRIM2021} extended the ECSW method to Petrov-Galerkin PROMs and demonstrated its success in constructing robust, accurate, and computationally efficient models for CFD applications, particularly those associated with convention-dominated viscous flows. Due to these properties, the ECSW method is the most suitable candidate for hyperreduction in this work.

Consider the discretization of a spatial domain into $N_e$ mesh entities making up the set $\mathcal{E}$, in the case of finite difference semi-discretization these would be vertices. The reduced-order residual for a steady-state problem can be  written as:
\begin{equation}
\label{eq:elem_res}
    \hat{\mathbf{R}}(\hat{\mathbf{w}};\boldsymbol{\mu}) = \sum_{e\in\mathcal{E}}\mathbf{W}^T\mathbf{L}_e^T\mathbf{R}_e(\mathbf{L}_{e^+}(\mathbf{w}_{\text{ref}} + \mathbf{V}\hat{\mathbf{w}});\boldsymbol{\mu}),
\end{equation}
where $\mathbf{L}_e \in \{0,1\}^{d_e \times N}$ is a boolean matrix localizing a FOM vector to the $d_e$ degrees of freedom (DOFs) associated with the mesh element $e$. $\mathbf{R}_e \in \mathbb{R}^{d_e}$ is the contribution of this element to the global FOM residual; the spatial stencil determines the set of entities $n_{e^+}$ required to evaluate $\mathbf{R}_e$, which can include the neighbours of $e$. Similar to $\mathbf{L}_e$, $\mathbf{L}_{e^+} \in \{0,1\}^{(d_e n_{e^+}) \times N}$ is the boolean matrix that localizes a FOM vector to the DOFs associated with the entities $n_{e^+}$. It is assumed that there exists a subset of all the mesh entities $\widetilde{\mathcal{E}} \subset \mathcal{E}$, with size $\widetilde{N}_e = |\widetilde{\mathcal{E}}| \ll N_e$, such that the reduced residual can be approximated with a smaller number of entities. The hyperreduced residual vector $\widetilde{\mathbf{R}}$ can be written as:
\begin{equation}
\label{eq:hyp_res}
    \hat{\mathbf{R}}(\hat{\mathbf{w}};\boldsymbol{\mu}) \approx \widetilde{\mathbf{R}} (\hat{\mathbf{w}};\boldsymbol{\mu}) = \sum_{e\in\widetilde{\mathcal{E}}}\xi_e\mathbf{W}^T\mathbf{L}_e^T\mathbf{R}_e(\mathbf{L}_{e^+}(\mathbf{w}_{\text{ref}} + \mathbf{V}\hat{\mathbf{w}});\boldsymbol{\mu}),
\end{equation}
and interpreted as a generalized quadrature rule with a set of mesh element weights $\{\xi_e \mid e \in \widetilde{\mathcal{E}}\}$. Note that this is an approximation of the projected residual $\hat{\mathbf{R}}$ not the FOM residual $\mathbf{R}$, therefore $\widetilde{\mathbf{R}} \in \mathbb{R}^{n}$. With no additional approximation, the hyperreduced Jacobian $\mathbf{\widetilde{J}}$ can be written as:
\begin{equation}
\label{eq:hyp_jac}
    \frac{\partial \hat{\mathbf{R}}}{\partial \mathbf{w}} (\hat{\mathbf{w}};\boldsymbol{\mu}) \approx \mathbf{\widetilde{J}} = \sum_{e\in\widetilde{\mathcal{E}}}\xi_e\mathbf{W}^T\mathbf{L}_e^T\mathbf{J_e}(\mathbf{L}_{e^+}(\mathbf{w}_{\text{ref}} + \mathbf{V}\hat{\mathbf{w}});\boldsymbol{\mu})\mathbf{L}_{e^+}\mathbf{V},
\end{equation}
where $\mathbf{J_e} \in \mathbb{R}^{d_e \times d_e^+}$ is the Jacobian matrix of $ \mathbf{R}_e$ with respect to the DOFs associated with the $n_{e^+}$ entities involved in its evaluation \cite{GRIMAIAA}.

\subsubsection{Reduced Mesh Set Selection}

\cite{TEZAUR2022} observed that, particularly in steady-state problems, residual values can often be very close to zero. This can result in poor training data and lead to inaccuracies in the reduced mesh and the associated weights. The challenges of applying the residual-based training strategy to steady-state problems are discussed in detail in \cite{TEZAUR2022}. A summary of these challenges is provided below:

The hyperreduction will be conducted in two stages: first, the offline computation of the reduced mesh and accompanying quadrature weights, followed by the online use of these approximations to evaluate the projected quantities within the LSPG framework. Two ECSW approaches were compared in the offline computation stage of this work. One approach trains the weights using the residual of a subset of the converged solution snapshots used to build the ROB, as in \cite{GRIM2021} and \cite{FARHAT2014}. The other approach uses the Jacobian of the snapshots instead. This was motivated by our own initial results, which identified challenges with the conditioning of the matrices produced by the residual approach. \cite{TEZAUR2022} observed that, particularly in steady-state problems, the residual values can often be very close to zero. This can result in poor training data and lead to inaccuracies in the reduced mesh and the associated weights. The challenges of applying the residual-based training strategy to steady-state problems are discussed in detail in \cite{TEZAUR2022}. A shortened list is provided here:
\begin{itemize}
    \item If the ROB is selected without truncation of the singular value decomposition, the snapshot solution will have an exact approximate representation in the reduced-order space. As a result, the residuals—and thus the training data—would nearly vanish for these computed solutions, which would be equivalent to training using numerical noise. This will likely lead to inaccurate approximations. Furthermore, \cite{TEZAUR2022} notes that this scenario would invalidate the early stopping criterion in Equation \ref{eq:opt_NNLS} (introduced later in the non-negative least squares (NNLS) problem) that is used to build the reduced mesh. 
    \item \cite{TEZAUR2022} also explains that parametric steady-state problems generally require fewer solution snapshots than unsteady ones, meaning that the available training data is more limited which can negatively impact the online performance of the hyperreduction approximations. This could create a challenge in the adaptive sampling procedure, as the purpose of the approach is to begin with as few snapshots as possible. If hyperreduction is then incorporated, there will be minimal training data available at earlier iterations. This could make the hyperreduction inaccurate, causing the error indicators at the ROM points to not capture the true error distribution, thereby affecting the placement of new FOM snapshots.
\end{itemize}

The next two sections outline the residual-based and Jacobian-based approaches for assembling the training data to determine the reduced mesh and associated weights. The final section describes the optimization problem arising from the training data.

\subsubsection{Residual-based ECSW Training Process}
To begin, each solution snapshot used to find the right ROB $\mathbf{V}$ in Equation \ref{eq:snap} can be represented on the lower dimension with the following projection \cite{GRIM2021}:
\begin{equation}
\label{eq:rom_y}
    \hat{\mathbf{w}}_s = \mathbf{V}^T(\mathbf{w}^s - \mathbf{w}_{\text{ref}}),
\end{equation}
which can then be converted into an approximation of the FOM vector of dimension $N$ through the equation:
\begin{equation}\begin{split}
\label{eq:rom_approx}
    \tilde{\mathbf{w}}_s &= \mathbf{w}_{\text{ref}} + \mathbf{V}\hat{\mathbf{w}}_s ,\\
    &= \mathbf{w}_{\text{ref}} + \mathbf{V}\mathbf{V}^T(\mathbf{w}^s - \mathbf{w}_{\text{ref}}).
    \end{split}
\end{equation}
The original set of snapshots $\mathcal{S}$ contains $n$ FOM solutions as shown in the matrix in Equation \ref{eq:snap}. We create a new set $\mathcal{S}_H$ of size $N_s$ which contains all or a subset of the snapshots in $S$, i.e. $\mathcal{S}_H \subseteq \mathcal{S}$. For each snapshot in $\mathcal{S}_H$, the contributions from each mesh element to the discrete ROM residual can be broken into:
\begin{equation}
\label{eq:ecsw}
\begin{split}
    \mathbf{c}_{se} &= \mathbf{W}^T\mathbf{L}_e^T\mathbf{R}_e(\mathbf{L}_{e^+}(\mathbf{w}_{\text{ref}}+\mathbf{V}\mathbf{V}^T(\mathbf{w}^s - \mathbf{w}_{\text{ref}}));\boldsymbol{\mu}) \in \mathbb{R}^n, \\
    \mathbf{d}_s &= \sum_{e\in\widetilde{\mathcal{E}}} c_{se}\in \mathbb{R}^n \;\;\; s = 1, \dots , N_s \;\;\; e = 1, \dots, N_e.
    \end{split}
\end{equation}
The training data can then be organized in block form and the exact assembly of this data will be as follows:
\begin{equation}\label{eq:CD}
\mathbf{C} \mathbf{1} = \mathbf{d},
\end{equation}
where,
\begin{equation}
\label{eq:ecsw_matrix}
\mathbf{C} = \begin{bmatrix}
c_{11}        & \dots & c_{1N_{e}} \\
\vdots  & \ddots & \vdots \\
c_{N_{s}1}       & \dots & c_{N_{s}N_{e}}
\end{bmatrix} \in \mathbb{R}^{(N_s n) \times N_e} \hspace{1cm} \text{and} \hspace{1cm} \mathbf{d} = \begin{bmatrix}
d_{1}         \\
\vdots  \\
d_{N_{s}}
\end{bmatrix} \in \mathbb{R}^{(N_s n)}
\end{equation}
and $\mathbf{1}$ is a $N_e$ length vector of ones. The above matrix and vector will then be fed into the NNLS problem to solve for the weights and mesh entities used for the hyperreduction approximation, which will be described in more detail in an upcoming section.

\subsubsection{Jacobian-based ECSW Training Process}
Similar to the residual approach, a set $\mathcal{S}_H$ of size $N_s$ is created containing all or a subset of the snapshots in $S$, i.e. $\mathcal{S}_H \subseteq \mathcal{S}$. For each snapshot in $\mathcal{S}_H$, the contributions from each mesh element to the discrete ROM Jacobian can be broken into:
\begin{equation}
\label{eq:ecsw_jac}
\begin{split}
    \mathbf{c}_{se} &= \mathbf{Q}\mathbf{W}^T\mathbf{L}_e^T\mathbf{J_e}(\mathbf{L}_{e^+}(\mathbf{w}_{\text{ref}} + \mathbf{V}\hat{\mathbf{w}});\boldsymbol{\mu})\mathbf{L}_{e^+}\mathbf{V} \in \mathbb{R}^{n^2} ,\\
    \mathbf{d}_s &= \sum_{e\in\widetilde{\mathcal{E}}} c_{se}\in \mathbb{R}^{n^2}, \;\;\; s = 1, \dots , N_s, \;\;\; e = 1, \dots, N_e, 
    \end{split}
\end{equation}
where $\mathbf{Q}$ is the matrix that converts an $n \times n$ matrix into a column vector of dimension $n^2$ by stacking its columns on top of each other. The exact assembly of the training data can be written as:
\begin{equation}\label{eq:CD_jac}
\mathbf{C} \mathbf{1} = \mathbf{d},
\end{equation}
where,
\begin{equation}
\label{eq:ecsw_mat_jac}
\mathbf{C} = \begin{bmatrix}
c_{11}        & \dots & c_{1N_{e}} \\
\vdots  & \ddots & \vdots \\
c_{N_{s}1}       & \dots & c_{N_{s}N_{e}}
\end{bmatrix} \in \mathbb{R}^{(N_s n^2) \times N_e} \hspace{1cm} \text{and} \hspace{1cm} \mathbf{d} = \begin{bmatrix}
d_{1}         \\
\vdots  \\
d_{N_{s}}
\end{bmatrix} \in \mathbb{R}^{(N_s n^2)}
\end{equation}
and $\mathbf{1}$ is a $N_e$ length vector of ones. It should be noted that while the number of rows in the $\mathbf{C}$ matrix scales with the dimension $n$ of the ROB $\mathbf{V}$ for residual-based training, for the Jacobian-based training it scales with $n^2$. This can lead to a significant increase in the training cost, which in this paper will be controlled by adjusting the number of snapshots in the set $\mathcal{S}_H$. Other methods for addressing this cost are also discussed in \cite{TEZAUR2022}.

\subsubsection{Non-Negative Least-Squares Problem}
Once the training data has been assembled using one of the two approaches above, the goal is then to find a set of weights and accompanying mesh entities that can reproduce Equation \ref{eq:CD} or \ref{eq:CD_jac} to a specified level of accuracy. Ideally, the vector of weights $\boldsymbol{\xi}$ replacing the vector of ones would be sparse, meaning the number of elements required to approximate the residual is much smaller than the total number of entities $N_e$. The result is an optimization problem of the following form \cite{GRIMAIAA}:
\begin{equation}
\label{eq:opt}
\begin{split}
\text{minimize} \; &\lVert \boldsymbol{\xi} \rVert_0 ,\\
\text{subject to} \; &\lVert \mathbf{C}\boldsymbol{\xi} -\mathbf{d}\rVert_2 \leq \epsilon \lVert \mathbf{d}\rVert_2 ,\\
& \boldsymbol{\xi} \geq \mathbf{0}^T,
\end{split}
\end{equation}
where $\epsilon$ is an error tolerance used to control the accuracy of the hyperreduction. However, this is known to be an NP-hard problem. Using the results from \cite{CHAP2016}, the non-convex problem above is replaced by a convex approximation known as a non-negative least-squares (NNLS) problem:
\begin{equation}
\label{eq:opt_NNLS}
\begin{split}
\text{minimize} \; &\lVert \mathbf{C}\boldsymbol{\xi} -\mathbf{d}\rVert_2^2 , \\
\text{subject to} \; & \boldsymbol{\xi} \geq \mathbf{0}^T ,\\
\text{with early stopping criteria} \; &\lVert \mathbf{C}\boldsymbol{\xi} -\mathbf{d}\rVert_2 \leq \epsilon \lVert \mathbf{d}\rVert_2  \text{ and } \boldsymbol{\xi} \geq \mathbf{0}^T .
\end{split}
\end{equation}
Note, in the NNLS problem, the closer $\epsilon$ is to unity, the more the solution approaches the results from the optimization problem in Equation \ref{eq:opt} and the more sparse the resulting weights vector is. The closer $\epsilon$ is to zero the more accurate the approximation is \cite{CHAP2016}. For completeness, Algorithm 1 summarizes the full procedure to find the weights and reduced mesh for the ECSW residual-based approach. The only difference for the Jacobian-based approach would be the entries of $\mathbf{c}_{se}$ and $\mathbf{d}_s$.

\begin{algorithm}[H]
\caption{Finding Weights for ECSW (Residual-based Training Data)}\label{alg:cap}
\begin{algorithmic}
\STATE{\textbf{Inputs:}}
\STATE{Set of Mesh Entities $\mathcal{E} = \{ e_1, e_2, \dots e_{n_e} \}$}
\STATE{POD Basis $\mathbf{V}$, Test Basis $\mathbf{W}$}
\STATE{Reference state $\mathbf{w_{\text{ref}}}$}
\STATE{Subset $S_H$ of snapshots $S$ (which were used to build the POD Basis)}
\STATE{NNLS Tolerance $\epsilon_{\text{NNLS}}$}
\STATE{\textbf{Outputs:}}
\STATE{Reduced Mesh Element Set $\widetilde{\mathcal{E}}$}
\STATE{Associate Set of Weights $\boldsymbol{\xi}_{\mathcal{E}} = \{\xi_{e} | e \in \widetilde{\mathcal{E}}\}$}
\STATE{-------------------------------------------------}

\FOR{$e \in \mathcal{E}$}
        \STATE{Construct $\mathbf{L}_e \in \{0,1\}^{d_e \times N}$ and $\mathbf{L}_{e^+} \in \{0,1\}^{(d_e n_{e^+}) \times N}$ }
        \FOR{$\mathbf{w}^s \in S_H$}
            \STATE{Let $\mathbf{c}_{se} = \mathbf{W}^T\mathbf{L}_e^T\mathbf{R}_e(\mathbf{L}_{e^+}(\mathbf{w}_{\text{ref}}+\mathbf{V}\mathbf{V}^T(\mathbf{w}^s - \mathbf{w}_{\text{ref}}));\boldsymbol{\mu})$, where $\mathbb{R}^{d_e}$ is the contribution of this element to the global FOM residual}
            \STATE{Let $\mathbf{d}_s = \sum_{e\in\widetilde{\mathcal{E}}}\mathbf{c}_{se}$}
        \ENDFOR
\ENDFOR
\STATE{Solve the NNLS optimization problem for $\widetilde{\mathcal{E}}$ and $\boldsymbol{\xi}_{\mathcal{E}}$:}
    \bindent
    \STATE{minimize $\lVert \mathbf{C}\boldsymbol{\xi} -\mathbf{d}\rVert_2^2$}
    \STATE{ subject to $\boldsymbol{\xi}  \geq \mathbf{0}^T$}
    \STATE{ with early termination criterion $\lVert \mathbf{C}\boldsymbol{\xi} -\mathbf{d}\rVert_2 \leq \epsilon_{\text{NNLS}} \lVert \mathbf{d} \rVert_2$}
    \eindent
\RETURN{$\boldsymbol{\xi}_{\mathcal{E}} = \{\xi_{e} | e \in \widetilde{\mathcal{E}}\}$}
\end{algorithmic}
\end{algorithm}

\subsection{Hyperreduced ROM Solution Evaluation with ECSW}

With the resulting set of weights and mesh elements, the total number of elements on which the residual needs to be computed to approximate its reduced-order representation is expected to be much smaller than the total number of mesh elements $N_e$. The ECSW results can then be applied within the Newton iterations of the LSPG projection shown in Equation \ref{eq:lspg_newton}. At each new iteration, both the test basis $\mathbf{W}$, which depends on the Jacobian, and the residual $\mathbf{R}$ must be re-evaluated. The process begins by solving for the test basis. Using the weights and mesh entities derived from the NNLS optimization solution, the Jacobian can then be approximated with a reduced number of mesh entities. For the FOM, the Jacobian at a specific Newton iteration $k$ can be broken down into the contributions from each element:
\begin{equation}
    \mathbf{J}^{(k)} = \frac{\partial \textbf{R}}{\partial \textbf{w}}^{(k)} = \sum_{e \in \mathcal{E}} \mathbf{L}_e^{T} \mathbf{J}_e^{(k)} \mathbf{L}_{e^+}, 
 \end{equation}
where $\mathbf{J}_e^{(k)}$ is the contribution of element $e$  to the global, FOM Jacobian for the solution $\tilde{\mathbf{w}}^{(k)}$. Therefore, the hyperreduced FOM Jacobian would be represented by:
\begin{equation}
    \mathbf{J}^{(k)} \approx \bar{\mathbf{J}}^{(k)} =  \sum_{e \in \tilde{\mathcal{E}}} \xi_e \mathbf{L}_e^{T} \mathbf{J}_e^{(k)} \mathbf{L}_{e^+}.  \label{hyp_jac}
 \end{equation}
Note that the above $\bar{\mathbf{J}}^{(k)}$ is the approximation of the FOM Jacobian, i.e. before projection, and differs from $\tilde{\mathbf{J}}^{(k)}$ which would be the hyperreduced reduced-order Jacobian at iteration $k$. This is required for the calculation of the test basis in the LSPG framework. The test basis at iteration $k$ would be:
\begin{equation}
    \mathbf{W}^{(k)} \approx \tilde{\mathbf{W}}^{(k)} =  \bar{\mathbf{J}}^{(k)}\mathbf{V} = (\sum_{e \in \tilde{\mathcal{E}}} \xi_e \mathbf{L}_e^{T} \mathbf{J}_e^{(k)} \mathbf{L}_{e^+})\mathbf{V}.\label{hyp_test_basis}
 \end{equation}
The right-hand-side of equation \ref{eq:lspg_newton} is the reduced-order residual which can also be found using the weights and the test basis found above:
\begin{equation}
    \mathbf{\hat{R}}^{(k)} \approx \tilde{\mathbf{R}}^{(k)} =  \sum_{e \in \tilde{\mathcal{E}}} \xi_e  \tilde{\mathbf{W}}^{(k)^{T}} \mathbf{L}_e^{T} \mathbf{R}_e^{(k)}, \label{eq:on_res}
 \end{equation}
where $\mathbf{R}_e^{(k)}$ is the contribution of element $e$  to the global, FOM residual for the solution $\tilde{\mathbf{w}}^{(k)}$. Thus the dependence on the high-dimensional model size $N$ has been removed from the LSPG projection approach and the cost should scale independently of this size. Below in Algorithm 2, the online computation of the reduced-order solution $\tilde{\mathbf{w}}_h$ is summarized.

\begin{algorithm}[H]
\caption{Online Computation of the reduced-order solution by the LSPG projection with Hyperreduced Quantities}\label{alg:online}
\begin{algorithmic}
\STATE{\textbf{Inputs:}}
\STATE{POD Basis $\mathbf{V}$, Reference state $\mathbf{w_{ref}}$, Convergence tolerance $\epsilon$}
\STATE{Reduced Mesh Element Set $\widetilde{\mathcal{E}}$}
\STATE{Associate Set of Weights $\boldsymbol{\xi}_{\mathcal{E}} = \{\xi_{e} | e \in \widetilde{\mathcal{E}}\}$}
\STATE{\textbf{Output:}}
\STATE{Converged approximate solution $\tilde{\mathbf{w}}_h$}
\STATE{-------------------------------------------------}
\STATE{Ensure consistency of initial guess: $    \tilde{\textbf{w}}_{\text{proj}}^{(0)} = \textbf{w}_{\text{ref}} + \textbf{V}(\textbf{V}^T(\textbf{w}^{(0)} - \textbf{w}_{\text{ref}}))$}
\STATE{Evaluate the hyperreduced Jacobian $\bar{\mathbf{J}}^{(0)}$ using Equation \ref{hyp_jac} and use this to find the test basis $\tilde{\mathbf{W}}^{(0)}$ using Equation \ref{hyp_test_basis}}
\STATE{Evaluate the hyperreduced residual $\tilde{\mathbf{R}}^{(0)}$ using Equation \ref{eq:on_res} and compute the $L_2$ norm of the initial hyperreduced residual $r = \left\Vert \tilde{\mathbf{R}}^{(0)} \right\Vert_2$}
\WHILE{$r > \epsilon$}
    \STATE{Evaluate the hyperreduced Jacobian $\bar{\mathbf{J}}^{(k)}$ using Equation \ref{hyp_jac} and use this to find the test basis $\tilde{\mathbf{W}}^{(k)}$ using Equation \ref{hyp_test_basis}}
    \STATE{Evaluate the hyperreduced residual $\tilde{\mathbf{R}}^{(k)}$ using Equation \ref{eq:on_res}}
    \STATE{Solve $\left[\tilde{\mathbf{W}}^{(k)^{T}}\tilde{\mathbf{W}}^{(k)}\right]\textbf{p}^{(k)} = -\tilde{\mathbf{R}}^{(k)}$ }
    \STATE{Compute the step length $a^{(k)}$ by a line-search procedure}
    \STATE{Update the approximate solution $\tilde{\textbf{w}}^{(k+1)} = \tilde{\textbf{w}}^{(k)} + \textbf{V}(a^{(k)}\textbf{p}^{(k)})$}
    \STATE{Update $\tilde{\mathbf{R}}^{(k)}$ and $r = \left\Vert \tilde{\mathbf{R}}^{(0)} \right\Vert_2$} 
\ENDWHILE
\end{algorithmic}
\end{algorithm}

%% file: 05_adaptSamp.tex
\section{Adaptive Sampling Procedure}\label{sec:adapt}

This section outlines the modifications required to integrate hyperreduction into the adaptive sampling procedure developed in \cite{blaisTHESIS, donJournal}. First, the approximation of the residual and Jacobian in the Newton iterations used to compute the ROM solutions, must be updated. Additionally, the dual-weighted error indicators which approximate the error at ROM points require adjustments to account for the additional approximation introduced by the hyperreduction. Finally, a method for calculating the work units associated with the various sampling procedures will be presented.

\subsection{Overview of the Adaptive Sampling Procedure}

In this work, we focus on parameterized problems. For this class of problems, the ROM is built by sampling the FOM at various parameter combinations $\boldsymbol{\mu}$ in the parameter domain $\mathcal{D}$. It is then used to predict the solution at new parameter locations where the ``true" solution is unknown. Without a reliable quantification of the error in the ROM, little value can be placed in these predictions \cite{blaisTHESIS, donJournal}. In CFD problems, the focus is not only on the accuracy of the solution approximation but in particular on the error introduced in the functional $\mathcal{J}$ by the ROM. For example, in aerodynamic shape optimization, the iterative design of the aircraft or airfoil is guided by an output of interest, such as lift or drag, therefore confidence in the functional value is additionally important.

Building a ROM for parameterized problems involves two key challenges: identifying an \textit{a posteriori} functional error estimate and determining optimal snapshot placement. The accuracy of the ROM depends heavily on the quality of the training data used to construct the reduced-order basis. Sampling too few points or placing them poorly within the parameter domain can lead to an inaccurate model. Conversely, oversampling the FOM incurs unnecessary computational costs. Therefore, it is essential to develop a method for selecting parameter values judiciously to collect snapshots efficiently.

In \cite{blaisTHESIS, donJournal}, a novel adaptive sampling procedure was developed with two major objectives. The first is to minimize the number of FOM solutions computed through efficient snapshot selection. The second is, through the use of an \textit{a posteriori} error estimate based on dual-weighted residual error indicators, ensure that a prescribed
output error tolerance is estimated to be satisfied across the entire parameter domain $\mathcal{D}$.

The following is a brief summary of the steps taken in the adaptive sampling procedure:
\begin{itemize}
    \item  An initial set of snapshots is evenly distributed across the parameter space, along with a set of ROM points placed between these snapshots. ROM points are used to probe the error between the ROM and FOM. Two dual-weighted residual errors are tracked at each ROM point.
    \item At each adaptive cycle, a radial basis function (RBF) interpolation is used to model the error across the parameter space. The next snapshot is placed at the extremum of the RBF. The error estimates at the ROM points are updated, and new ROM points are added near the newly selected snapshot.
    \item This process is repeated until the specified error tolerance is met at all ROM points across the parameter domain.
\end{itemize}
Further details on these steps, along with an algorithm describing the procedure, can be found in \cite{blaisTHESIS, donJournal}.

\subsection{Dual-Weighted Residual Errors}

Two distinct dual-weighted residual (DWR) errors are used to approximate the error at a ROM point: $\epsilon_f$, which measures the error between the FOM and the ROM, and $\epsilon_r$, which quantifies the error between a ``coarse" ROM and a ``fine" ROM. The coarseness of a ROM refers to the dimension of the basis and which iteration of the sampling cycle a ROM is from. Meaning, that a ``coarse" ROM will have a coarse basis $V_H$ which will have fewer columns than the fine basis $V_h$. We then say that for a given parameter location $\boldsymbol{\mu}$ with full-order solution $\mathbf{w}$, there are two ROM solutions, one computed on the fine basis $\mathbf{w}_{\text{ref}_h} + \mathbf{V}_h\hat{\mathbf{w}}_{h}$ and one computed on the coarse basis $\mathbf{w}_{\text{ref}_H} + \mathbf{V}_H\hat{\mathbf{w}}_{H}$. Therefore, the total error $\epsilon$ in the functional at this ROM is approximately \cite{blaisTHESIS}:
\begin{equation}
\begin{split}
\epsilon = \mathcal{J}(\mathbf{w}) - \mathcal{J}(\mathbf{w}_{\text{ref}_h} + \mathbf{V}_h\hat{\mathbf{w}}_{h})  &= \left(\mathcal{J}(\mathbf{w}) - \mathcal{J}(\mathbf{w}_{\text{ref}_H} + \mathbf{V}_H\hat{\mathbf{w}}_{H}) \right) \\
& \;\;\; + \left( \mathcal{J}(\mathbf{w}_{\text{ref}_H} + \mathbf{V}_H\hat{\mathbf{w}}_{H}) - \mathcal{J}(\mathbf{w}_{\text{ref}_h} + \mathbf{V}_h\hat{\mathbf{w}}_{h}) \right)
\\
& \approx \epsilon_f + \epsilon_r.
\end{split}
\end{equation}
$\epsilon_r$ will be the DWR error impacted by hyperreduction as it requires the evaluation of the reduced residual and Jacobian, therefore we will re-derive it in the context of a HROM.

\subsubsection{Error Estimate Between Coarse and Fine Hyperreduced Reduced-Order Models}
The error estimate between a coarse and fine ROM derived in \cite{blaisTHESIS, donJournal} will be extended here to account for a hyperreduced ROM. As before, let the coarse ROM solution be denoted by $\mathbf{w}_{\text{ref}_H} + \mathbf{V}_H\hat{\mathbf{w}}_{H}$, where the subscript $H$ indicates a coarse-space representation. Similarly, let the fine ROM solution be $\mathbf{w}_{\text{ref}_h} + \mathbf{V}_h\hat{\mathbf{w}}_{h}$, where the subscript $h$ represents a fine-space representation. For brevity, these solutions will also be referred to as $\tilde{\mathbf{w}}_H$ and $\tilde{\mathbf{w}}_h$, respectively. Then the solutions must satisfy the following:
\begin{equation}
\begin{split}
\tilde{\mathbf{R}}_h(\tilde{\mathbf{w}}_h) &= 0,
\\
\tilde{\mathbf{R}}_H(\tilde{\mathbf{w}}_H) &= 0,
\end{split}
\end{equation}
where $\tilde{\mathbf{R}}_*$ represents the hyperreduced reduced-order residual from each HROM. Using a first-order Taylor series expansion of the fine-dimension hyperreduced reduced-order residual about the coarse reduced-order solution, we find:
\begin{equation}
\begin{split}
\tilde{\mathbf{R}}_h(\tilde{\mathbf{w}}_h) &= 0  \approx \tilde{\mathbf{R}}_h(\tilde{\mathbf{w}}_H) + \left[ \frac{\partial \tilde{\mathbf{R}}_h}{\partial \hat{\mathbf{w}}_h} \bigg|_{\hat{\mathbf{w}}_H} \right] (\tilde{\mathbf{w}}_h - \tilde{\mathbf{w}}_H). \label{eq:res_exp}
\end{split}
\end{equation}

Note that the derivative is taken with respect to the reduced-order solution representation $\hat{\mathbf{w}}_h$ rather than $\tilde{\mathbf{w}}h$. This choice simplifies the analysis and is justified because the errors arising from the derivative of the residual with respect to the reference state $\mathbf{w}_{\text{ref}_h}$ or the POD basis $\mathbf{V}_h$ are primarily influenced by the accuracy of the DG approach not the ROM.

Conducting the chain rule on this derivative, we find:
\begin{equation}
\begin{split}
\frac{\partial \tilde{\mathbf{R}}_h}{\partial \hat{\mathbf{w}}_h} \bigg|_{\hat{\mathbf{w}}_H} & = \frac{\partial \tilde{\mathbf{R}}_h}{\partial \tilde{\mathbf{w}}_h} \bigg|_{\tilde{\mathbf{w}}_H}  \frac{\partial \tilde{\mathbf{w}}_h}{\partial \hat{\mathbf{w}}_h} \\
& = \frac{\partial \tilde{\mathbf{R}}_h}{\partial \tilde{\mathbf{w}}_h} \bigg|_{\tilde{\mathbf{w}}_H}  \mathbf{V}_h .\label{eq:der_hat}
\end{split}
\end{equation}

Using the definition of the hyperreduced residual in Equation \ref{eq:on_res}, we can then write the remaining derivative as follows:
\begin{equation}
\begin{split}
\frac{\partial \tilde{\mathbf{R}}_h}{\partial \tilde{\mathbf{w}}_h} & = \frac{\partial}{\partial \tilde{\mathbf{w}}_h}
\left( \sum_{e \in \tilde{\mathcal{E}}} \xi_e  \tilde{\mathbf{W}}_h^{T} \mathbf{L}_e^{T} \mathbf{R}_e \right) \\
& = 
 \sum_{e \in \tilde{\mathcal{E}}} \xi_e  \frac{\partial \tilde{\mathbf{W}}_h^{T}}{\partial \tilde{\mathbf{w}}_h} \mathbf{L}_e^{T} \mathbf{R}_e  + \sum_{e \in \tilde{\mathcal{E}}} \xi_e  \tilde{\mathbf{W}}_h^{T}{\partial \tilde{\mathbf{w}}_h} \mathbf{L}_e^{T} \frac{\partial \mathbf{R}_e} {\partial \tilde{\mathbf{w}}_h} .\label{eq:der_hyp_res}
\end{split}
\end{equation}
Given that a LSPG framework was used, second-order sensitivities will arise when taking the derivative of the test basis with respect to the solution, as it contains the Jacobian as well. \cite{zahr2015} notes that this can be safely neglected with a minor loss in accuracy. Therefore, the first term on the right-hand side of Equation \ref{eq:der_hyp_res} will be ignored.

In a similar fashion, the first-order Taylor series expansions of the functional is:
\begin{equation}
\begin{split}
\mathcal{J}(\tilde{\mathbf{w}}_h) \approx \mathcal{J}(\tilde{\mathbf{w}}_H) + \left[ \frac{\partial \mathcal{J}}{\partial \hat{\mathbf{w}}_h} \bigg|_{\hat{\mathbf{w}}_H} \right] (\tilde{\mathbf{w}}_h - \tilde{\mathbf{w}}_H).
\end{split}
\end{equation}



We then define $\tilde{\psi}_h $ as the solution of the hyperreduced reduced-order dual problem:
\begin{equation}
\label{eq:hyp_DWR_pre}
\left[  \frac{\partial \tilde{\mathbf{R}}_h}{\partial \hat{\mathbf{w}}_h} \bigg|_{\hat{\mathbf{w}}_H} \right]^T  \tilde{\psi}_h = -\left[ \frac{\partial \mathcal{J}}{\partial \hat{\mathbf{w}}_h} \bigg|_{\hat{\mathbf{w}}_H} \right]^T .
\end{equation}

Using Equations \ref{eq:der_hat} and \ref{eq:der_hyp_res}, and applying the chain rule in a similar manner to the derivative of the functional, we find:
\begin{equation}
\label{eq:hyp_DWR}
\left[ \left\{ \sum_{e\in\widetilde{\mathcal{E}}}\xi_e\mathbf{W}_h^T\mathbf{L}_e^T \frac{\partial \mathbf{R}_e}{\partial \tilde{\mathbf{w}}_h} \bigg|_{\tilde{\mathbf{w}}_H}  \right\}\mathbf{V}_h \right]^T  \tilde{\psi}_h = -\left[ \frac{\partial \mathcal{J}}{\partial \tilde{\mathbf{w}}_h} \bigg|_{\tilde{\mathbf{w}}_H} \mathbf{V}_h\right]^T .
\end{equation}
Note that the set of weights and the test basis will be from the fine ROM, meaning they are computed in the most recent iteration. Using the hyperreduced reduced-order adjoint in Equation \ref{eq:hyp_DWR} and the Taylor series expansion in Equation \ref{eq:res_exp}, an error metric between the fine and coarse space can be
defined:
\begin{equation}\label{eq:hyp_DWR_final}
\begin{split}
\epsilon_r &= \mathcal{J}(\tilde{\mathbf{w}}_H) - \mathcal{J}(\tilde{\mathbf{w}}_h) \\
&\approx - \left[ \frac{\partial \mathcal{J}}{\partial \tilde{\mathbf{w}}_h} \bigg|_{\tilde{\mathbf{w}}_H} \right](\tilde{\mathbf{w}}_h - \tilde{\mathbf{w}}_H)
\\
& \approx - \left[ \frac{\partial \mathcal{J}}{\partial \tilde{\mathbf{w}}_h} \bigg|_{\tilde{\mathbf{w}}_H} \right]\left[ \frac{\partial \tilde{\mathbf{R}}_h}{\partial \tilde{\mathbf{w}}_h} \bigg|_{\tilde{\mathbf{w}}_H} \right]^{-1}(\tilde{\mathbf{R}}_h(\tilde{\mathbf{w}}_h) - \tilde{\mathbf{R}}_h(\tilde{\mathbf{w}}_H))
\\
& \approx - \tilde{\psi}_h^T \tilde{\mathbf{R}}_h(\tilde{\mathbf{w}}_H) .
\end{split}
\end{equation}

\subsection{Updated Goal-Oriented Adaptive Sampling}

The adaptive sampling procedure summarized in Section 5.5 and Algorithm 6 of \cite{blaisTHESIS} is updated here with modifications to accommodate the current work. While not all components are discussed in detail, we refer the reader to \cite{blaisTHESIS} for comprehensive descriptions. The major changes come from the need to train the ECSW weights at every iteration before solving for the ROM probing points, as well as the introduction of the updated dual-weighted residual (DWR) error indicator for measuring the error between coarse and fine ROMs. For clarity, algorithms from \cite{blaisTHESIS} are referenced as D-* (where * corresponds to the original numbering), while algorithms introduced in this work are referred to by their sequence number in this document.

\begin{algorithm}[H]
\caption{Goal-Oriented Adaptive Sampling}\label{alg:adapt_samp}
\begin{algorithmic}
\STATE{\textbf{Inputs:}}
\STATE{Error tolerance $\epsilon$, NNLS tolerance $\epsilon_{\text{NNLS}}$}
\STATE{Initial snapshot parameter locations $\mathcal{S}^{\boldsymbol{\mu}} = \{ \boldsymbol{\mu}_1^s, \boldsymbol{\mu}_2^s, \dots, \boldsymbol{\mu}_k^s \}$}
\STATE{\textbf{Output:}}
\STATE{POD basis $\mathbf{V}$, Maximum error $\epsilon_{\text{max}}$}
\STATE{-------------------------------------------------}
\STATE{Compute initial set of $k$ snapshot solutions $\mathcal{S}^{\mathbf{w}}$ at $\mathcal{S}^{\boldsymbol{\mu}}$}
\STATE{Compute initial POD basis $\mathbf{V}_H$ from $\mathcal{S}^{\mathbf{w}}$ using Algo. D-2}
\STATE{Determine initial set of ROM probing parameter locations $\mathcal{P}^{\boldsymbol{\mu}} = \{ \boldsymbol{\mu}_1^p, \boldsymbol{\mu}_2^p, \dots, \boldsymbol{\mu}_k^p \}$ using Algo. D-3}
\STATE{Solve for the reduced mesh set and weights using the ECSW approach in Algo. \ref{alg:cap} or the Jacobian-based equivalent}
\FOR{$\boldsymbol{\mu}^p \in \mathcal{P}^{\boldsymbol{\mu}}$}
    \STATE{Compute the HROM solution $\tilde{\mathbf{w}}(\boldsymbol{\mu}^p)$ with POD basis $\mathbf{V}_H$ using Algo. \ref{alg:online}}
    \STATE{Compute error between FOM and coarse ROM, $\epsilon_f(\boldsymbol{\mu}^p) = - \psi^T\mathbf{R}(\tilde{\mathbf{w}})$}
    \STATE{Update total error, $\epsilon(\boldsymbol{\mu}^p) \leftarrow \epsilon_f(\boldsymbol{\mu}^p$)}
\ENDFOR
\STATE{Compute RBF and obtain maximum error $\epsilon_{\text{max}}$ at $\boldsymbol{\mu}_{\text{max}}$ using Algo. D-4}
\WHILE{$\epsilon_{\text{max}} > \epsilon$}
    \STATE{Compute new snapshot solution $\mathbf{w}(\boldsymbol{\mu}_{\text{max}})$, update $\mathcal{S}^{\mathbf{w}} = \mathcal{S}^{\mathbf{w}} \cup \{ \mathbf{w}(\boldsymbol{\mu}_{\text{max}})\}$}
    \STATE{Compute $\mathbf{V}_h$ from $\mathcal{S}^{\mathbf{w}}$ using Algo. D-2}
    \STATE{Solve for the reduced mesh set and weights using the ECSW approach in Algo. \ref{alg:cap} or the Jacobian-based equivalent}
    \FOR{$\boldsymbol{\mu}^p \in \mathcal{P}^{\boldsymbol{\mu}}$}
        \STATE{Compute error between coarse ROM and fine ROM, $\epsilon_r(\boldsymbol{\mu}^p) = - \tilde{\psi}_h^T \tilde{\mathbf{R}}_h(\tilde{\mathbf{w}}_H)$}
        \STATE{Update total error, $\epsilon(\boldsymbol{\mu}^p) \leftarrow \epsilon_f(\boldsymbol{\mu}^p) + \epsilon_r(\boldsymbol{\mu}^p)$}
    \ENDFOR
    \STATE{Recompute solution at select ROM probing points in $\mathcal{P}^{\boldsymbol{\mu}}$ using Algo. D-5}
    \STATE{Add new ROM probing points $\mathcal{P}_*^{\boldsymbol{\mu}}$ at midpoint of $\boldsymbol{\mu}_{\text{max}}$ and $n_p + 1$ nearest neighbours in $\mathcal{S}^{\boldsymbol{\mu}}$}
    \FOR{$\boldsymbol{\mu}_*^p \in \mathcal{P}_*^{\boldsymbol{\mu}}$}
        \STATE{Update $\mathcal{P}^{\boldsymbol{\mu}} = \mathcal{P}^{\boldsymbol{\mu}} \cup \{ \boldsymbol{\mu}_*^p \}$}
        \STATE{Let $\mathbf{V}_H = \mathbf{V}_h$}
        \STATE{Compute ROM solution $\tilde{\textbf{w}} (\boldsymbol{\mu}_*^p)$ with POD basis $\mathbf{V}_H$ using Algo. \ref{alg:online}}
        \STATE{Compute error between FOM and coarse ROM, $\epsilon_f(\boldsymbol{\mu}^p) = - \psi^T\mathbf{R}(\tilde{\mathbf{w}})$}
        \STATE{Update total error, $\epsilon(\boldsymbol{\mu}_*^p) \leftarrow \epsilon_f(\boldsymbol{\mu}_*^p$)}
    \ENDFOR
    \STATE{Compute RBF and obtain maximum error $\epsilon_{\text{max}}$ at $\boldsymbol{\mu}_{\text{max}}$ using Algo. D-4}
\ENDWHILE
\end{algorithmic}
\end{algorithm}

\subsection{Analysis of Computational Cost Savings through Work Units}\label{sec:work}

Currently, due to the implementation of the DG method in \texttt{PHiLiP} \cite{philip_2022}, it is not possible to compute the residual or Jacobian on a specific mesh element. TAs a result, a pseudo-implementation of the ECSW method is employed. wIn this approach, the residual and Jacobian are evaluated over the entire mesh, but only the elements in the reduced mesh are utilized in the HROMs. This means the CPU time cannot be used to compare the efficiency of the models. Instead, an approximation of work units, based on matrix dimensions and sparsity, is used to compare the cost of the different sampling procedures.

The number of work units at each cycle comes from two major contributions; the first is from the cost of solving for the reduced-order solutions at the new ROM points and the second is from re-evaluating the DWR error indicator between the coarse and fine ROM at all the previously added ROM points. The cost of solving for the FOM snapshots is not considered as it is not impacted by the hyperreduction. The cost of computing the FOM snapshots is excluded from these considerations as it is not impacted by the hyperreduction, and thus the cost incurred will be the same for any sampling procedure.

For ROM point solutions, the number of total non-linear iterations in a sampling iteration is multiplied by the cost of solving the linear system in Equation \ref{eq:lspg_newton}. Work units are estimated under the following assumptions: evaluating an entry in the residual or Jacobian costs one work unit; each floating-point operation (FLOP) in a matrix multiplication counts as one work unit; and the cost of solving the linear system is approximated using the worst-case operation count for the generalized minimal residual method (GMRES), which is $\mathcal{O}(n^3)$.

In the original adaptive sampling procedure which builds the ROM with no hyperreduction, the total number of work units to solve one non-linear iteration is given by:
\begin{equation}\label{eq:ROM_w_non_lin}
\begin{split}
    W_{\text{non lin}} = & (\text{units to evaluate the residual}) + \\
    & (\text{units to evaluate the Jacobian}) + \\
    & (\text{units to assemble LHS matrix and RHS vector in Equation \ref{eq:lspg_newton}}) + \\
    & (\text{units to solve for solution update}) \\ 
    =  & (N) + (N^2) + [(2Nn_i + n_i^2 + N + n_i)(2N-1)] + (n_i^3).
\end{split}
\end{equation}
where $N$ is the dimension of the FOM and $n_i$ is the number of bases in the POD or reduced-order subspace dimension at adaptive sampling cycle $i$. In the sampling procedures where the hyperreduction is used to evaluate the residual and Jacobian to solve Equation \ref{eq:lspg_newton} at the ROM points, each nonlinear iteration has work units determined by:
\begin{equation}\label{eq:HROM_w_non_lin}
\begin{split}
    W_{\text{non lin}} = & (\text{units to evaluate the hyperreduced residual}) + \\
    & (\text{units to evaluate the hyperreduced test basis}) + \\
    & (\text{units to assemble LHS matrix and RHS vector in Equation \ref{eq:lspg_newton}}) + \\
    & (\text{units to solve for solution update}) \\ 
    = & [n_{e_{i}}  (d_e + 2n_i d_e + n_i)] + [2n_{e_{i}} d_e d_e^+  +
        2 d_e d_e^+ n_{e_{i}} n_i] +  [n_i^2 (2N -1)] + (n_i^3).
\end{split}
 \end{equation}
where $n_{e_{i}}$ the number of elements in the reduced mesh at adaptive sampling cycle $i$. Details of how both Equations \ref{eq:ROM_w_non_lin} and \ref{eq:HROM_w_non_lin} are derived are included in \ref{sec:appendix}.

The other contribution is from the evaluation of the second DWR error indicator $\epsilon_r$ at all of the ROM points from the previous sampling cycles. Similar to the previous case, the work units for this evaluation can be broken down into a couple of components. For the original adaptive sampling procedure which builds the ROM with no hyperreduction, the work units required to solve for the error indicator are:
\begin{equation}\label{eq:ROM_w_eps}
\begin{split}
    W_{\text{DWR error}} = & (\text{units to evaluate the residual}) + \\
    & (\text{units to evaluate the Jacobian}) + \\
    & (\text{units to evaluate the derivative of the functional w.r.t. $\tilde{\mathbf{w}}$}) + \\
    & (\text{units to assemble the adjoint problem}) + \\
    & (\text{units to solve for the adjoint}) + \\ 
    & (\text{units to solve $\epsilon_r$}) \\
    =  & (N) + (N^2) + (N) + [(2Nn_i + n_i^2 + n_i)(2N-1)] + \\
    & (n_i^3) + [(N + n_i)(2N -1) + (2n_i -1)].
\end{split}
\end{equation}
In the sampling procedures where the hyperreduction is applied in the DWR error indicator, such as in Equation \ref{eq:hyp_DWR} and \ref{eq:hyp_DWR_final}, the work unit required to solve for the error indicator are:
\begin{equation}\label{eq:HROM_w_eps}
\begin{split}
    W_{\text{DWR error}} = & (\text{units to evaluate the hyperreduced residual}) + \\
    & (\text{units to evaluate the hyperreduced test basis}) + \\
    & (\text{units to evaluate the derivative of the functional w.r.t. $\tilde{\mathbf{w}}$}) + \\
    & (\text{units to assemble thee adjoint problem}) + \\
    & (\text{units to solve for the adjoint}) + \\ 
    & (\text{units to solve $\epsilon_r$}) \\
    =  & [n_{e_{i}}  (d_e + 2n_i d_e + n_i)] + [n_{e_{i}} d_e d_e^+  +
        2 d_e d_e^+ n_{e_{i}} n_i] + (N) \\
        & +  [n_{e_{i}}(d_e n_i + n_i) + 
        (n_i^2)(2N-1) + n_i(2N-1)] + (n_i^3) + (2n_i -1).
\end{split}
\end{equation}
Details of how both Equations \ref{eq:ROM_w_eps} and \ref{eq:HROM_w_eps} are derived are included in \ref{sec:appendix}.

Once the sampling procedure is completed, we can gather the number of nonlinear iterations at each cycle. We will refer to the total number sampling of cycles as $c$. This, along with the knowledge of the dimension of the reduced-order subspace $n_i$ and the number of previously added ROM points at each iteration, can be used to approximate the work units for a given adaptive sampling procedure. The algorithms for computing the work units for the three approaches used in the results section are included below. Algorithm \ref{alg:ROM_w} corresponds to the original sampling procedure, which does not include hyperreduction. Algorithm \ref{alg:HROM_no_DWR_w}  represents the sampling procedure that incorporates hyperreduction in the ROM point solution evaluations but excludes the updated DWR error. Algorithm \ref{alg:HROM_w_DWR_w} represents the sampling procedure that includes both hyperreduction in the ROM point solution evaluations and the hyperreduced DWR error.

\begin{algorithm}[H]
\caption{Evaluating Work Units for ROM}\label{alg:ROM_w}
\begin{algorithmic}
\STATE{\textbf{Inputs:}}
\STATE{Number of sampling cycles $c$, Number of design parameters $n_p$}
\STATE{Total number of nonlinear iterations in each sampling cycle ($\text{nnl} = {\mathcal{N}_1,\mathcal{N}_2, \dots, \mathcal{N}_c }$)}
\STATE{\textbf{Outputs:}}
\STATE{Work Units for each sampling cycle  $W_{\text{tot}}(i)$ for $i = 1, 2, \dots, c$}
\STATE{-------------------------------------------------}

\FOR{i = 1, 2, \dots, c}
        \STATE{Evaluate $W_{\text{non lin}}$ using  Equation \ref{eq:ROM_w_non_lin}}
        \STATE{Solve for the work units from ROM points solutions : $W_{\text{all ROM}} = (\mathcal{N}_i) * W_{\text{non lin}}$}
        \STATE{Evaluate $W_{\text{DWR error}}$ using \ref{eq:ROM_w_eps}}
        \STATE{Solve for the work units needed to find the DWR error at all previous ROM points: $W_{\text{all DWR}} = (n_p + 1)(i-1) * W_{\text{DWR error}}$}
        \STATE{$W_{\text{tot}}(i) = W_{\text{all ROM}} + W_{\text{all DWR}}$}
\ENDFOR
\end{algorithmic}
\end{algorithm}

\begin{algorithm}[H]
\caption{Evaluating Work Units for HROM without Hyper-DWR}\label{alg:HROM_no_DWR_w}
\begin{algorithmic}
\STATE{\textbf{Inputs:}}
\STATE{Number of sampling cycles $c$, Number of design parameters $n_p$}
\STATE{Total number of nonlinear iterations in each sampling cycle ($\text{nnl} = {\mathcal{N}_1,\mathcal{N}_2, \dots, \mathcal{N}_c }$)}
\STATE{\textbf{Outputs:}}
\STATE{Work Units for each sampling cycle  $W_{\text{tot}}(i)$ for $i = 1, 2, \dots, c$}
\STATE{-------------------------------------------------}

\FOR{i = 1, 2, \dots, c}
        \STATE{Evaluate $W_{\text{non lin}}$ using  Equation \ref{eq:HROM_w_non_lin}}
        \STATE{Solve for the work units from ROM points solutions : $W_{\text{all ROM}} = (\mathcal{N}_i) * W_{\text{non lin}}$}
        \STATE{Evaluate $W_{\text{DWR error}}$ using \ref{eq:ROM_w_eps}}
        \STATE{Solve for the work units needed to find the DWR error at all previous ROM points: $W_{\text{all DWR}} = (n_p + 1)(i-1) * W_{\text{DWR error}}$}
        \STATE{$W_{\text{tot}}(i) = W_{\text{all ROM}} + W_{\text{all DWR}}$}
\ENDFOR
\end{algorithmic}
\end{algorithm}

\begin{algorithm}[H]
\caption{Evaluating Work Units for HROM with Hyper-DWR}\label{alg:HROM_w_DWR_w}
\begin{algorithmic}
\STATE{\textbf{Inputs:}}
\STATE{Number of sampling cycles $c$, Number of design parameters $n_p$}
\STATE{Total number of nonlinear iterations in each sampling cycle ($\text{nnl} = {\mathcal{N}_1,\mathcal{N}_2, \dots, \mathcal{N}_c }$)}
\STATE{\textbf{Outputs:}}
\STATE{Work Units for each sampling cycle  $W_{\text{tot}}(i)$ for $i = 1, 2, \dots, c$}
\STATE{-------------------------------------------------}

\FOR{i = 1, 2, \dots, c}
        \STATE{Evaluate $W_{\text{non lin}}$ using  Equation \ref{eq:HROM_w_non_lin}}
        \STATE{Solve for the work units from ROM points solutions : $W_{\text{all ROM}} = (\mathcal{N}_i) * W_{\text{non lin}}$}
        \STATE{Evaluate $W_{\text{DWR error}}$ using \ref{eq:HROM_w_eps}}
        \STATE{Solve for the work units needed to find the DWR error at all previous ROM points: $W_{\text{all DWR}} = (n_p + 1)(i-1) * W_{\text{DWR error}}$}
        \STATE{$W_{\text{tot}}(i) = W_{\text{all ROM}} + W_{\text{all DWR}}$}
\ENDFOR
\end{algorithmic}
\end{algorithm}

%% file: 06_results.tex
\section{Results}\label{sec:results}

This section presents results from the updated adaptive sampling framework. Verfication is conducted for the implementation of the ECSW hyperreduction seperate from the sampling procedure in \ref{sec:verf}. As noted in section \ref{sec:model_red}, \cite{TEZAUR2022} observed that in steady-state problems, residual values can often be very close to zero which cause issues when using it as training data for the ECSW hyperreduction approach. Furthermore, based on the findings in \ref{sec:verf}, Jacobian-based training data has been shown to produce smaller reduced mesh sizes, resulting in more computationally efficient approximations of the residual and Jacobian. For these reasons, all HROMs constructed in this section will exclusively use Jacobian-based training data with an NNLS tolerance of $1 \times 10^{-6}$. Various flow conditions and design parameter combinations for the NACA0012 airfoil will make up the different test cases. HROMs will be compared to ROMs built with the original sampling procedure presented in \cite{blaisTHESIS}. To evaluate the effectiveness of the updated DWR error indicator in capturing the additional approximation errors introduced by the ECSW hyperreduction approach, two types of HROMs are studied: those without the updated DWR error indicator (referred to as ``HROM w/o Hyper-DWR") and those with the updated indicator (referred to as ``HROM w/ Hyper-DWR").

\begin{figure}[H]
\centering
\includegraphics[ width=0.65\textwidth]{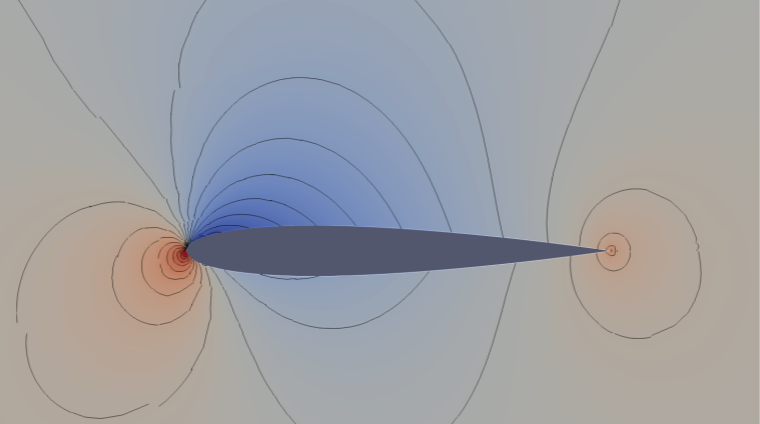}
\captionsetup{justification=centering}
\caption{Example Solution of NACA0012 Airfoil Test Case}
\label{fig:ex_NACA}
\end{figure}

All of the test cases will use the NACA0012 airfoil in inviscid flow. The Euler equations governing two-dimensional inviscid flow past an airfoil are:
\begin{equation}
       \nabla_x \mathbf{F}_x + \nabla_y \mathbf{F}_y = 0,
    \end{equation}
    where,
    \begin{equation}
\label{eq:euler}
\mathbf{F}_x = \begin{bmatrix}
\rho v_1        \\
\rho v_1 v_1 + p \\
\rho v_1 v_2 \\
v_1 (\rho e + p) \\
\end{bmatrix} \hspace{1cm} \text{and} \hspace{1cm} \mathbf{F}_y = \begin{bmatrix}
\rho v_2        \\
\rho v_1 v_2 \\
\rho v_2 v_2 + p \\
v_2 (\rho e + p) \\
\end{bmatrix},
\end{equation}
where $\rho$ is the density, $v_1$ and $v_2$ are velocity components, $p$ is the pressure and $e$ the internal energy. The far-field boundary conditions are specified by the Mach number, angle of attack $\alpha$, static pressure, and density. The airfoil surface is subject to no-slip boundary conditions and is adiabatic. Figure \ref{fig:ex_NACA} shows the pressure coefficient distribution of an example solution taken at an angle of attack of $2^{\circ}$.

\subsection{One Parameter NACA0012 Airfoil in Inviscid Subsonic Flow}\label{sec:one_param}

This test case will consider subsonic flow around the NACA0012 when the Mach number is fixed at 0.5, and the angle of attack is varied between $[0, 4]^{\circ}$. The full-order model is run on a mesh with 560 cells and a polynomial of order 0, resulting in 2240 DOFs. The adaptive sampling tolerance is set to $1E-4$ and the sampling process is initialized with 3 snapshots distributed evenly across the parameter space. As stated previously, three ROMs will be built and compared. The first will use the unchanged adaptive sampling procedure without hyperreduction incorporated, referred to in plots as the ``ROM". The second will be a hyperreduced ROM which uses the ECSW hyperreduction method to evaluate the residual and Jacobian to find reduced-order solutions but does include the update to the DWR error indicator discussed in section \ref{sec:adapt}. This HROM is referred to as the ``HROM w/o Hyper-DWR". The final HROM is built using the update adaptive sampling procedure shown in Algorithm \ref{alg:adapt_samp} which includes the new DWR error indicator which incorporates hyperreduction. This HROM is referred to as the ``HROM w/ Hyper-DWR". 

\begin{figure}[H]
\centering
\includegraphics[trim = 40 300 60 250, clip, width=0.7\textwidth]{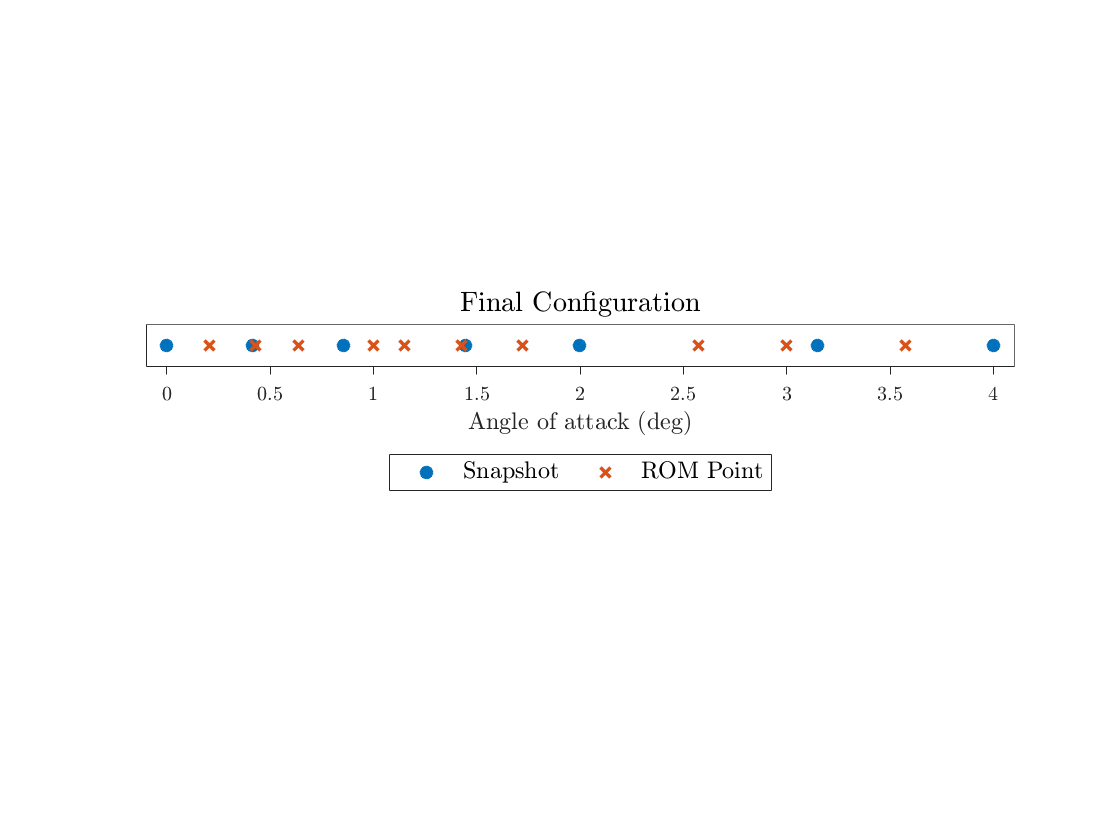}
\includegraphics[trim = 40 5 60 10, clip, width=0.7\textwidth]{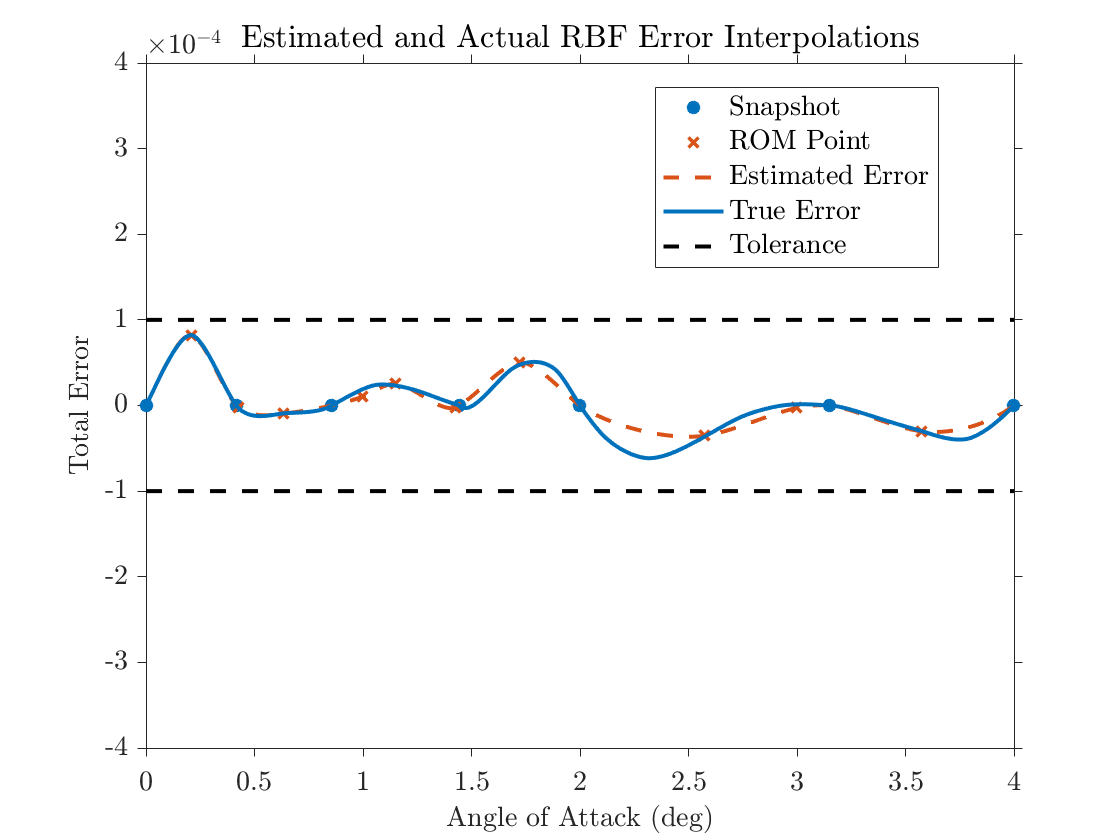}
\captionsetup{justification=centering}
\caption{Snapshot and ROM Points with Estimated and True Error for the ROM}\label{fig:one_subsonic_ROM}
\end{figure}

\begin{figure}[H]
\centering
\includegraphics[trim = 40 300 60 250, clip, width=0.7\textwidth]{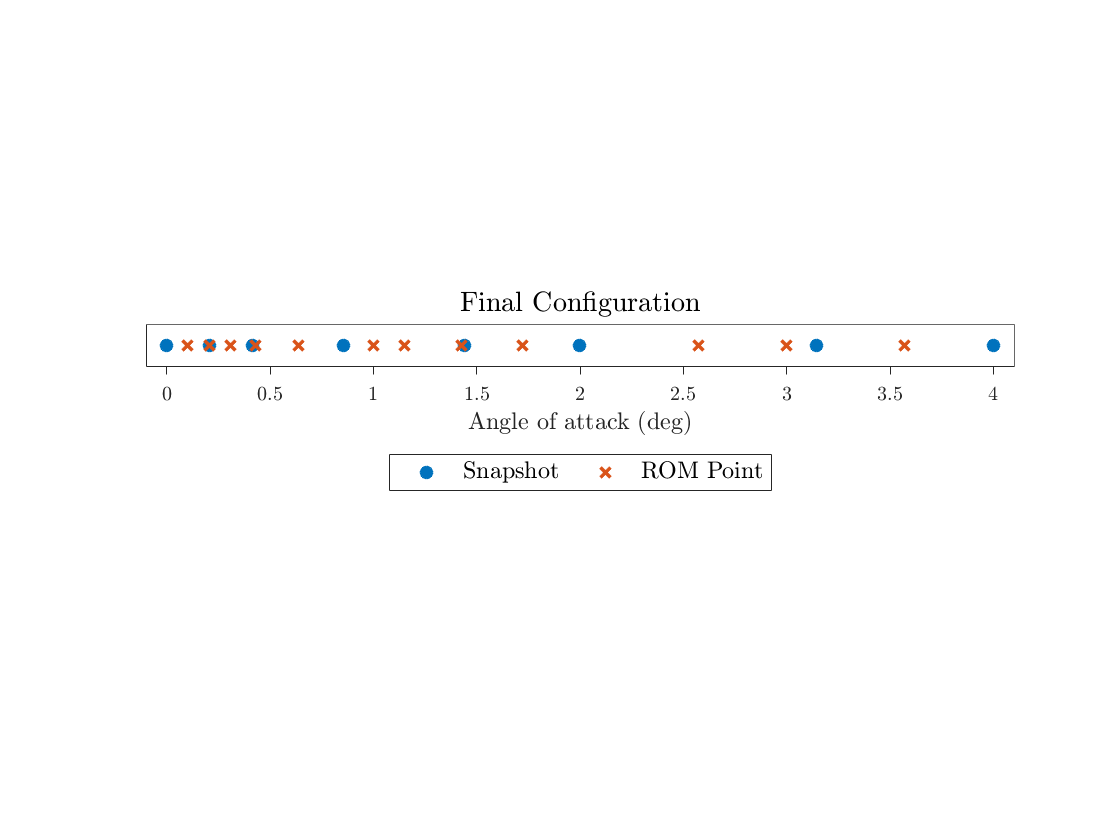}
\includegraphics[trim = 40 5 60 10, clip, width=0.7\textwidth]{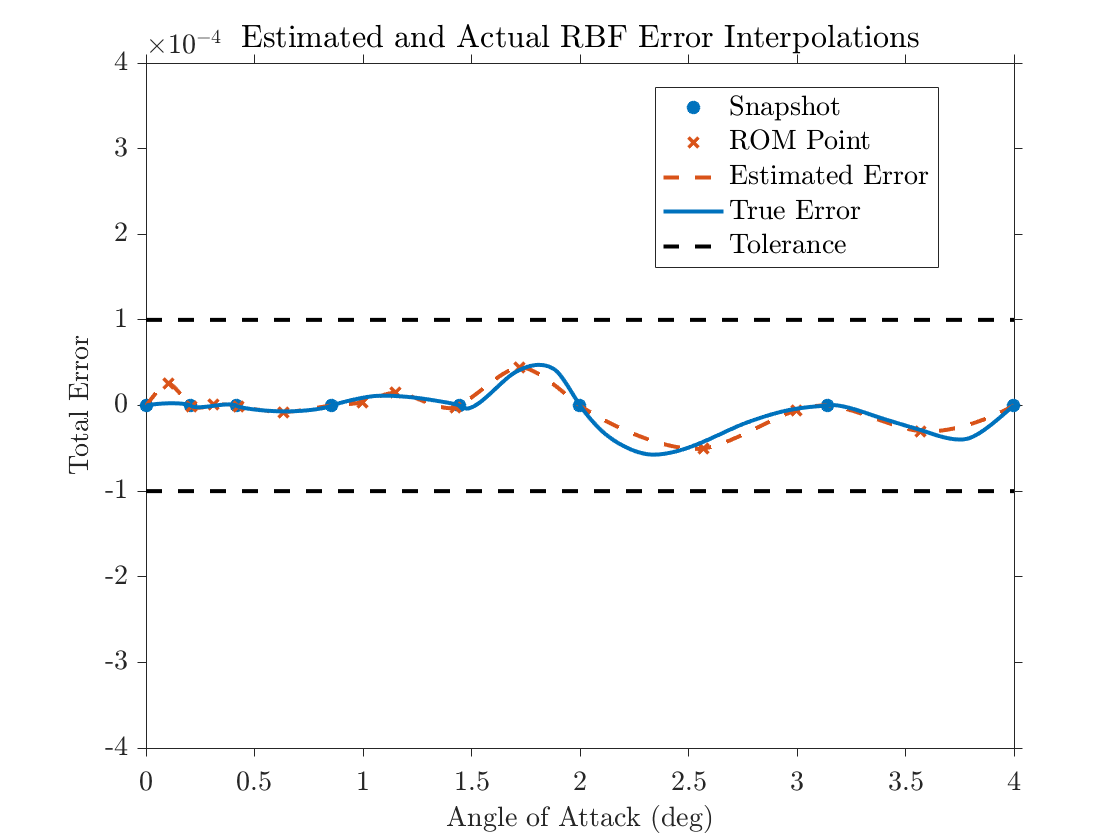}
\captionsetup{justification=centering}
\caption{Snapshot and ROM Points with Estimated and True Error for the HROM w/o Hyper-DWR}\label{fig:one_subsonic_HROM}
\end{figure}

\begin{figure}[H]
\centering
\includegraphics[trim = 40 300 60 250, clip, width=0.7\textwidth]{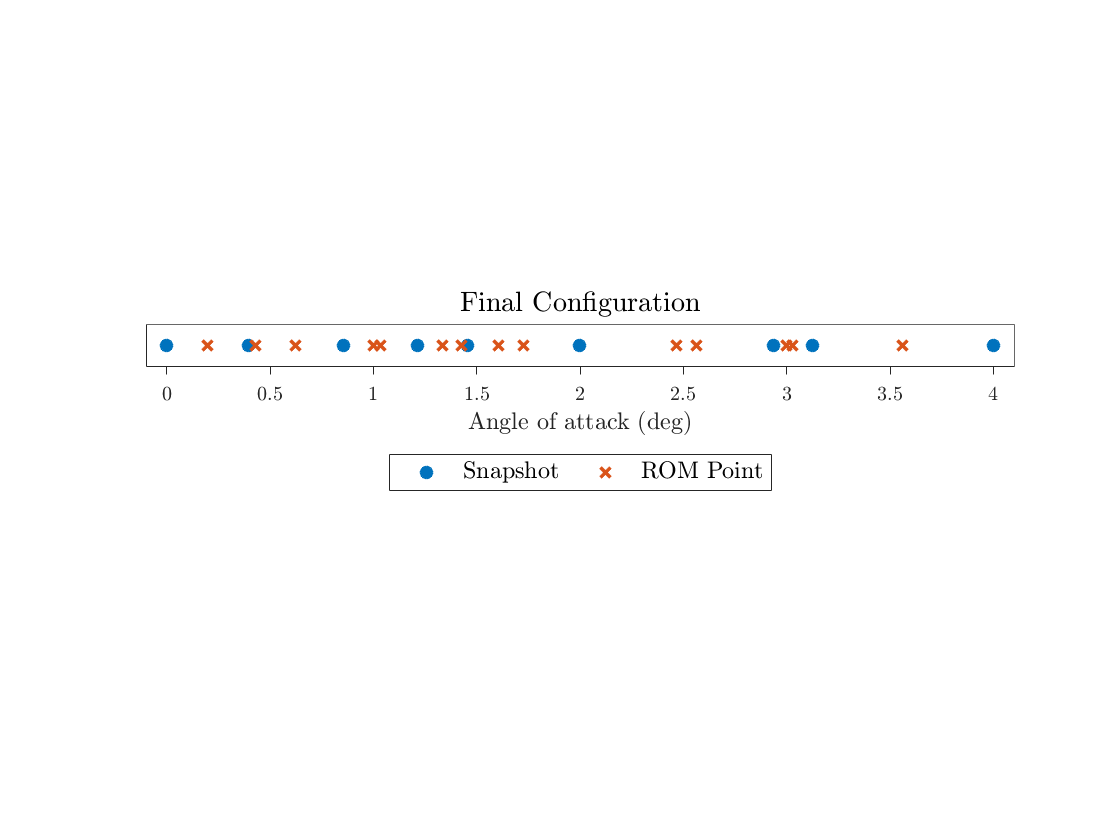}
\includegraphics[trim = 40 5 60 10, clip, width=0.7\textwidth]{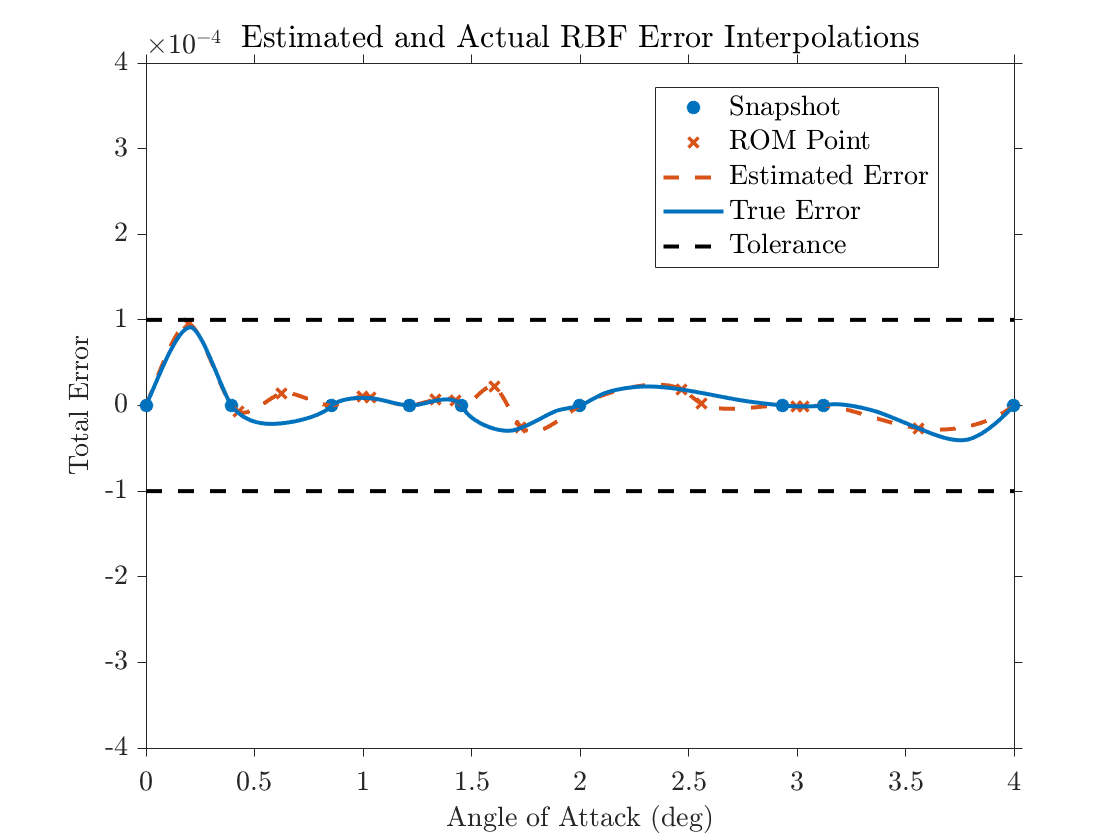}
\captionsetup{justification=centering}
\caption{Snapshot and ROM Points with Estimated and True Error for the HROM w/ Hyper-DWR}\label{fig:one_subsonic_HROM_w_DWR}
\end{figure}

Figures \ref{fig:one_subsonic_ROM} to \ref{fig:one_subsonic_HROM_w_DWR} show the results of the three models in the order they were just introduced. Each figure includes a plot of the final configuration of the parameter domain which shows the placement of the FOM snapshots used to build the reduced basis and the ROM points used to estimate the error distribution. They also include a plot of the ``estimated" and ``true" error distributions for each ROM. The estimated error is produced via an RBF interpolation using the snapshot locations and ROM points from the adaptive sampling cycle. The true error is found by sampling the error between each ROM and the FOM at 20 evenly distributed points across the parameter space at the end of the sampling procedure. The functional error tolerance bound is also plotted. It can be seen that all three processes successfully predicted the online behaviour of the ROMs, as the ``true" error remained inside of the tolerance bound for all of the plots.

From Figure \ref{fig:one_subsonic_HROM}, we can see that for the HROM without the hyperreduced error indicator the distribution of the snapshots is similar to the ROM, with one additional point added at $0.02^{\circ}$. This is likely because the error indicator is what dictates the placement of the snapshot locations, so without consideration of the hyperreduction in the DWR error indicator, the placement of the snapshots will remain similar unless the ROM point solutions are also inaccurate due to hyperreduction of the residual and Jacobian in the Newton iterations in Equation \ref{eq:lspg_newton}.

\begin{figure}[H]
\centering
\includegraphics[trim = 40 10 60 10, clip, width=0.7\textwidth]{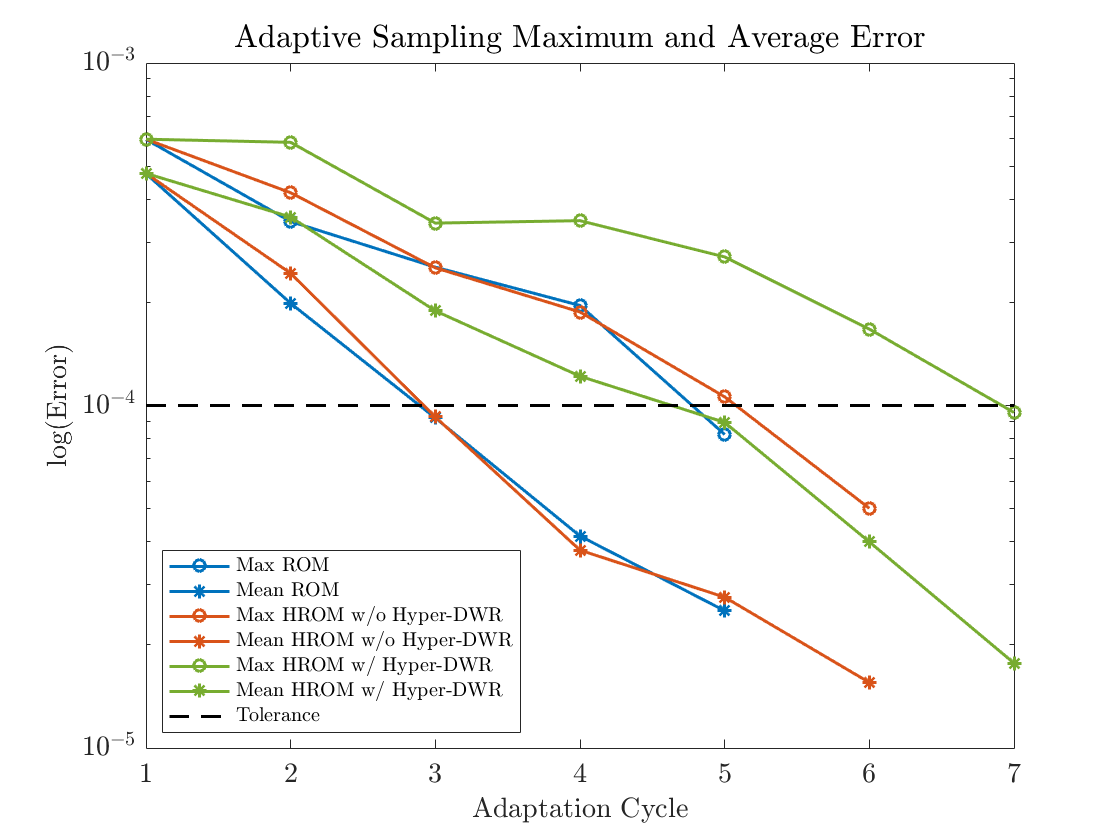}
\captionsetup{justification=centering}
\caption{Maximum and Mean Error Estimate at the ROM Points during each Sampling Iteration for the Three Models}\label{fig:one_param_iter}
\end{figure}

Figure \ref{fig:one_subsonic_HROM_w_DWR} shows that once hyperreduction is introduced into the error indicator, more snapshots are required to achieve the same error tolerance as in the ROM. In Figure \ref{fig:one_param_iter}, the maximum and average error at all of the ROM points for each model are plotted over the adaptive sampling iterations. With the introduction of hyperreduction and the additional layer of approximation, one additional adaptive sampling cycle is required for the HROM w/o hyper-DWR and two additional sampling cycles are required for the HROM w/ hyper-DWR to achieve the same tolerance that the ROM reached in 5 cycles.

\begin{figure}[H]
\centering
\includegraphics[trim = 300 0 0 70, clip, width=0.7\textwidth]{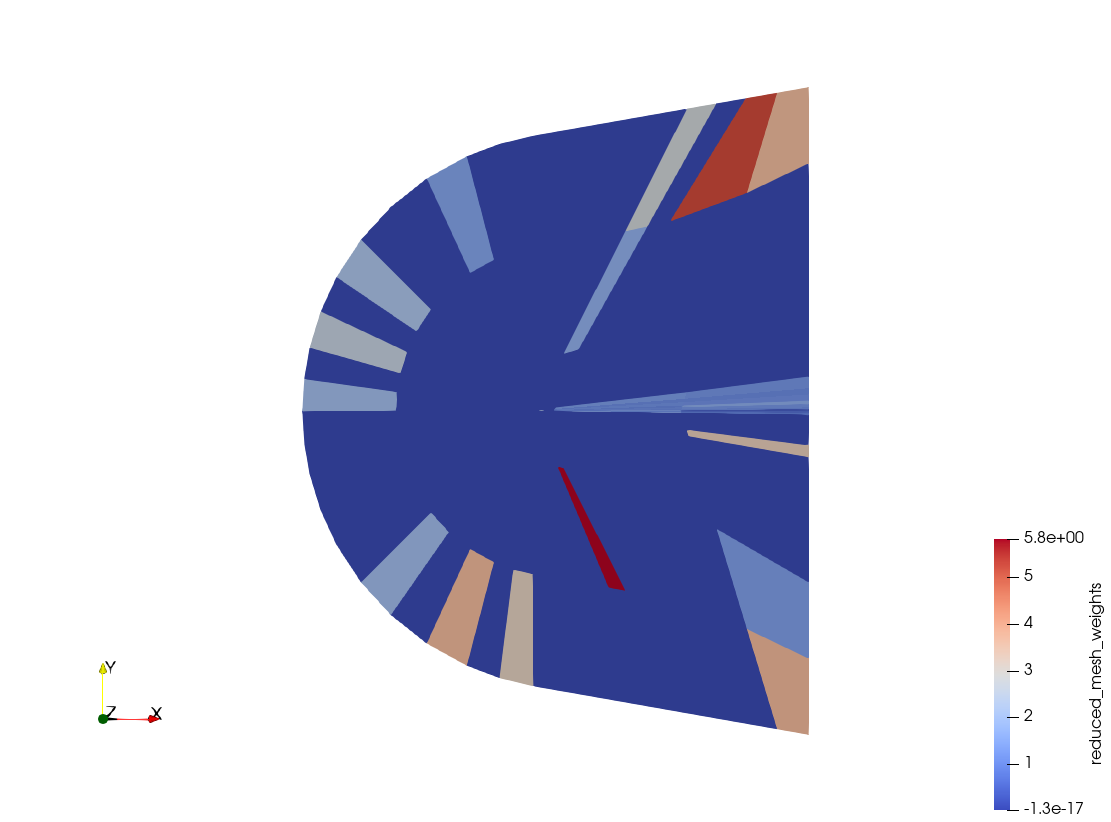}
\caption{ECSW Reduced Mesh Set for the HROM w/o Hyper-DWR}\label{fig:one_param_HROM_mesh}
\end{figure}

\begin{figure}[H]
\centering
\includegraphics[trim = 300 0 0 70, clip, width=0.7\textwidth]{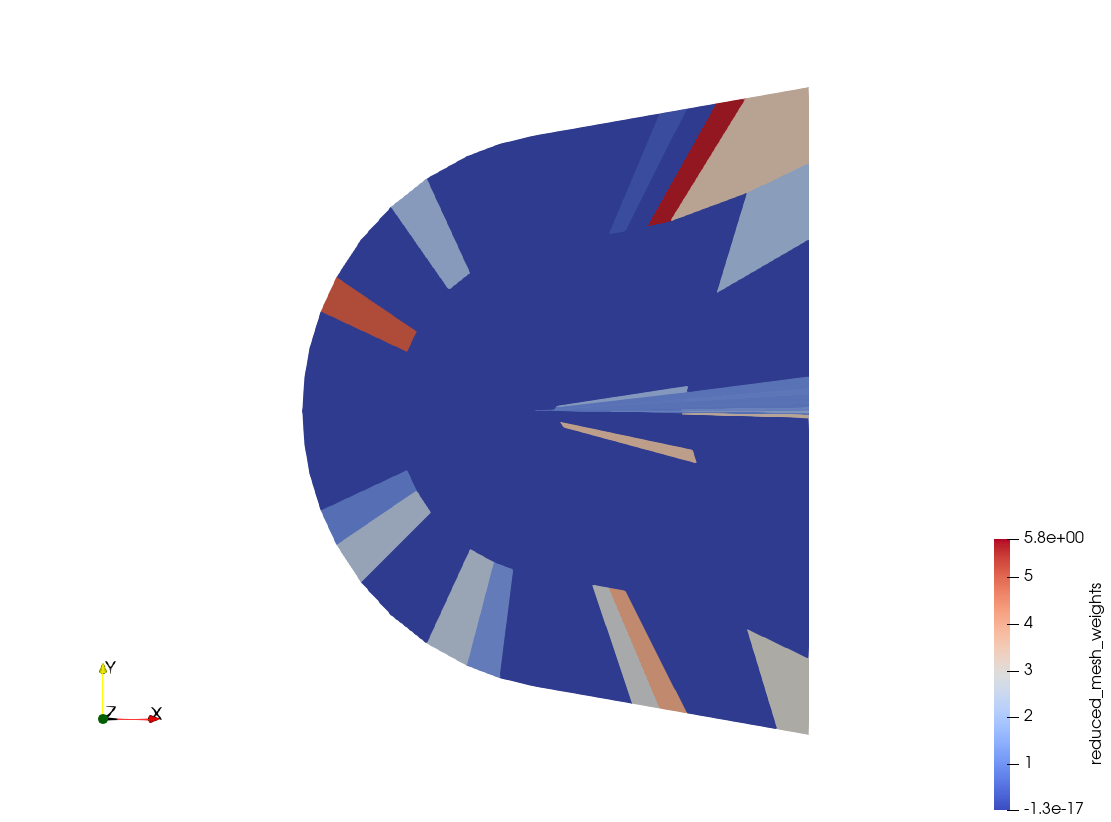}
\caption{ECSW Reduced Mesh Set for the HROM w/ Hyper-DWR}\label{fig:one_param_HROM_DWR_mesh}
\end{figure}

In Table \ref{tab:one_param}, key characteristics and quantities for each of the models are summarized. We can see that the size of the reduced space $n$ increases as hyperreduction is introduced first in the ROM solutions and then into the DWR error indicators, however, all are still significantly smaller than the FOM dimension $N$ of 2240. The reduction in the number of elements used to evaluate the residual and Jacobian is also shown, for the HROM w/o Hyper-DWR only 39 of the 560 mesh elements are required and 45 are used in the HROM w/ hyper-DWR. Figures \ref{fig:one_param_HROM_mesh} and \ref{fig:one_param_HROM_DWR_mesh} show the selected elements that are included in the reduced mesh for the HROM w/o and w/ Hyper-DWR, respectively. It can be seen that both approaches include some of the far field mesh elements ahead of the airfoil as well as some of the elements aft. Both include a majority of the cell directly behind the trailing edge of the airfoil. All three models have accuracy within the same orders of magnitude at the ROM points from the adaptive sampling cycle as well as the 20 points used to evaluate the models ``online". Both HROMs also outperform the ROM online, likely due to their larger POD dimension.

\begin{table}[H]
\caption{One Parameter NACA0012 Airfoil ROM and HROM Important Dimensions, Average ROM Point Error and Average Online Error} \label{tab:one_param}
\centering
\begin{tabular}{|p{4.5cm}||p{0.5cm}|p{1cm}|p{2.5cm}|p{2.5cm}|}
        \hline
         Model & $n$  & $\lVert \boldsymbol{\xi} \rVert_0$ & ROM Error & Online Error\\
        \hline \hline
        ROM & 7 & - &  $2.5162E-5$ & $2.2850E-5$ \\ \hline
        HROM w/o Hyper-DWR & 8 & 39 & $1.5532E-5$ & $1.7671E-5$\\ \hline
        HROM w/ Hyper-DWR & 9 & 45 & $1.7692E-5$ & $1.5497E-5$ \\ \hline
    \end{tabular}
\end{table}

\subsection{Two Parameter NACA0012 Airfoil in Transonic Flow}

\begin{figure}[H]
\centering
\includegraphics[trim = 40 60 60 10, clip, width=0.49\textwidth]{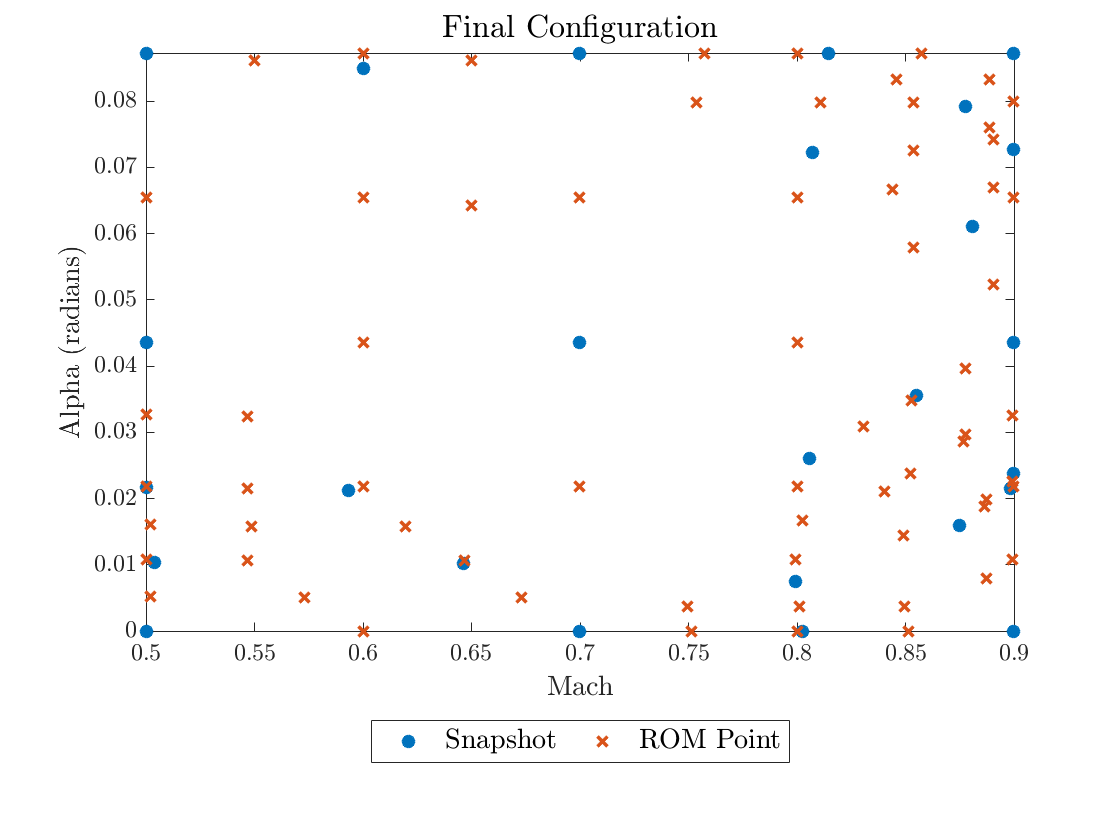}
\includegraphics[trim = 40 60 60 10, clip, width=0.49\textwidth]{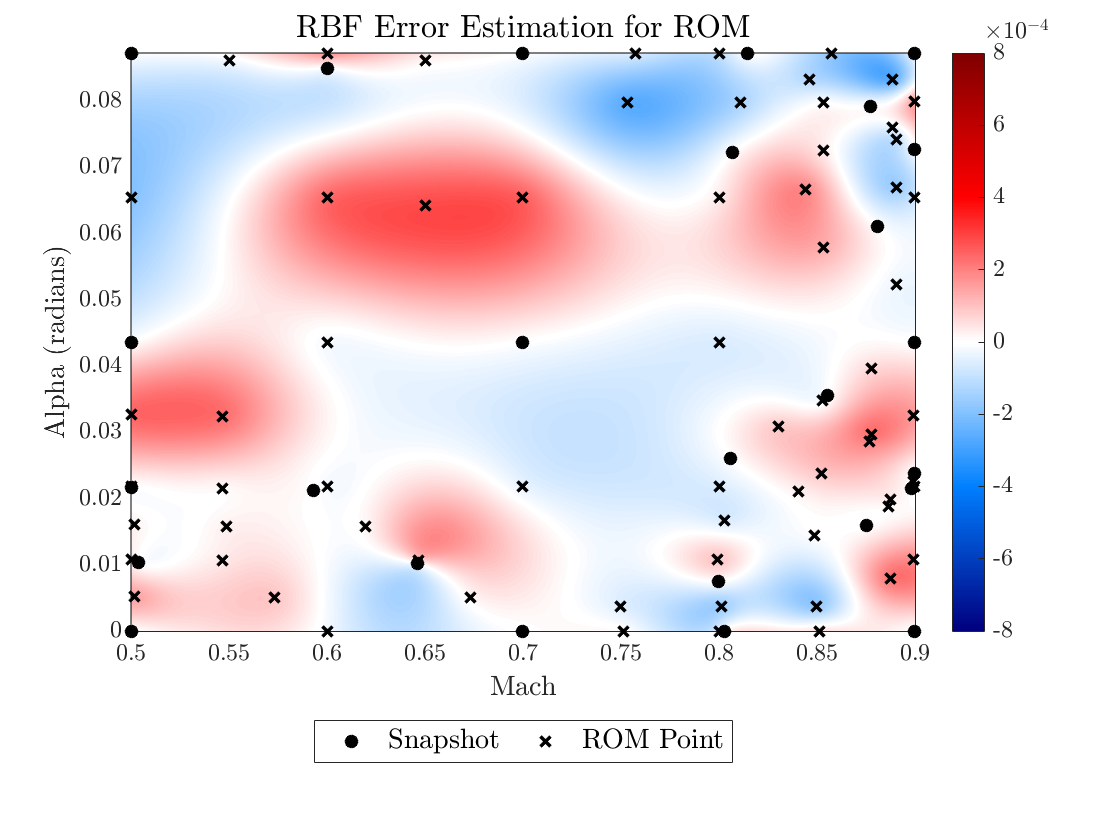}
\captionsetup{justification=centering}
\caption{Final Configuration of the Parameter Space for the ROM and the Interpolated Estimated Error Distribution from the ROM Points}
\label{fig:ROM_trans_adapt}
\end{figure}

\begin{figure}[H]
\centering
\includegraphics[trim = 40 60 60 10, clip, width=0.49\textwidth]{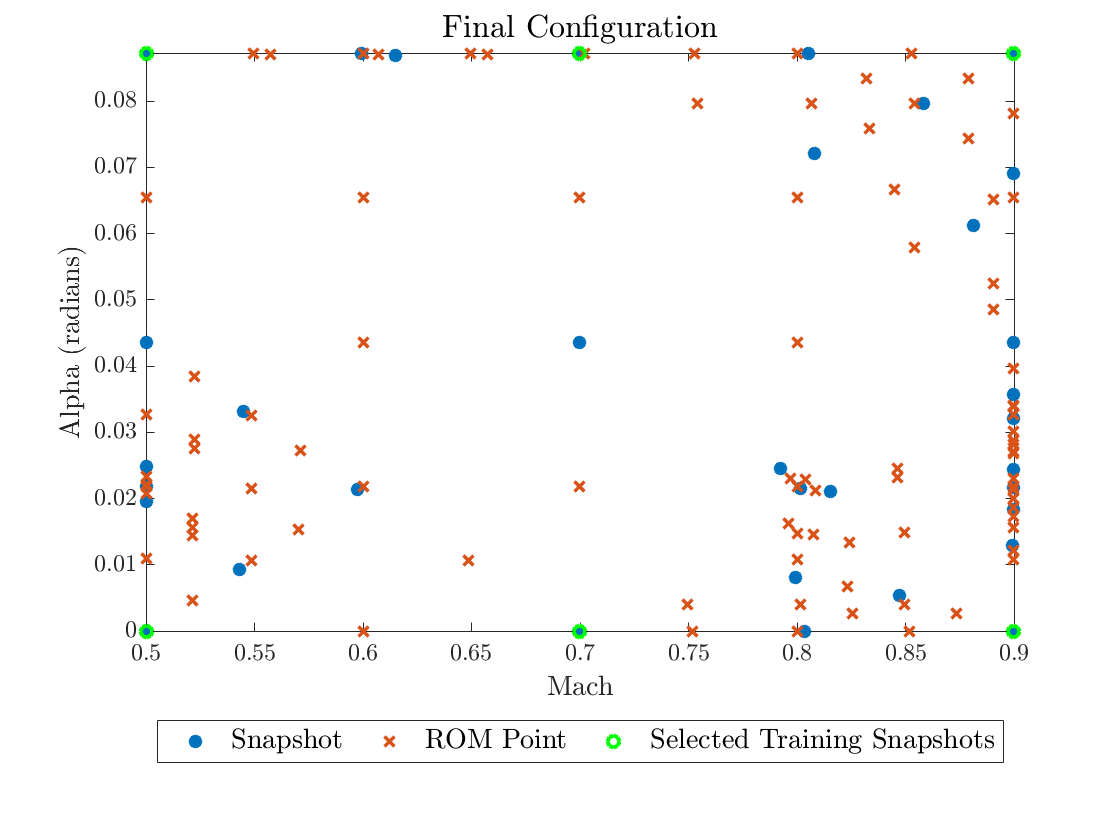}
\includegraphics[trim = 40 60 60 10, clip, width=0.49\textwidth]{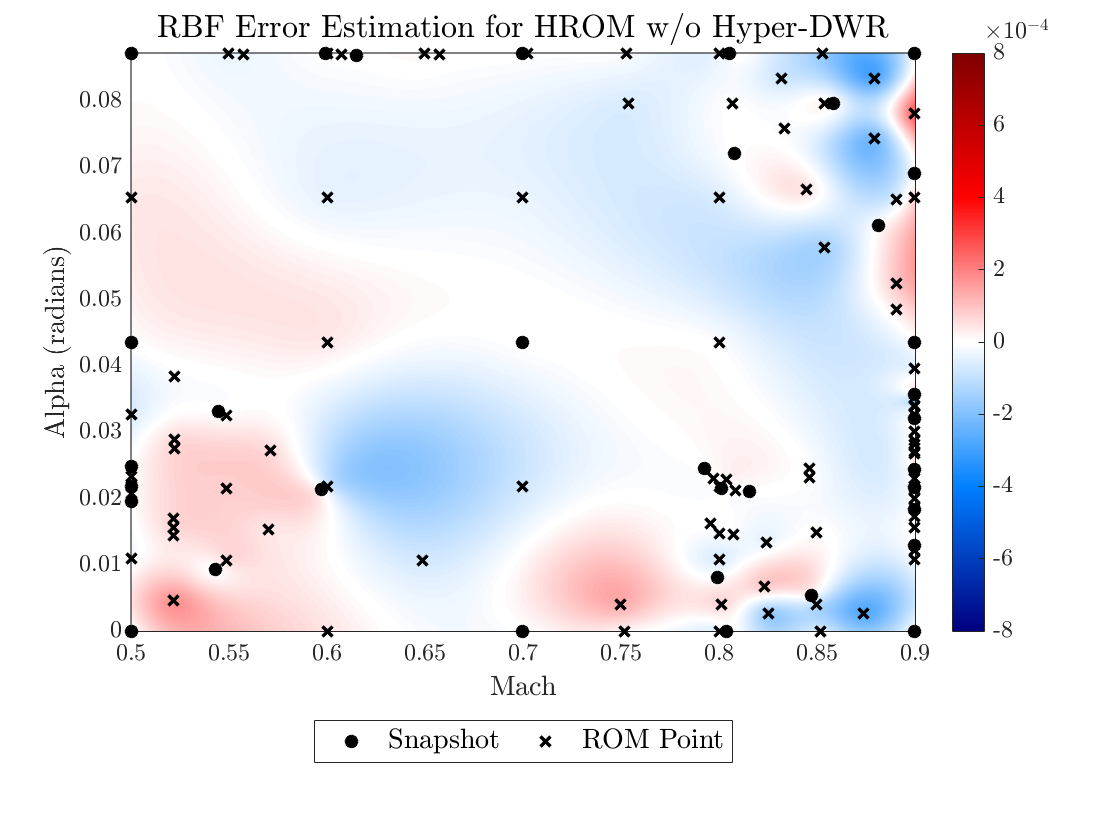}
\captionsetup{justification=centering}
\caption{Final Configuration of the Parameter Space for the HROM w/o Hyper-DWR and the Interpolated Estimated Error Distribution from the ROM Points}
\label{fig:HROM_trans_adapt}
\end{figure}

\begin{figure}[H]
\centering
\includegraphics[trim = 40 60 60 10, clip, width=0.49\textwidth]{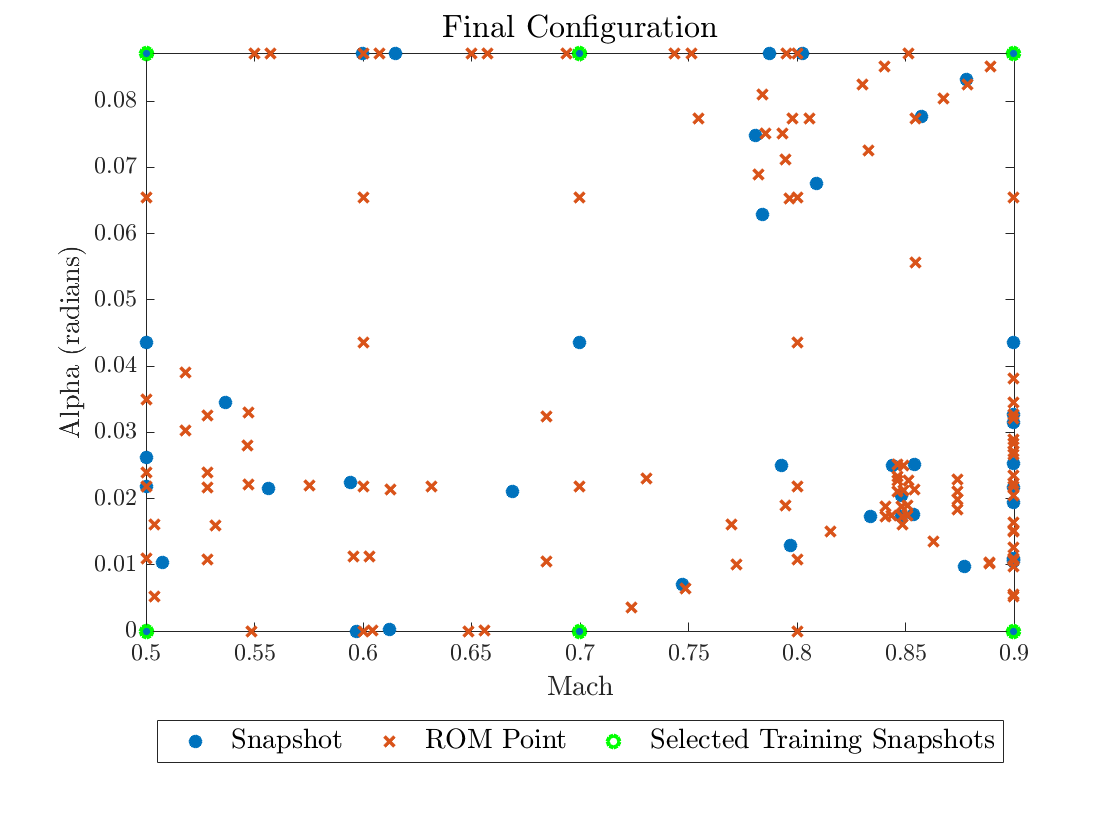}
\includegraphics[trim = 40 60 60 10, clip, width=0.49\textwidth]{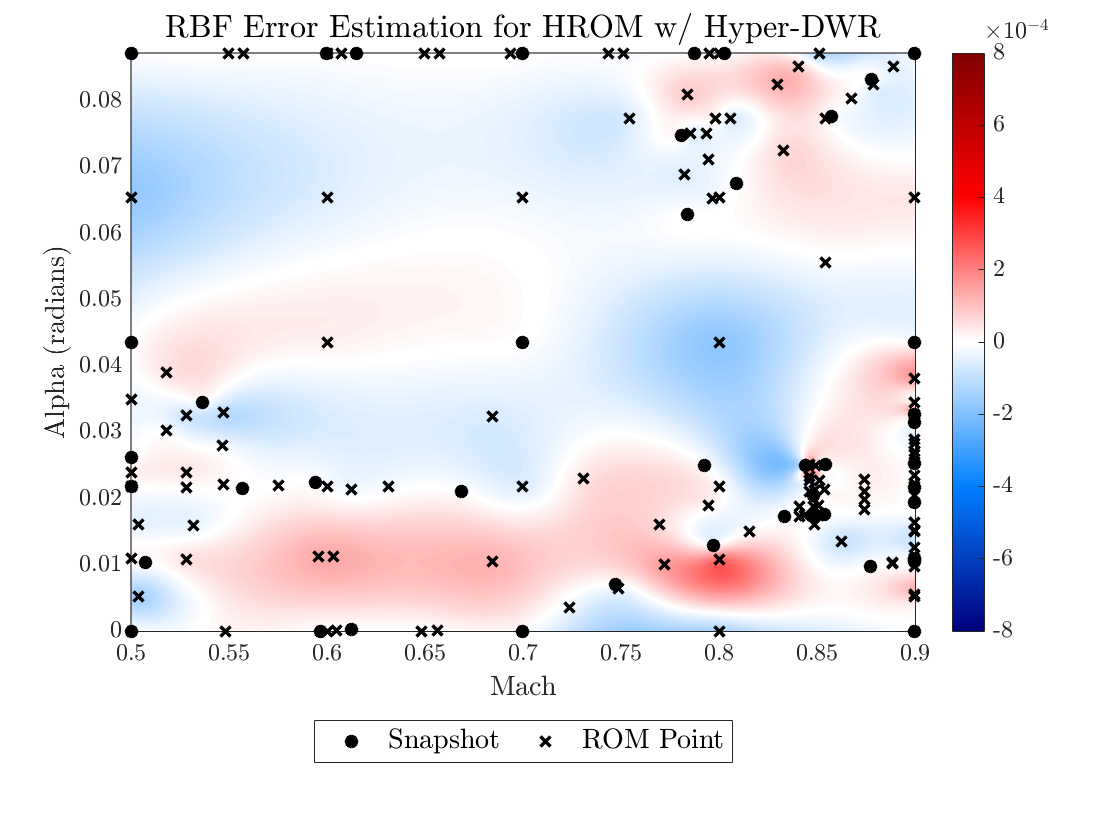}
\captionsetup{justification=centering}
\caption{Final Configuration of the Parameter Space for the HROM w/ Hyper-DWR and the Interpolated Estimated Error Distribution from the ROM Points}
\label{fig:HROM_w_DWR_trans_adapt}
\end{figure}

This test case considers a NACA0012 airfoil in transonic flow with two design parameters. The design parameters are Mach number in the transonic range $[0.5, 0.9]$ and angle of attack $\alpha$ between $[0, 5]^{\circ}$, and the functional of interest is the lift coefficient. The FOM setup is the same as the previous section, meaning there are 560 mesh elements. The adaptive sampling tolerance is increased to $3E-4$ and the sampling procedure is initialized with 9 snapshots arranged in a grid pattern across the parameter space. Figures \ref{fig:ROM_trans_adapt} to \ref{fig:HROM_w_DWR_trans_adapt} show the final configuration of the parameter space for each ROM as well as the error distribution found from the ROM points at the last adaptive sampling iteration. For the ROMs which included hyperreduction in some capacity, the FOM snapshots used to train the ECSW hyperreduction approach are highlighted using green circles. These were selected for training as they are solved for at the very first iteration of the sampling procedure and can therefore be used to train the reduced mesh at every iteration, and are evenly distributed in the parameter domain. Further investigation into how to select these points is left as future work.

It can be seen that both hyperreduction in the ROM solution and the DWR error indicator change the distribution of the interpolated error, and therefore the final snapshot placements for all three models differ. The results show that both the HROM w/o Hyper-DWR and HROM w/ Hyper-DWR sampling procedures tend to favor placing snapshot locations at higher Mach numbers and on the boundaries of the parameter space.

\begin{figure}[H]
\centering
\includegraphics[trim = 40 10 60 10, clip, width=0.7\textwidth]{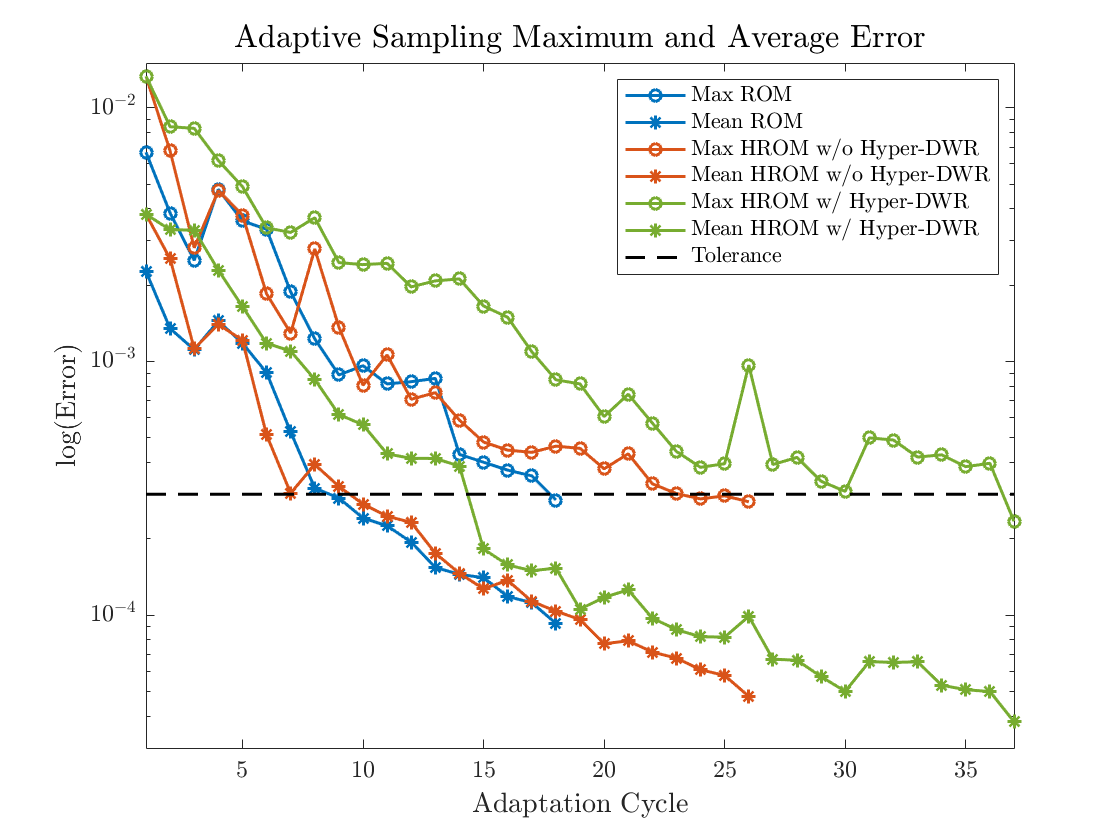}
\captionsetup{justification=centering}
\caption{Maximum and Mean Error Estimate at the ROM Points during each Sampling Iteration for the Three Models} \label{fig:trans_adapt}
\end{figure}

Figure \ref{fig:trans_adapt} shows the maximum and average error at the ROM points during each adaptation cycle for the three models. The standard ROM terminates first and results in the smallest $n$; size of the reduced space. However, we will see later that despite having a smaller reduced basis and thus less unknowns to solve for, the ROM is not less expensive to construct or to use online than the HROMs. This is due to the dimension of the residual and Jacobian which is unaddressed in the original sampling approach. The HROM w/o Hyper-DWR sampling procedure finishes 8 cycles later, due to the hyperreduction, which increases the errors at the ROM points requiring more snapshots to achieve the same tolerance. Since the second DWR is not updated, some of the error introduced by hyperreduction is still not captured in this approach. The HROM w/ Hyper-DWR tends to have higher errors and runs for more adaptation cycles, which is expected as this additional approximation introduced by the hyperreduction is now being tracked via the updated DWR error indicator. Although this produces a larger $n$ dimension, later results will show that this sampling procedure is both more cost-effective and more accurately predicts the online behaviour of the HROM.

\begin{figure}[H]
\centering
\includegraphics[trim = 40 60 60 10, clip, width=0.49\textwidth]{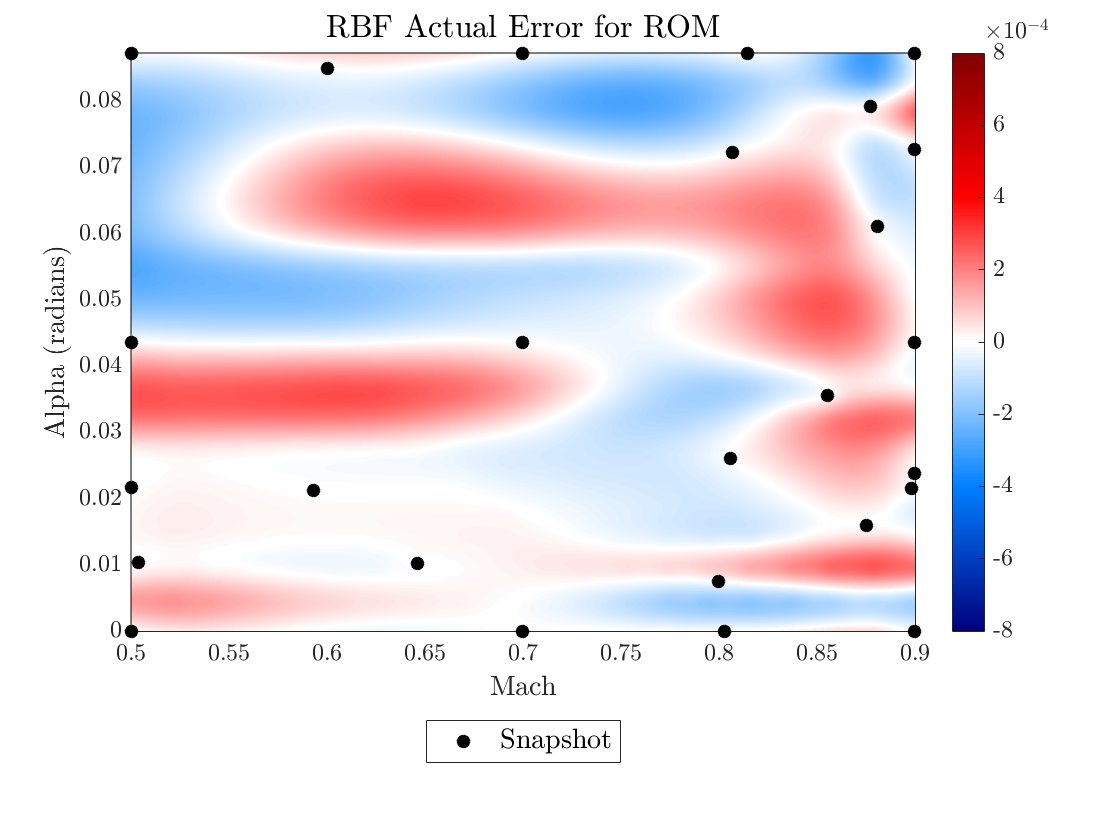}
\includegraphics[trim = 40 60 60 10, clip, width=0.49\textwidth]{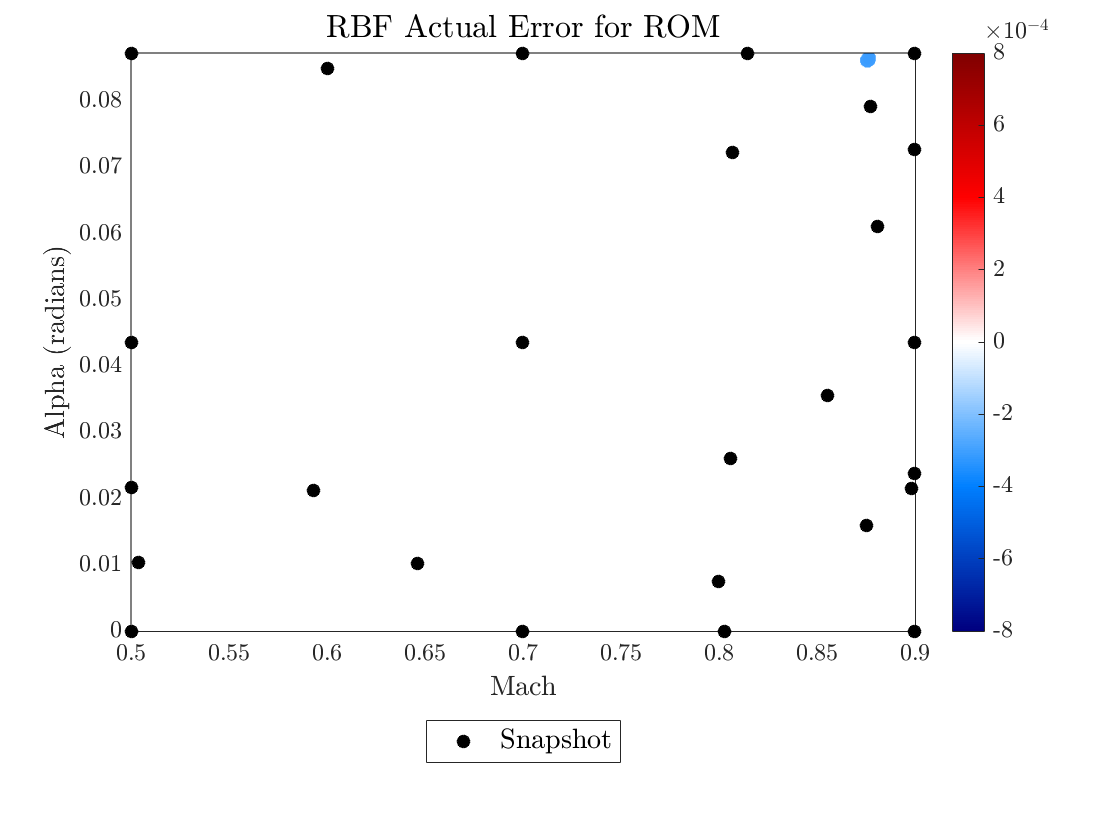}
\caption{True Error for the ROM and Zones of Tolerance Violation} \label{fig:online_ROM_trans}
\end{figure}

\begin{figure}[H]
\centering
\includegraphics[trim = 40 60 60 10, clip, width=0.49\textwidth]{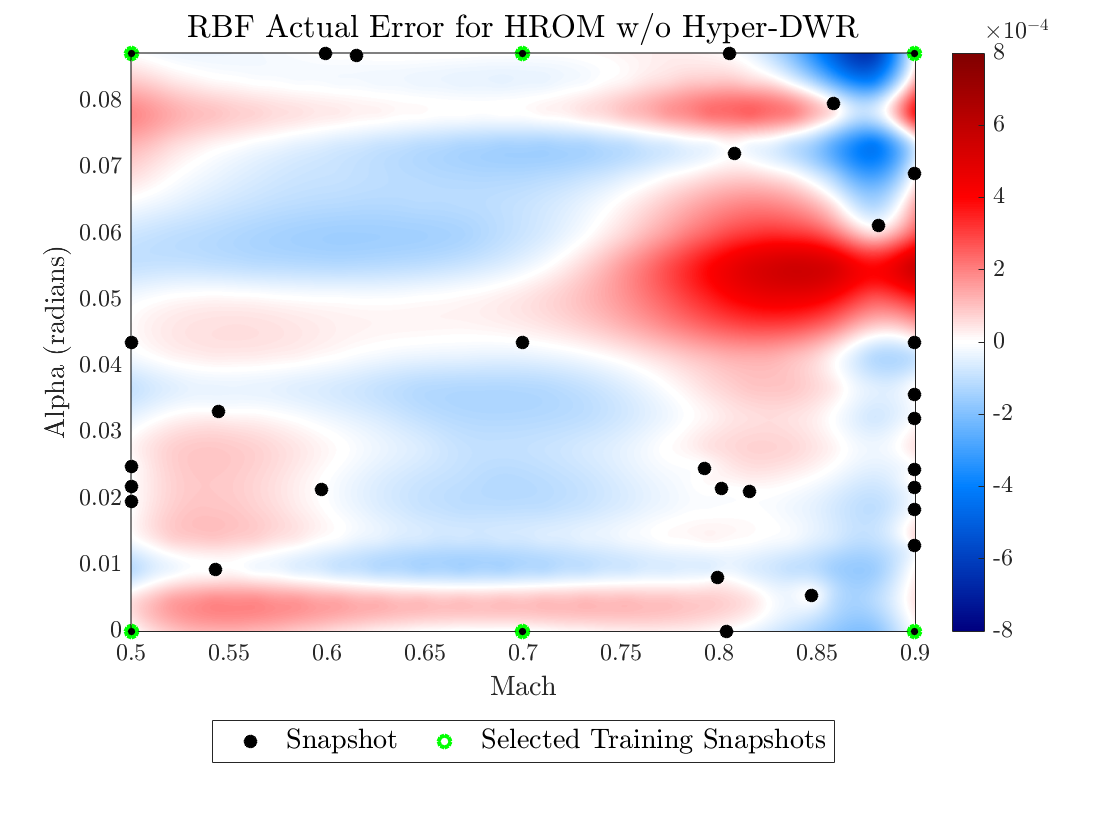}
\includegraphics[trim = 40 60 60 10, clip, width=0.49\textwidth]{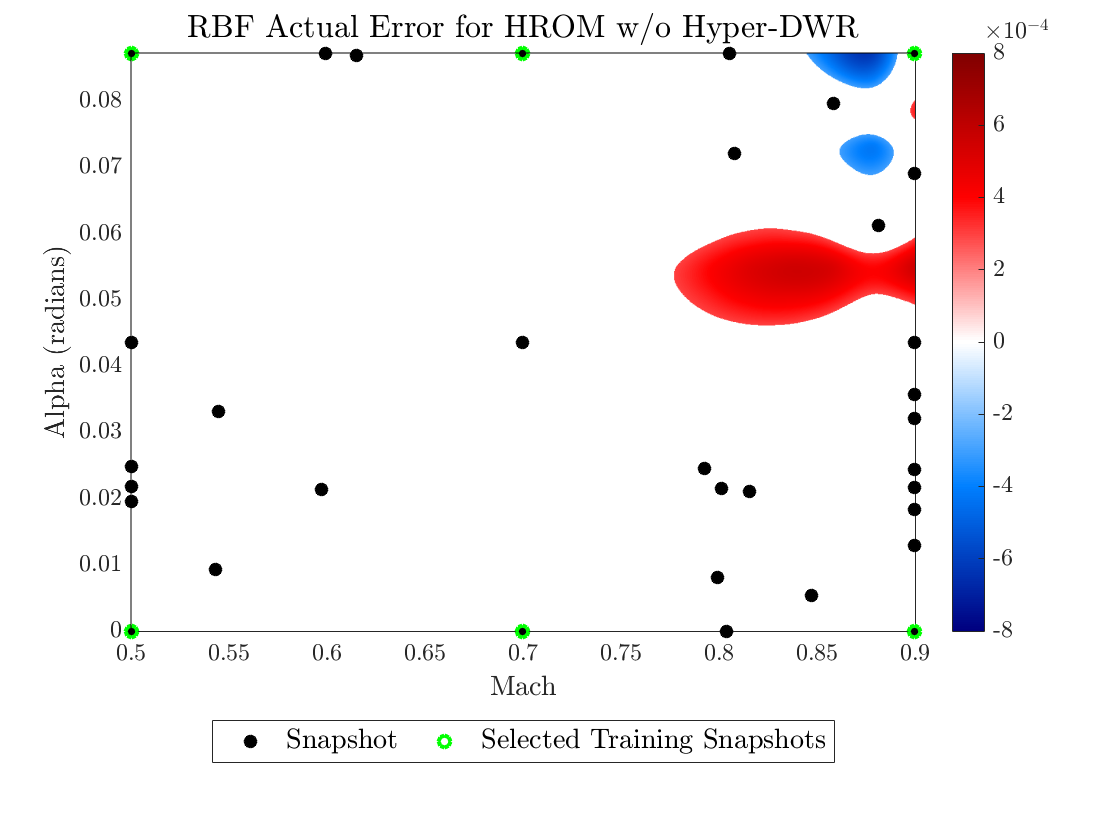}
\caption{True Error for the HROM w/o Hyper-DWR and Zones of Tolerance Violation} \label{fig:online_HROM_trans}
\end{figure}

\begin{figure}[H]
\centering
\includegraphics[trim = 40 60 60 10, clip, width=0.49\textwidth]{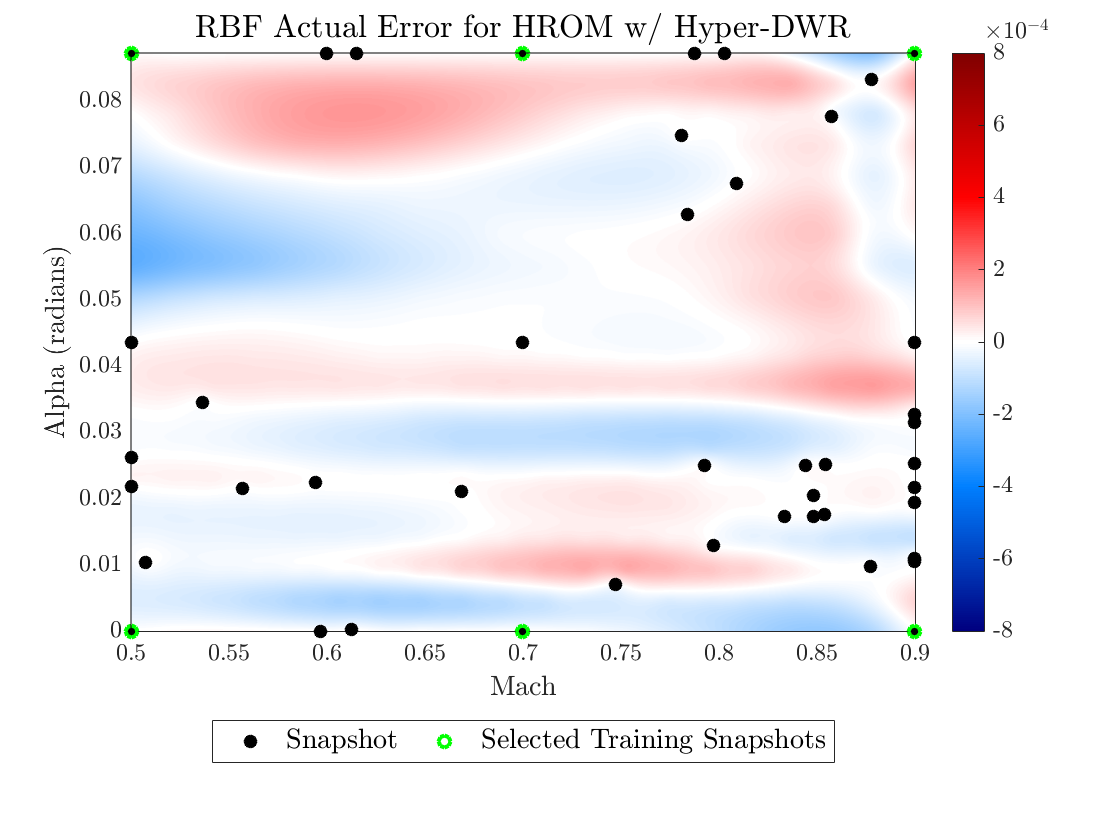}
\includegraphics[trim = 40 60 60 10, clip, width=0.49\textwidth]{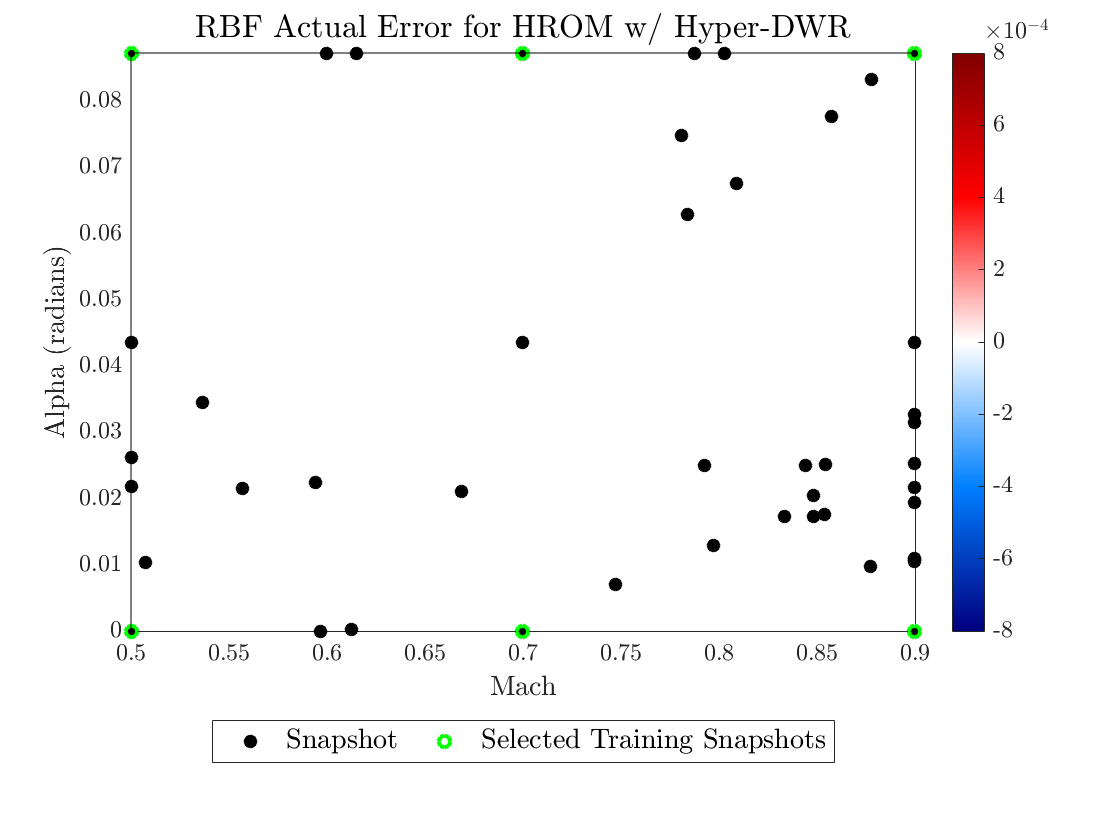}
\caption{True Error for the HROM w/ Hyper-DWR and Zones of Tolerance Violation} \label{fig:online_HROM_w_DWR_trans}
\end{figure}

Figures \ref{fig:online_ROM_trans} to \ref{fig:online_HROM_w_DWR_trans} show the expected online behaviour or ``true" error for each of the models as well as the areas of the domain where the functional error tolerance of $3E-4$ is not met. This distribution is found by evenly distributing 400 points over the parameter space and using an RBF interpolation to approximate the error everywhere. Here we can see the impact of updating the DWR error indicator to include the approximations introduced by hyperreduction. In Figure \ref{fig:online_ROM_trans}, we can see the ROM meets the tolerance bound over almost the entire domain, expect for a very small region of high Mach numbers and high angles of attack. This suggests the adaptive sampling procedure and error estimates closely predicted the online behaviour of the model. In contrast, in Figure \ref{fig:online_HROM_trans}, the HROM w/o Hyper-DWR has multiple larger regions where the tolerance is not met online. The error behaviours in these regions were not predicted by the RBF interpolation from the sampling procedure in Figure \ref{fig:HROM_trans_adapt} and were therefore not considered during the sampling procedure. Once the hyperreduced DWR is used, the tolerance is met over the entire parameter space, as seen in Figure \ref{fig:online_HROM_w_DWR_trans}. This suggests the updated error indicator is correctly capturing the additional error introduced by the hyperreduction, allowing the adaptive sampling procedure to build an accurate HROM with a controllable functional error prediction.

\begin{figure}[H]
\centering
\includegraphics[trim = 300 5 0 80, clip, width=0.66\textwidth]{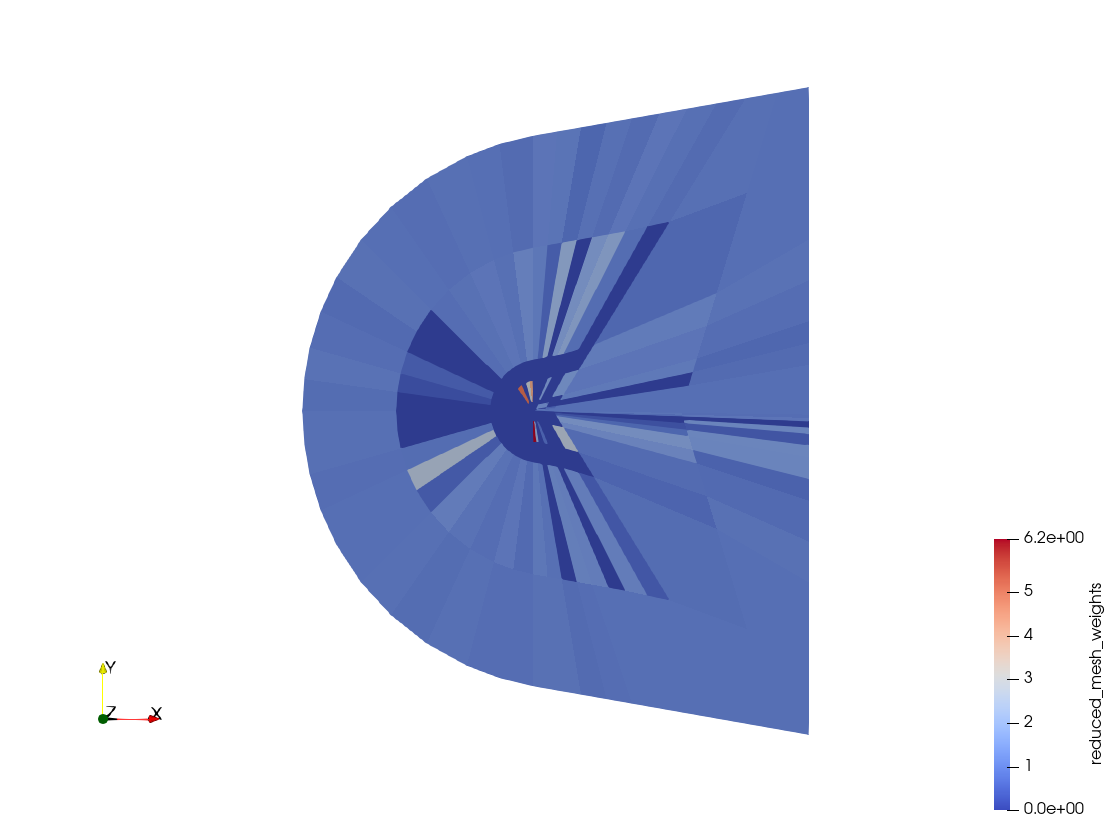}
\caption{ECSW Reduced Mesh Set for the HROM w/o Hyper-DWR}\label{fig:trans_HROM_mesh}
\end{figure}

\begin{figure}[H]
\centering
\includegraphics[trim = 300 5 0 80, clip, width=0.66\textwidth]{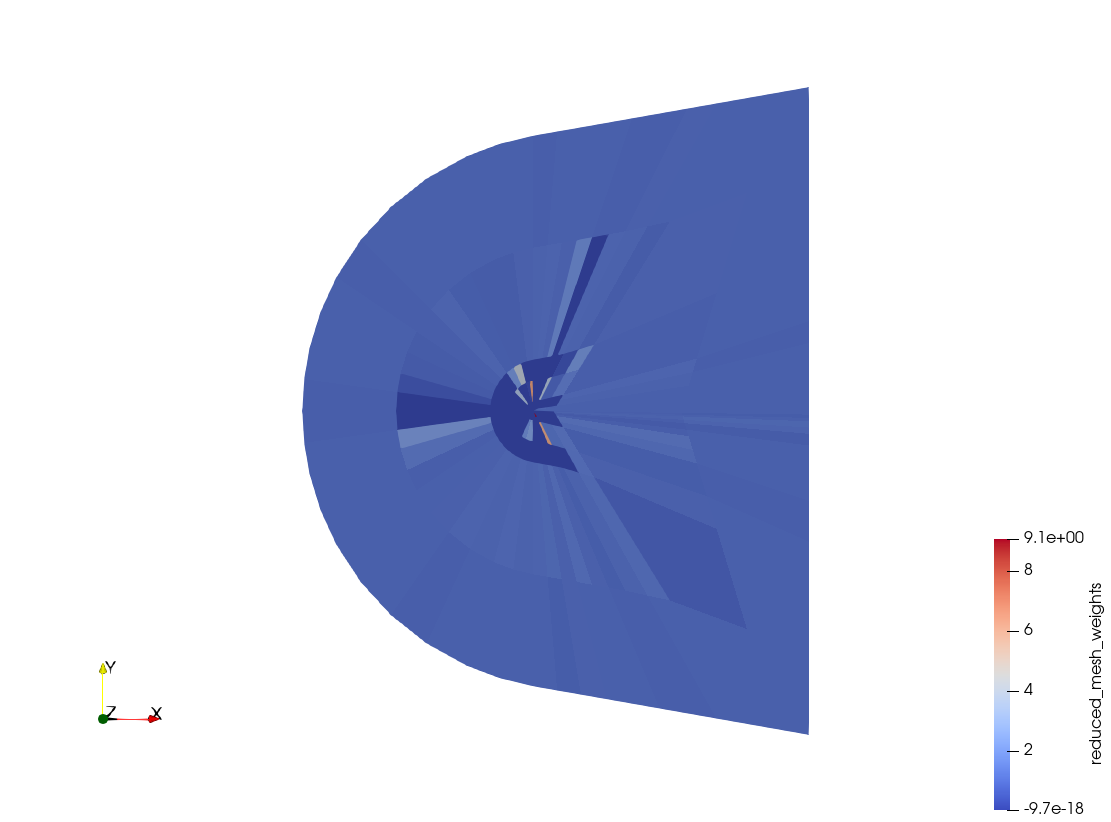}
\caption{ECSW Reduced Mesh Set for the HROM w/ Hyper-DWR}\label{fig:trans_HROM_w_DWR_mesh}
\end{figure}

Table \ref{tab:transonic} summarizes some of the key dimensions of the three models as well as the average error at the ROM points and at the 400 points used to estimate the online behaviour. The HROM w/o Hyper-DWR has reduced-space dimension $n$ of 34 which is 9 larger than the ROM dimension, and 152 elements in the reduced mesh. The HROM w/ Hyper-DWR has a larger $n$ of 40, note however it is still true that $n \ll N$ given that $N$ is 2240, and a reduced mesh with 179 elements down from the 560 FOM mesh. This results in a 68 percent reduction in the number of elements used to evaluate the residual and Jacobian. The reduced mesh of the two HROMs is shown in Figures \ref{fig:trans_HROM_mesh} and \ref{fig:trans_HROM_w_DWR_mesh}. It can be seen that in comparison to the one parameter subsonic case, these reduced meshes tend to weigh more heavily elements closer to the airfoil, which may be because in transonic flow these elements are more relevant to changes in flow solution, and in turn the residual and Jacobian.

\begin{table}[H]
\caption{Two Parameter NACA0012 Airfoil in Transonic Flow ROM and HROM Important Dimensions, Average ROM Point Error and Average Online Error} \label{tab:transonic}
\centering
\begin{tabular}{|p{4.5cm}||p{0.5cm}|p{1cm}|p{2.5cm}|p{2.5cm}|}
        \hline
         Model & $n$  & $\lVert \boldsymbol{\xi} \rVert_0$ & ROM Error & Online Error\\
        \hline \hline
        ROM & 26 & - &  $9.2964E-5$ & $9.4169E-5$ \\ \hline
        HROM w/o Hyper-DWR & 34 & 152 & $4.7918E-5$ & $8.7406E-5$\\ \hline
        HROM w/ Hyper-DWR & 45 & 179 & $3.8079E-5$ & $4.8941E-5$ \\ \hline
    \end{tabular}
\end{table}

\begin{figure}[H]
    \centering         \includegraphics[trim = 30 10 10 10, clip, width=0.8\textwidth]{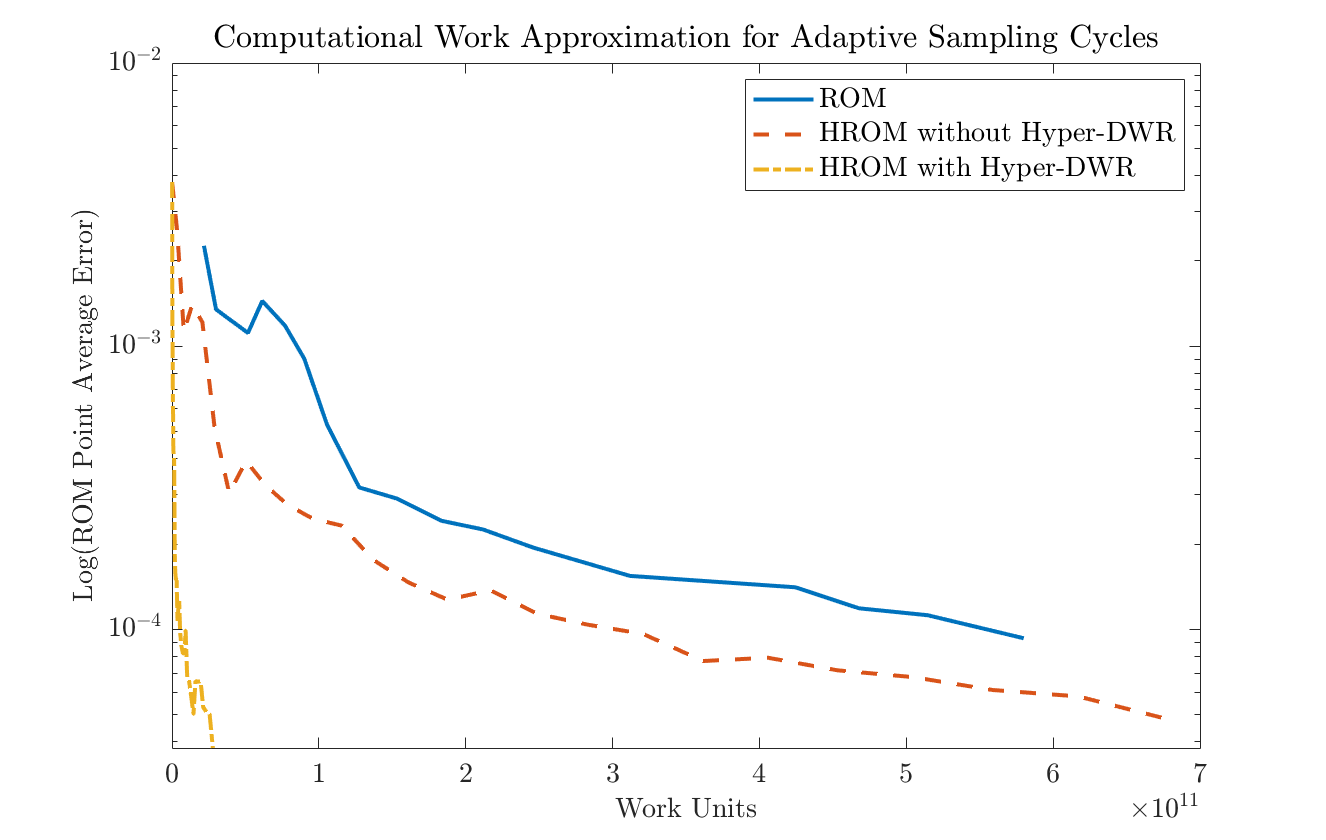}
    \captionsetup{justification=centering}
    \caption{Work Units for the Adaptive Sampling Procedure to Build the Three Models} \label{fig:work_transonic}
\end{figure}

Next, the impact of the hyperreduction on the computational cost of the adaptive sampling procedure can be studied. The work units discussed in section \ref{sec:work} for each of the models are plotted against the average ROM point error, as seen in Figure \ref{fig:work_transonic}. It can be seen that although the POD dimension $n$ and reduced mesh are slightly larger for the HROM w/ Hyper-DWR, each iteration is much less computationally expensive than the other two models due to the savings in evaluating the residual and Jacobian. The addition of hyperreduction into the second DWR indicator results in significant savings as well, since this must be evaluated at all the ROM points from previous iterations, the number of which only grows as the adaptive sampling cycle proceeds. Overall, it appears that by incorporating hyperreduction into the ROM the true computational savings can be realized and with the updated error indicator the online behaviour can be accurately predicted.

\begin{figure}[H]
\centering
\includegraphics[trim = 30 10 60 10, clip, width=0.49\textwidth]{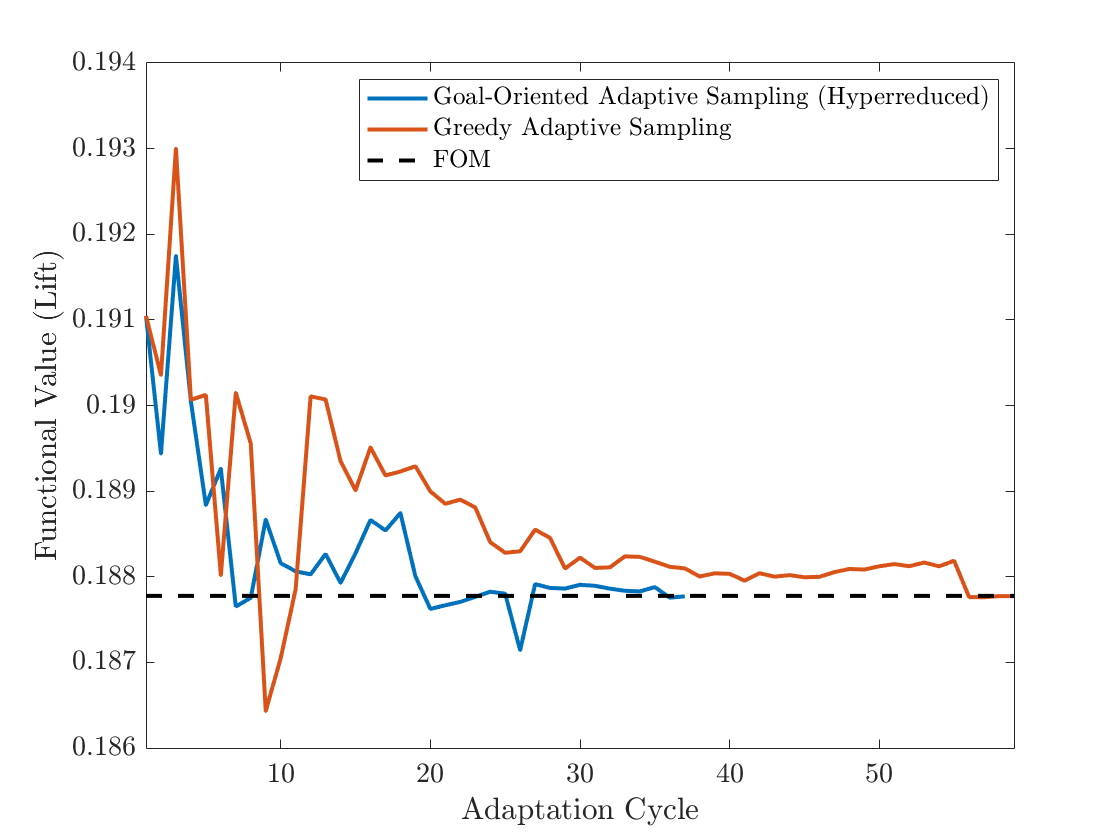}
\includegraphics[trim = 30 10 60 10, clip, width=0.49\textwidth]{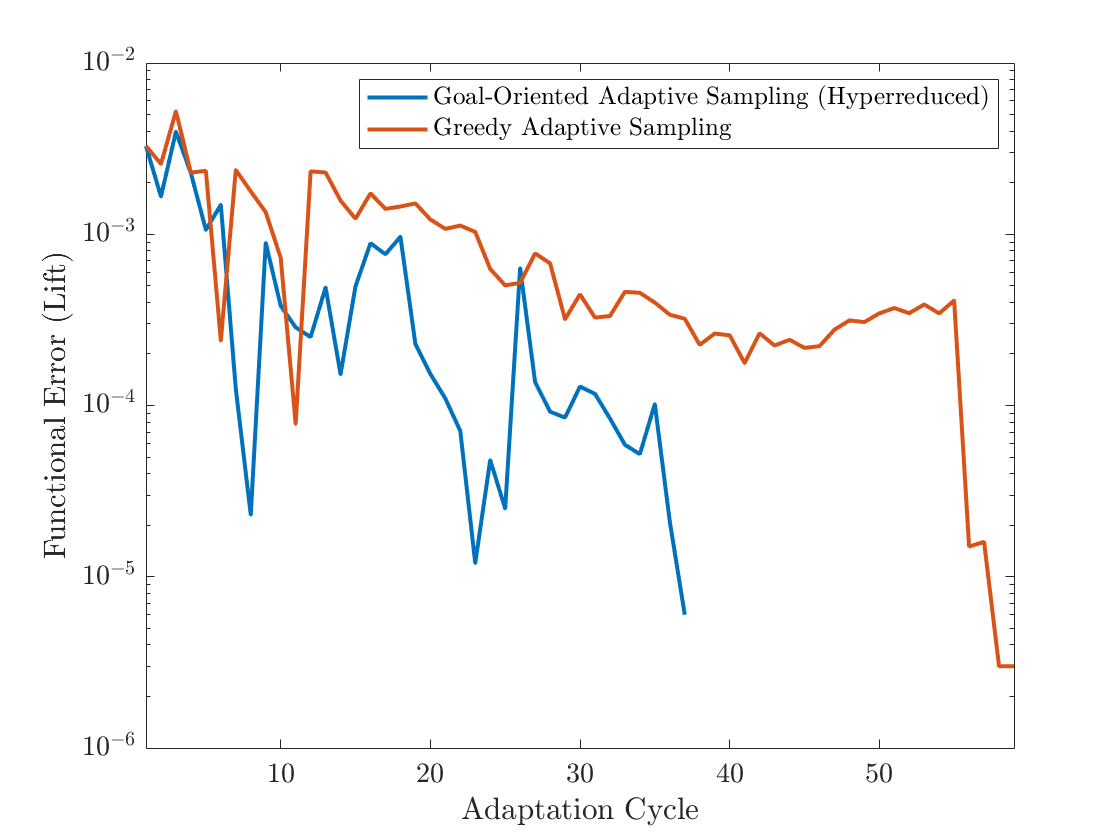}
\includegraphics[trim = 30 10 60 10, clip, width=0.49\textwidth]{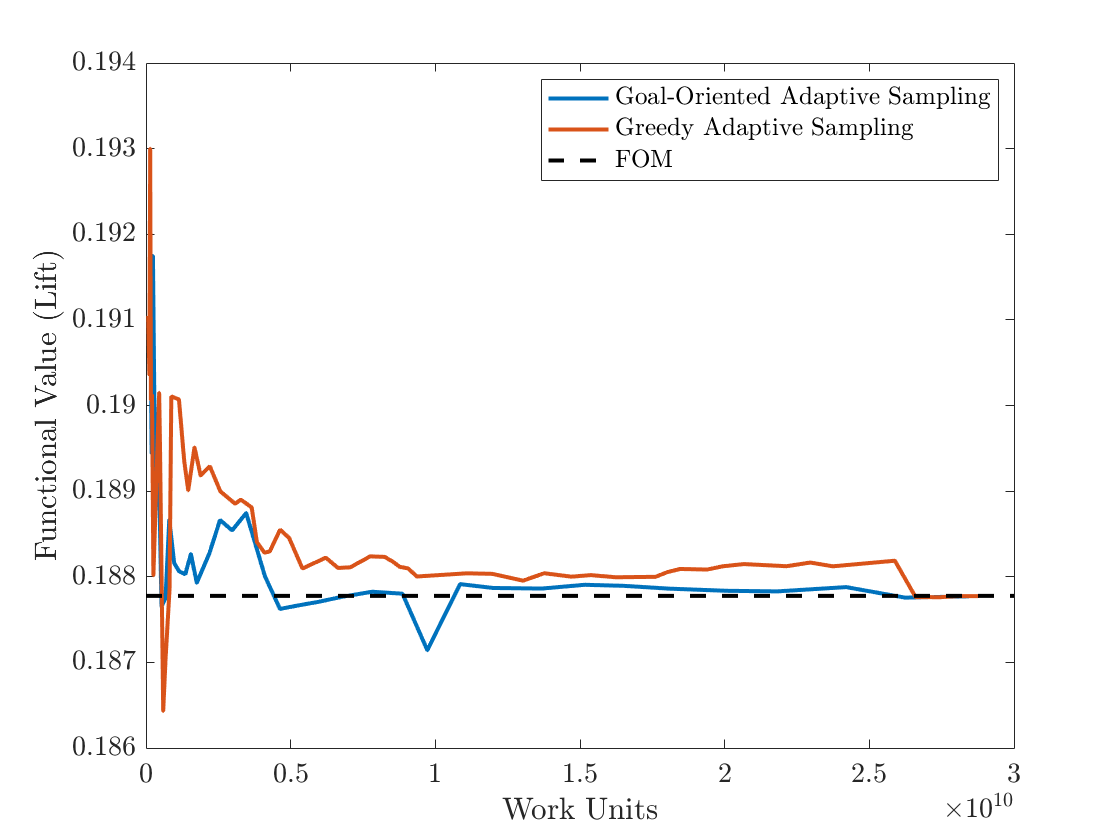}
\includegraphics[trim = 30 10 60 10, clip, width=0.49\textwidth]{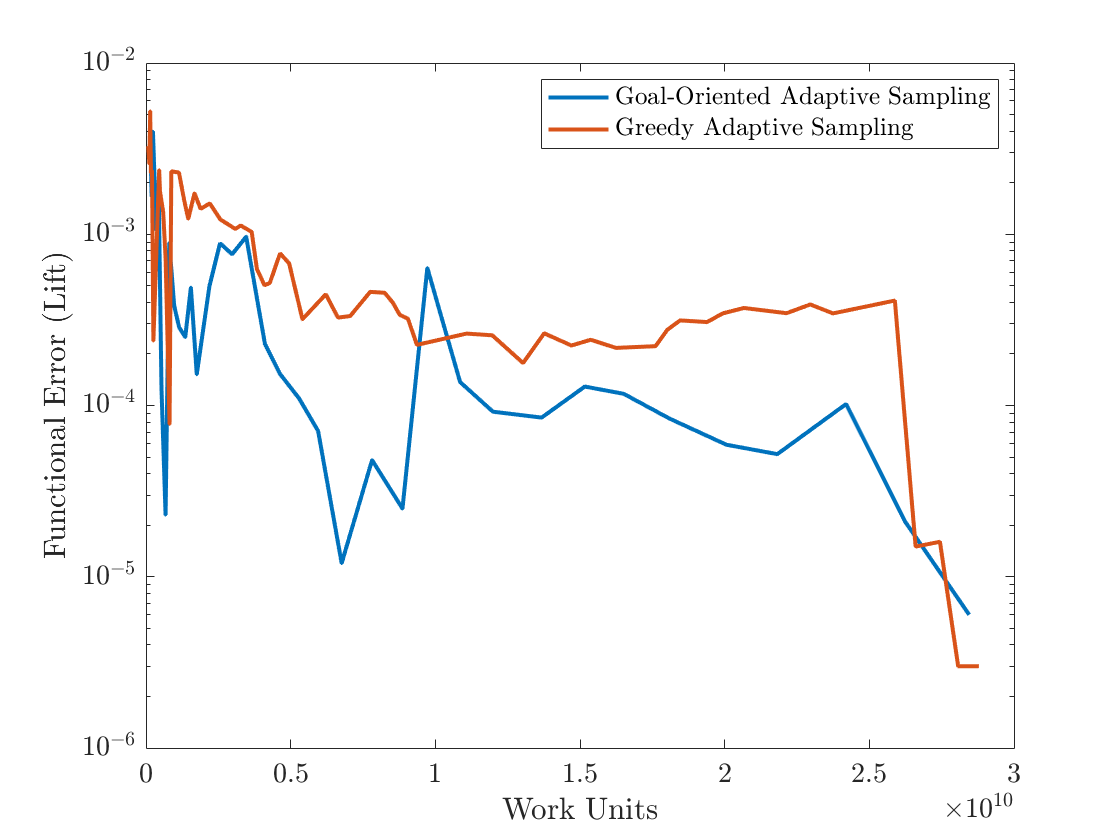}
\captionsetup{justification=centering}
\caption{Comparison of the Functional Value and Error in the Functional between Goal-Oriented and Greedy Adaptive Sampling when Plotted Against Adaptation Cycle and Works Units at Mach = 0.85 and Angle of Attack = 2$^\circ$}
\label{fig:func_tracking}
\end{figure}

Finally, similar to \cite{donJournal}, we compare the hyperreduced goal-oriented approach to a hyperreduced greedy method. To mimic a greedy sampling approach, instead
of sampling the parameter space based on the two DWR error indicators, we instead use the $L_2$-norm of the residual while keeping the rest of the framework unchanged. We update the exit criteria for the sampling procedure via the error tolerance at the ROM points such that the two sampling procedures will run for approximately the same total work units. This is done so that we can compare models that were of equal costs to assemble offline. We select a point in the parameter space, with a Mach of 0.85 and angle of attack of 2$^\circ$, to track the changes in the functional value and the functional error between the HROM and FOM over the sampling procedure. Figure \ref{fig:func_tracking} shows both of these values plotted against the adaptation cycle as well as the work units for the two sampling approaches. It can be seen that the greedy sampling procedure takes more cycles than the goal-oriented approach to converge to the FOM functional value, however the work units incurred by both are similar. This is because each cycle of the goal-oriented approach is more computationally expensive, as it requires assembling the adjoint linear system and the associated linear solves. While in the offline phase the goal-oriented approach is generally comparable to or more expensive than the greedy method, it can more efficiently assemble the reduced-order basis thanks to its intelligent snapshot placements, which results in fewer FOM solutions. The final dimension of the reduced basis for the goal-oriented approach was 45, while that for the greedy was 67. This would mean the reduced-subspace is smaller, and therefore, the resulting HROM would be more cost-effective online.

\subsection{Two Parameter NACA0012 Airfoil in Inviscid Subsonic Flow on a Fine Mesh}

The final test case will be for a two-design parameter NACA0012 airfoil in inviscid subsonic flow. The parameters are the Mach number in the range $[0.3, 0.55]$ and the angle of attack in the range $[0, 3]^{\circ}$. A finer mesh used is used for this test case, which has 2240 cells and 8960 DOFs. The functional of interest will still be the lift coefficient, for which the estimated error tolerance is set to $1E-4$.

\begin{figure}[H]
\centering
\includegraphics[trim = 30 60 60 10, clip, width=0.49\textwidth]{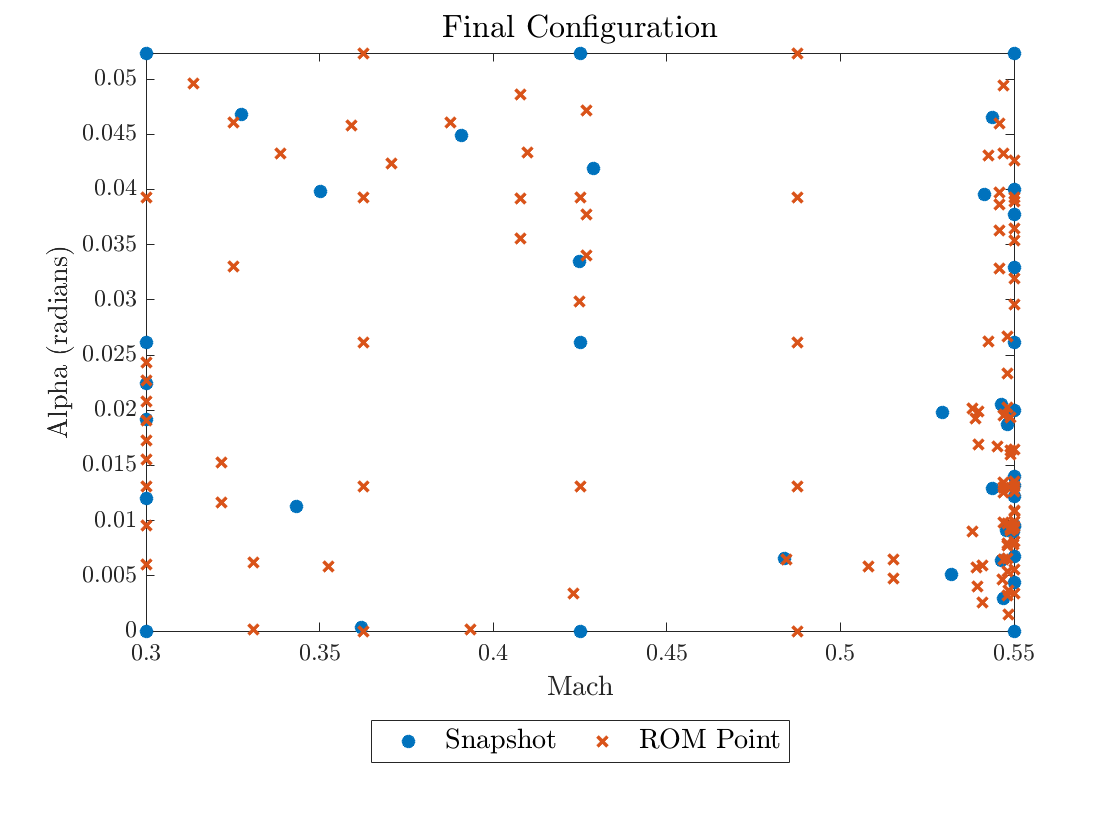}
\includegraphics[trim = 30 60 60 10, clip, width=0.49\textwidth]{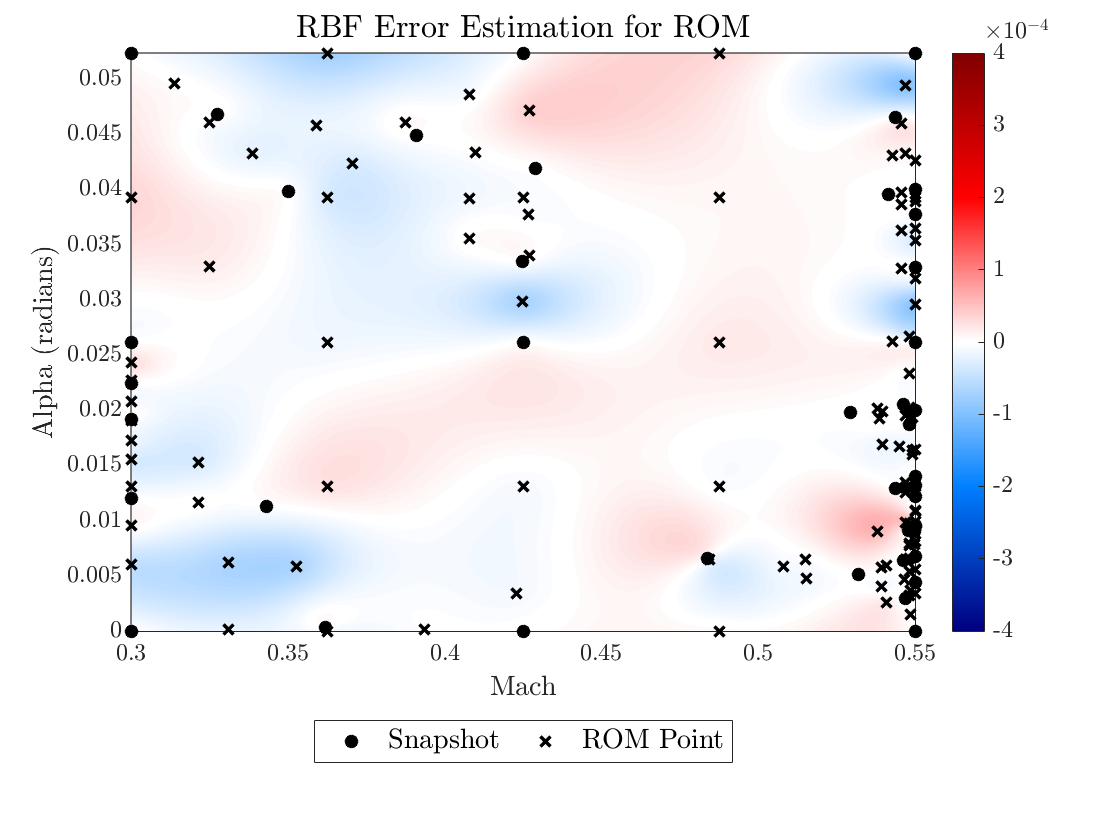}
\captionsetup{justification=centering}
\caption{Final Configuration of the Parameter Space for the ROM and the Interpolated Estimated Error Distribution from the ROM Points}
\label{fig:ROM_sub_fine}
\end{figure}

\begin{figure}[H]
\centering
\includegraphics[trim = 30 60 60 10, clip, width=0.49\textwidth]{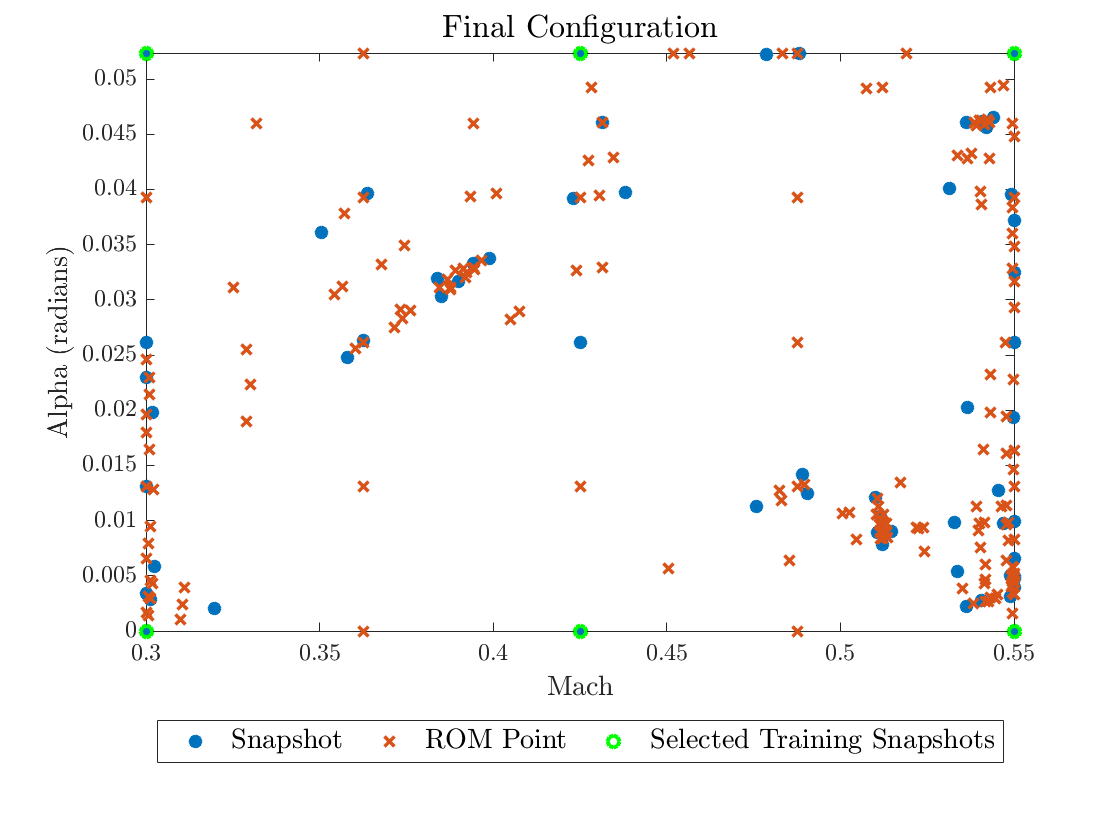}
\includegraphics[trim = 30 60 60 10, clip, width=0.49\textwidth]{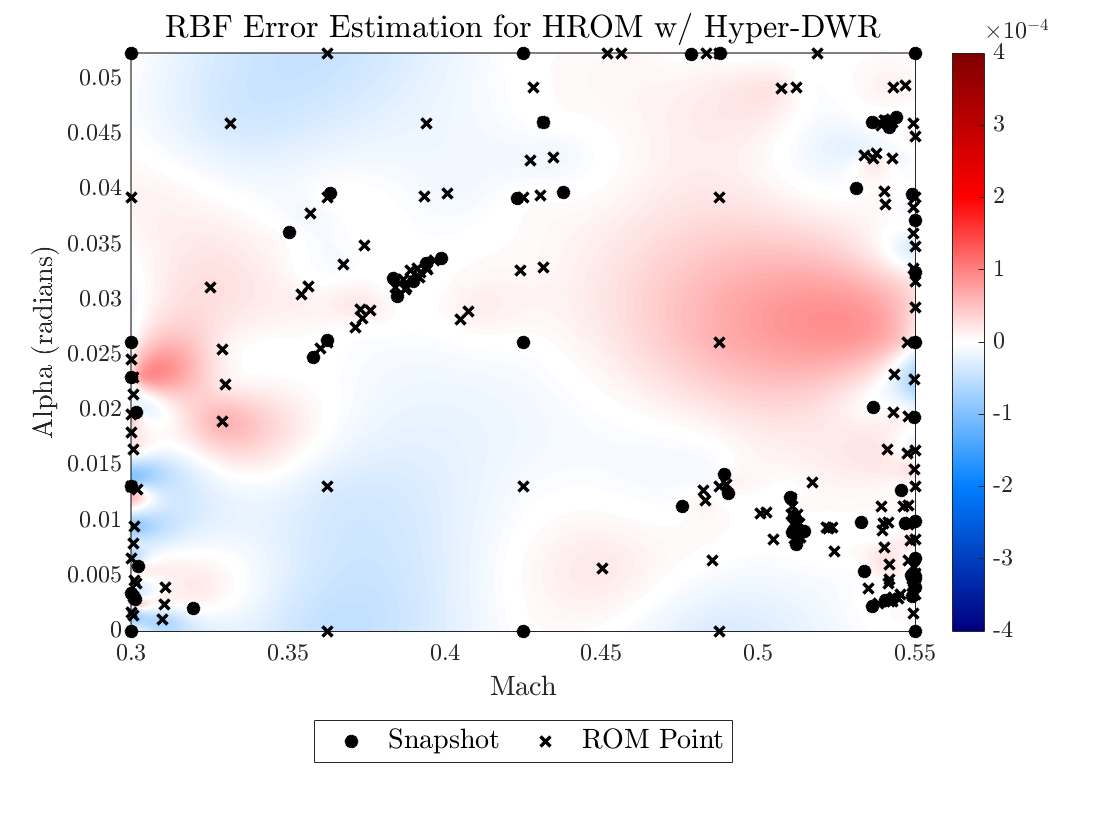}
\captionsetup{justification=centering}
\caption{Final Configuration of the Parameter Space for the HROM w/ Hyper-DWR and the Interpolated Estimated Error Distribution from the ROM Points}
\label{fig:HROM_DWR_sub_fine}
\end{figure}

Figures \ref{fig:ROM_sub_fine} and \ref{fig:HROM_DWR_sub_fine} show the final parameter space configurations and error distributions for the ROM and HROM. It can be seen that both tend to place more snapshots at the highest Mach number, this is likely because more variation is seen in the lift coefficient and flow solutions at higher Mach numbers. 

\begin{figure}[H]
\centering
\includegraphics[trim = 40 10 60 20, clip, width=0.7\textwidth]{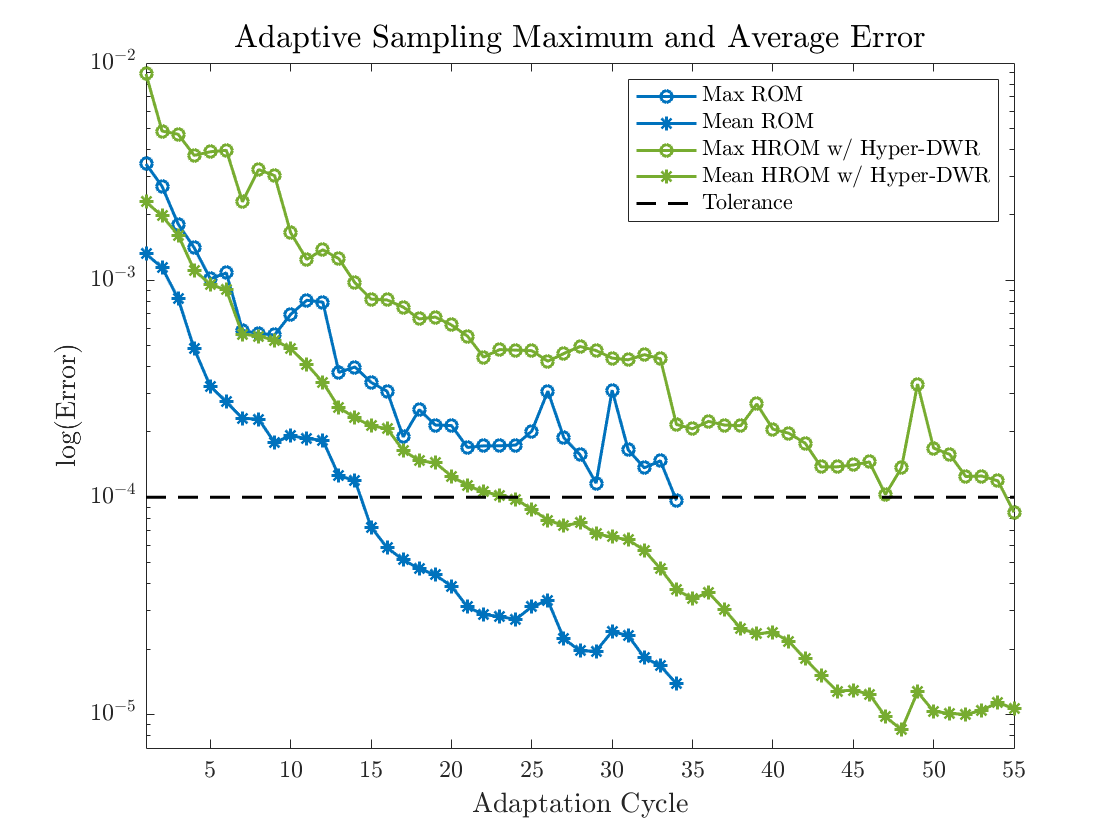}
\captionsetup{justification=centering}
\caption{Maximum and Mean Error Estimate at the ROM Points during each Sampling Iteration for the ROM and the HROM with the Hyperreduced Error} \label{fig:two_fine_adapt}
\end{figure}

Figure \ref{fig:two_fine_adapt} shows the average and maximum ROM point error from the two models over the sampling cycles. We can see that the HROM takes 21 more iterations to achieve the same tolerance. It is interesting to note, however, that the average error dips below the tolerance quite quickly but the maximum ROM point error takes longer to drop and even plateaus at some of the later iterations. This may be due to the training data used to find the ECSW reduced mesh set and weights. The error introduced by hyperreduction cannot be overcome by adding more snapshots into the POD, as the same six points highlighted in Figure \ref{fig:HROM_DWR_sub_fine} are always used for training. Reviewing the snapshot distribution, we can also see the HROM ends up placing more snapshots in the center of the parameter space than the ROM. This may be because the procedure is trying to compensate for the lack of points in this region used for hyperreduction. Future work could look into how this training data may be intelligently selected as well to result in the most accurate reduced mesh set.

\begin{table}[H]
\caption{Two Parameter NACA0012 Airfoil in Subsonic Flow ROM and HROM Important Dimensions, Average ROM Point Error and Average Online Error} \label{tab:two_sub_fine}
\centering
\begin{tabular}{|p{4.5cm}||p{0.5cm}|p{1cm}|p{2.5cm}|p{2.5cm}|}
        \hline
         Model & $n$  & $\lVert \boldsymbol{\xi} \rVert_0$ & ROM Error & Online Error\\
        \hline \hline
        ROM & 42 & - &  $1.3841E-5$ & $2.4963E-5$ \\ \hline
        HROM w/ Hyper-DWR & 63 & 221 & $1.0600E-5$ & $2.0884E-5$ \\ \hline
    \end{tabular}
\end{table}

\begin{figure}[H]
\centering
\includegraphics[trim = 300 0 0 70, clip, width=0.7\textwidth]{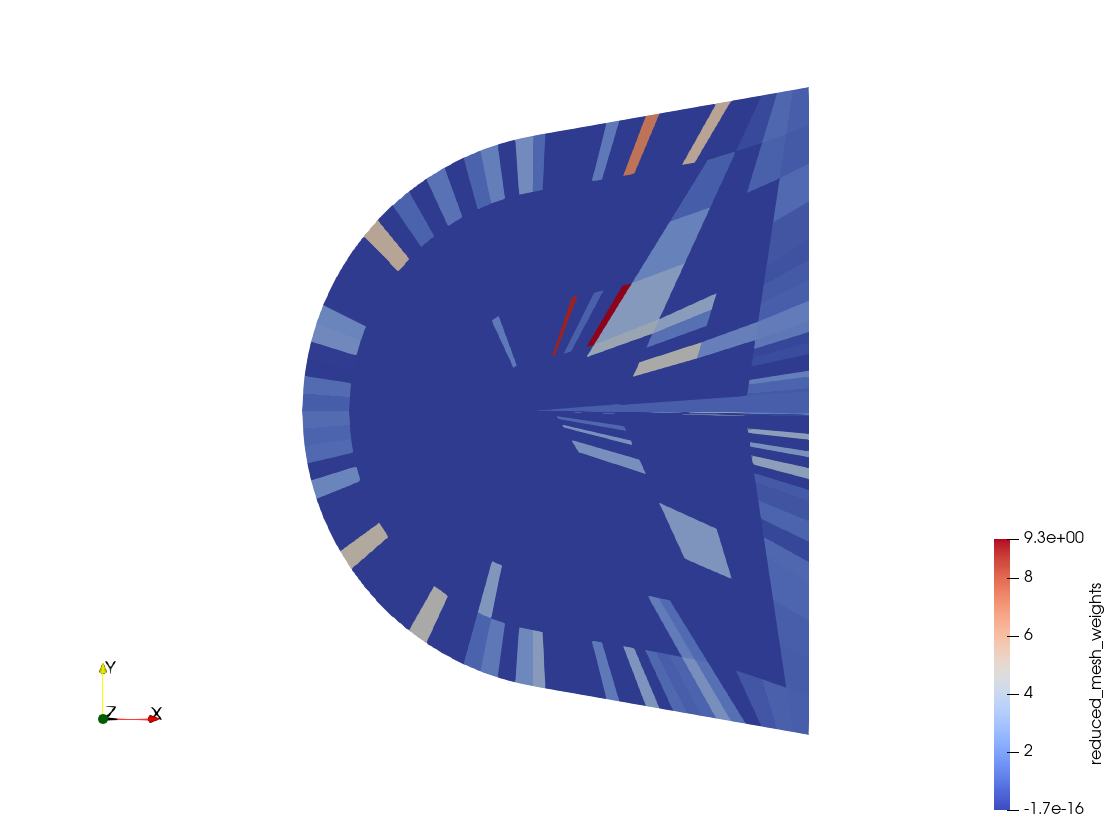}
\caption{ECSW Reduced Mesh set for the HROM w/ Hyper-DWR}\label{fig:sub_HROM_mesh}
\end{figure}

Table \ref{tab:two_sub_fine} summarizes the key dimensions as well as the average ROM point and online error for the two models. This case shows a significant reduction in the number of elements required to accurately evaluate the residual and Jacobian. The FOM mesh contains 2240 elements, while the HROM only requires 221 to solve for these quantities in the Newton iteration in Equation \ref{eq:lspg_newton} and in the DWR error indicator in Equation \ref{eq:hyp_DWR_final}. This is a 90.1\% reduction in the number of mesh elements. We can note that the savings in the subsonic case are even greater than the two-parameter transonic case. This may be because of the mesh size, having more cells resulting in an even greater reduction in the mesh, or the complexity of the problem. Figure \ref{fig:sub_HROM_mesh} shows the reduced mesh set for the HROM.

\begin{figure}[H]
\centering
\includegraphics[trim = 30 60 60 10, clip, width=0.49\textwidth]{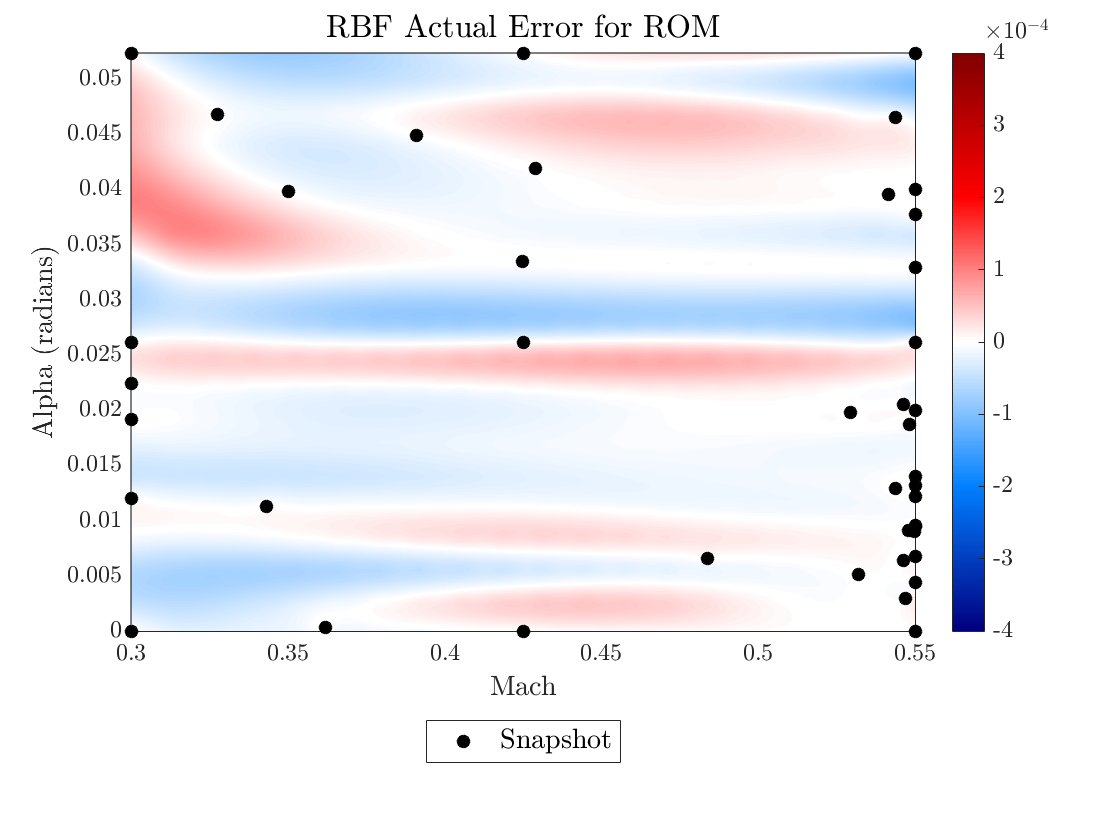}
\includegraphics[trim = 30 60 60 10, clip, width=0.49\textwidth]{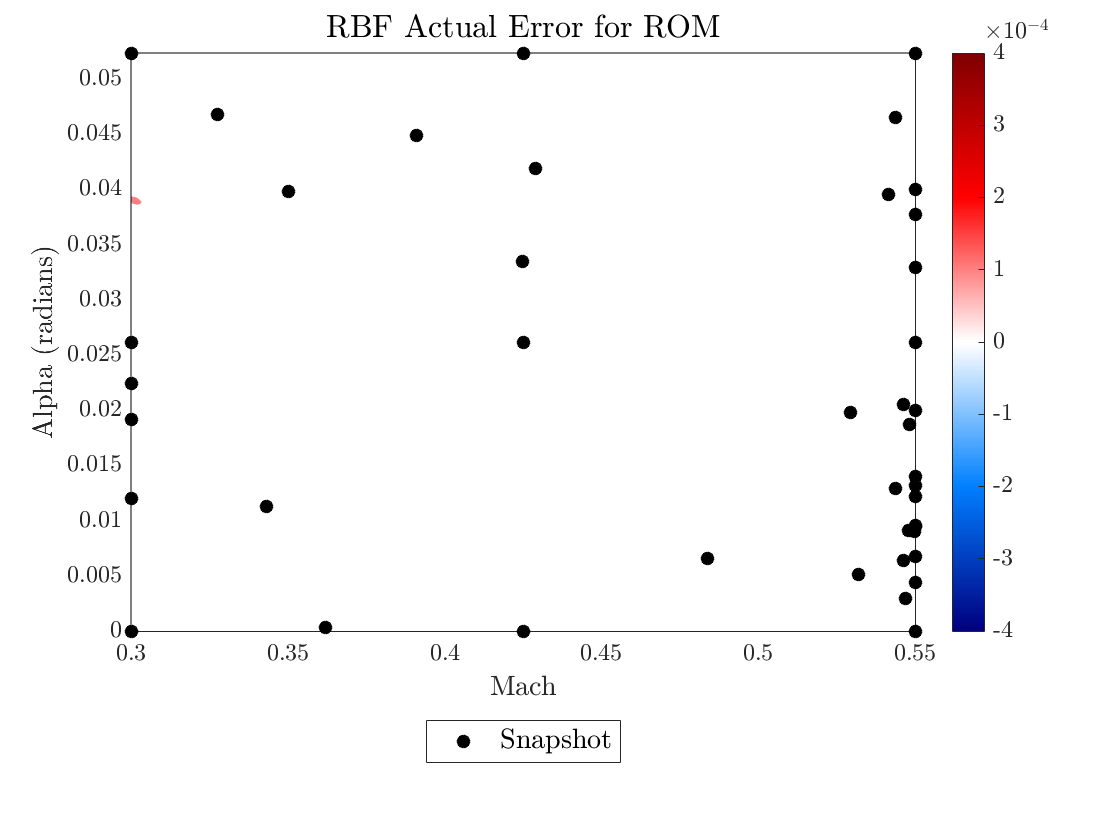}
\caption{True Error for the ROM and Zones of Tolerance Violation} \label{fig:online_ROM_sub}
\end{figure}

\begin{figure}[H]
\centering
\includegraphics[trim = 30 60 60 10, clip, width=0.49\textwidth]{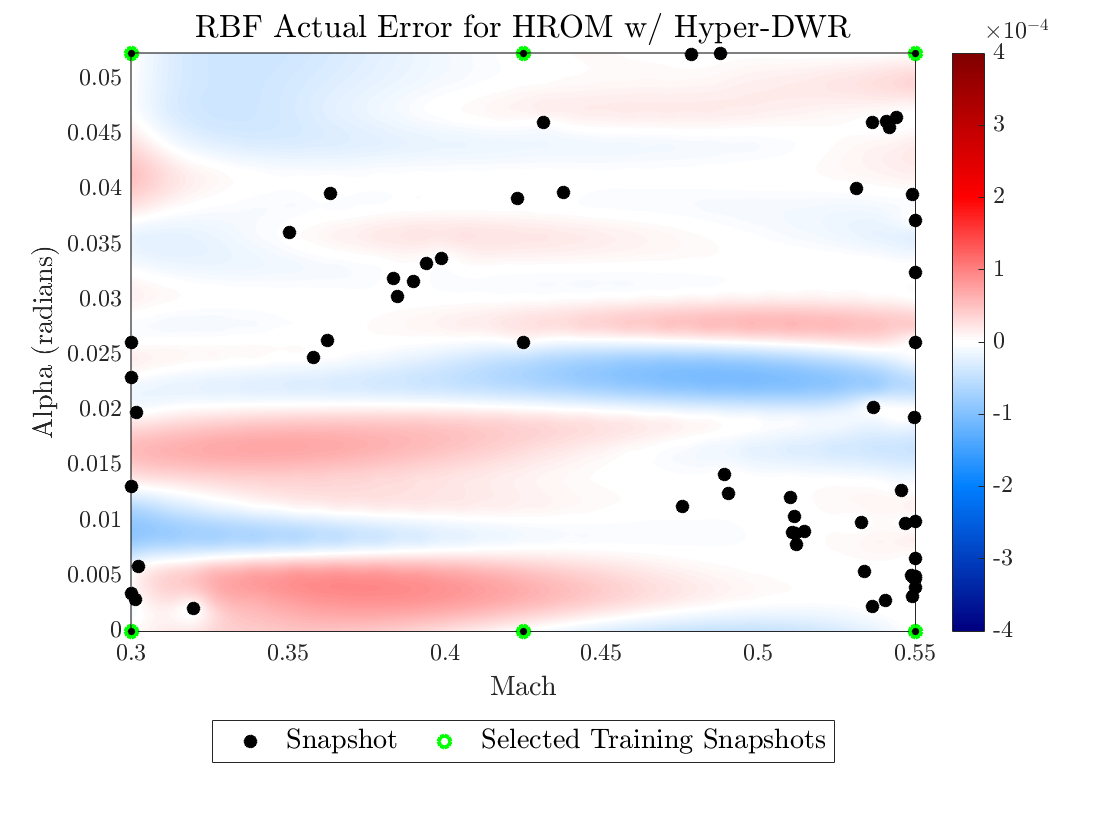}
\includegraphics[trim = 30 60 60 10, clip, width=0.49\textwidth]{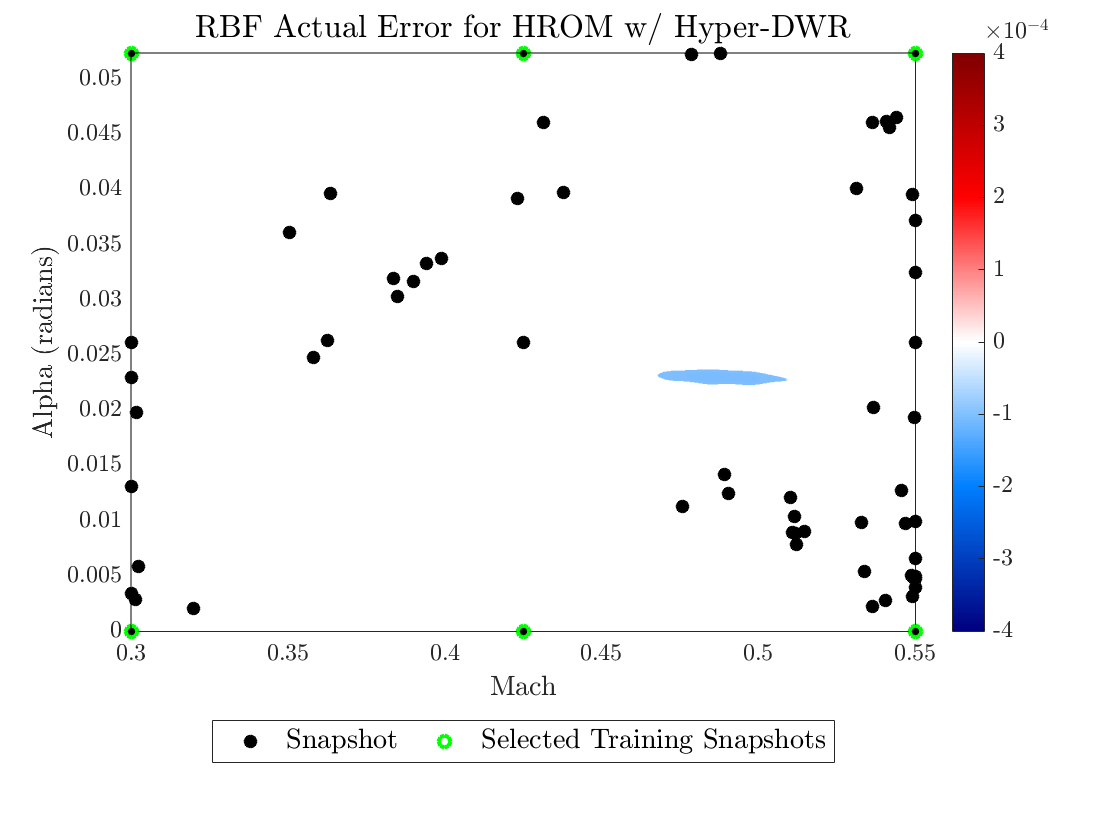}
\caption{True Error for the HROM w/ Hyper-DWR and Zones of Tolerance Violation} \label{fig:online_HROM_sub}
\end{figure}

Figures \ref{fig:online_ROM_sub} and \ref{fig:online_HROM_sub} show the true error distributions from the ROM and HROM, respectively. We can see both models perform quite well onlinee and only have small regions wherethe tolerance bound is not met.

\begin{figure}[H]
    \centering         \includegraphics[trim = 30 10 10 10, clip, width=0.8\textwidth]{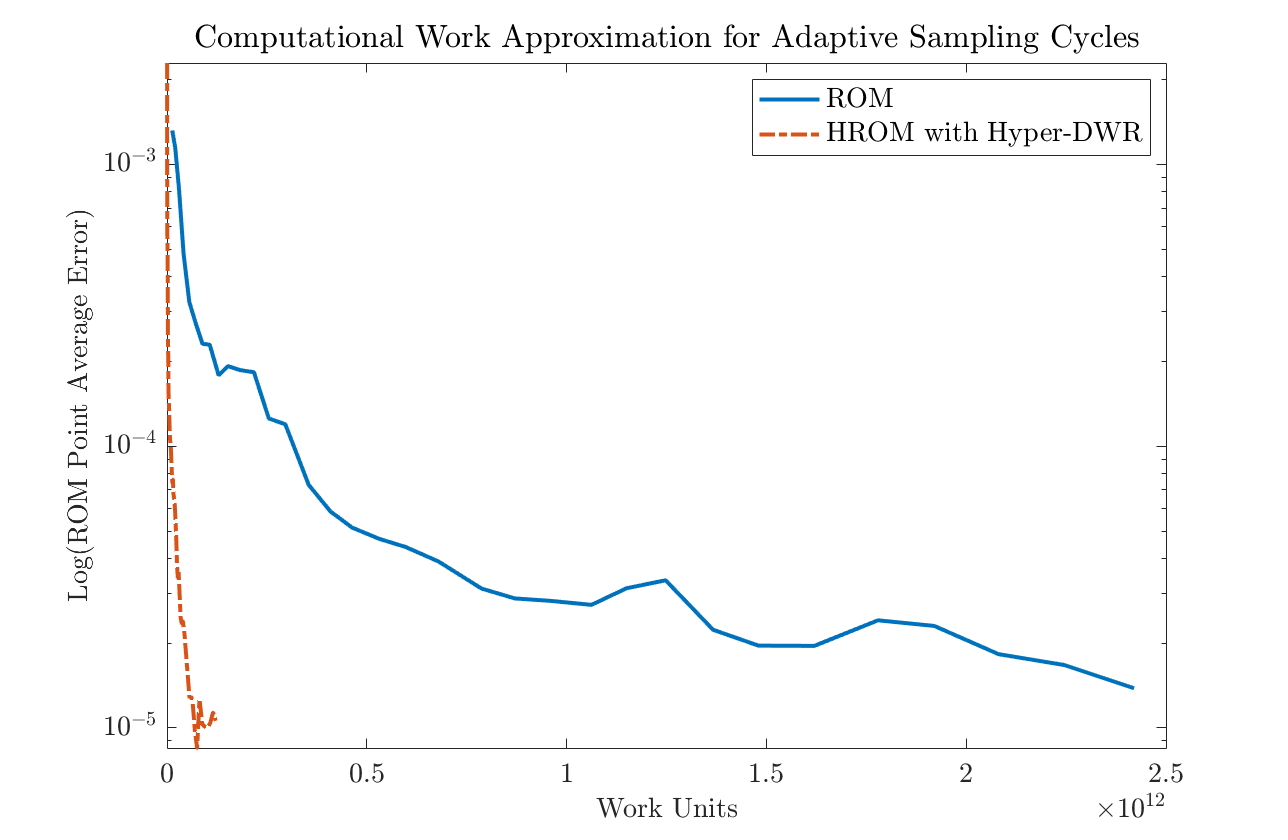}
    \captionsetup{justification=centering}
    \caption{Work Units for the Adaptive Sampling Procedure to Build the Three Models} \label{fig:work_subsonic}
\end{figure}

The work units for each model are plotted in Figure \ref{fig:work_subsonic}. We can see the HROM with the hyperreduced DWR error indicator costs much less to achieve the same average ROM point error. This means the HROM is both less expensive to assemble offline and more computationally efficient during online evaluations through the approximation of the residual and Jacobian. 

\begin{figure}[H]
\centering
\includegraphics[trim = 40 60 60 10, clip, width=0.6\textwidth]{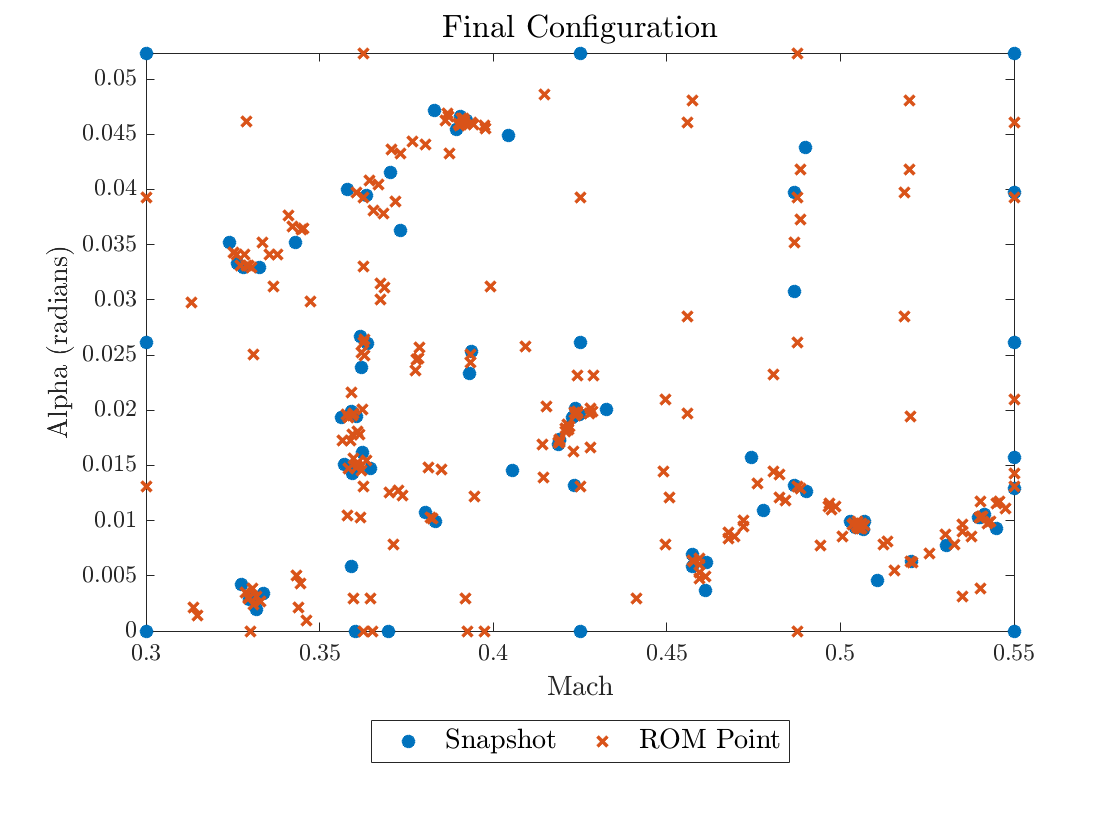}
\captionsetup{justification=centering}
\caption{Final Configuration of the Parameter Space for the HROM Built Using the Greedy Adaptive Sampling Approach}
\label{fig:greedy_config}
\end{figure}

\begin{figure}[H]
\centering
\includegraphics[trim = 40 60 60 10, clip, width=0.6\textwidth]{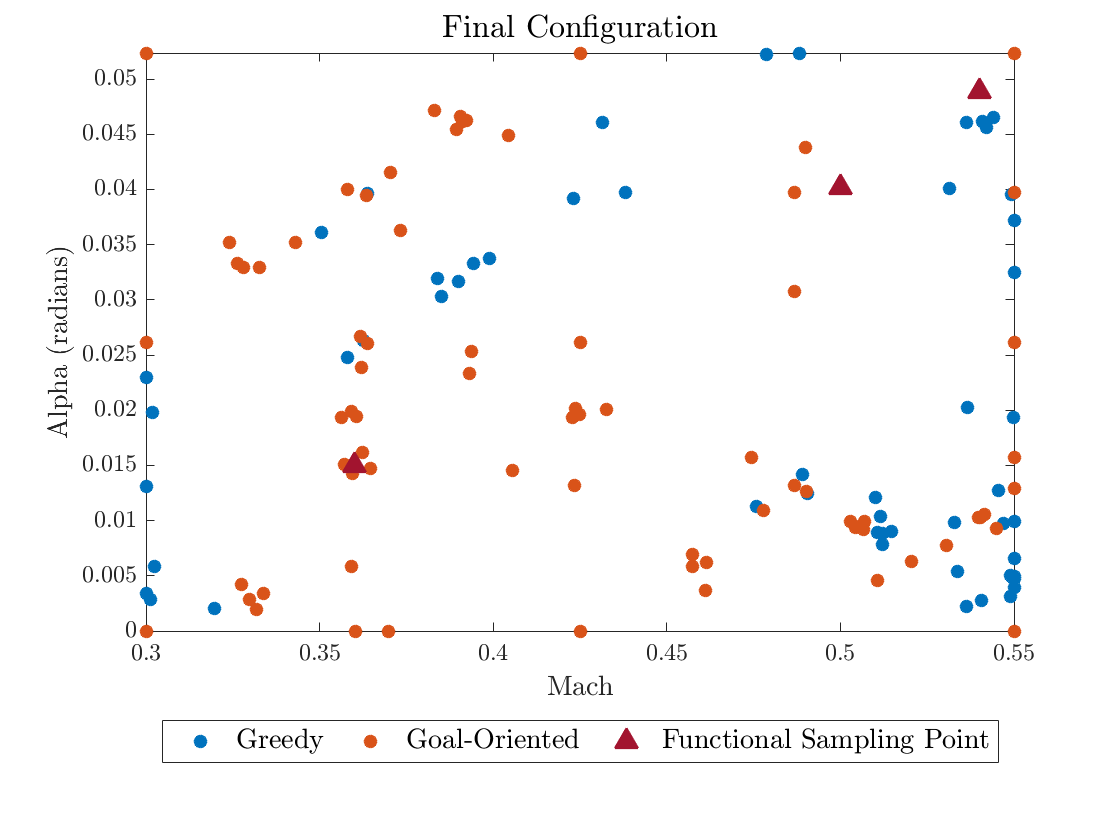}
\captionsetup{justification=centering}
\caption{Locations of the Snapshots for the Reduced Basis Built Using the Greedy Adaptive Sampling and Goal-Oriented Adaptive Sampling as well as the Points where the Functional was Sampled}
\label{fig:sampling_points}
\end{figure}

Similarly to the previous test case, the functional value (lift) at a design parameter location is used to compare the goal-oriented sampling procedure with hyperreduction and a greedy sampling procedure. Figure \ref{fig:greedy_config} displays the final configuration of the snapshot and ROM locations for the HROM constructed using the greedy sampling procedure. It is evident that altering the error criterion driving snapshot selection leads to a significantly different distribution compared to Figure \ref{fig:HROM_DWR_sub_fine}. Also given that the sampling procedure is run to the same work units but each sampling cycle is cheaper than the goal-oriented, this approaches places more total snapshots. It has a final reduced-subspace dimension of 75. The goal-oriented approach tends to cluster more snapshots in the highest Mach number region, with a large concentration between 0.5 and 0.55. In contrast, the greedy sampling approach distributes the snapshots more evenly across all Mach and alpha values. Notably, this approach places fewer snapshots in the upper-right region of the parameter space compared to the goal-oriented method. Based on these results, three parameter locations will be used to compare the construction of the reduced-order basis. The first location, Mach 0.54 and an angle of attack of 0.04885 radians, lies near the upper range of both parameters. This region was identified by the goal-oriented approach as having higher functional errors, leading to a greater number of snapshots being added compared to the greedy approach. Examining this location will help assess whether the goal-oriented method is more effective in identifying regions with higher functional error. The second location, Mach 0.5 and an angle of attack of 0.04014 radians, falls within a region where goal-oriented sampling added very few snapshots. This suggests that the functional error tolerance was already met at nearby ROM points, and additional snapshots in this area were deemed less critical for improving HROM accuracy. Tracking the lift at this location will help evaluate how well the reduced-order basis performs in regions where fewer snapshots were added. Finally, both the functional estimations from both models are sampled at Mach 0.36 and an angle of attack of 0.015 radians, a region where the greedy sampling approach concentrates many of its snapshots where the goal-oriented does not. All three of these locations as well as the final snapshots from each of the sampling approaches are plotted in Figure \ref{fig:sampling_points}. The points used to query the functional are indicated by red triangles.

\begin{figure}[H]
\centering
\includegraphics[trim = 20 10 60 10, clip, width=0.49\textwidth]{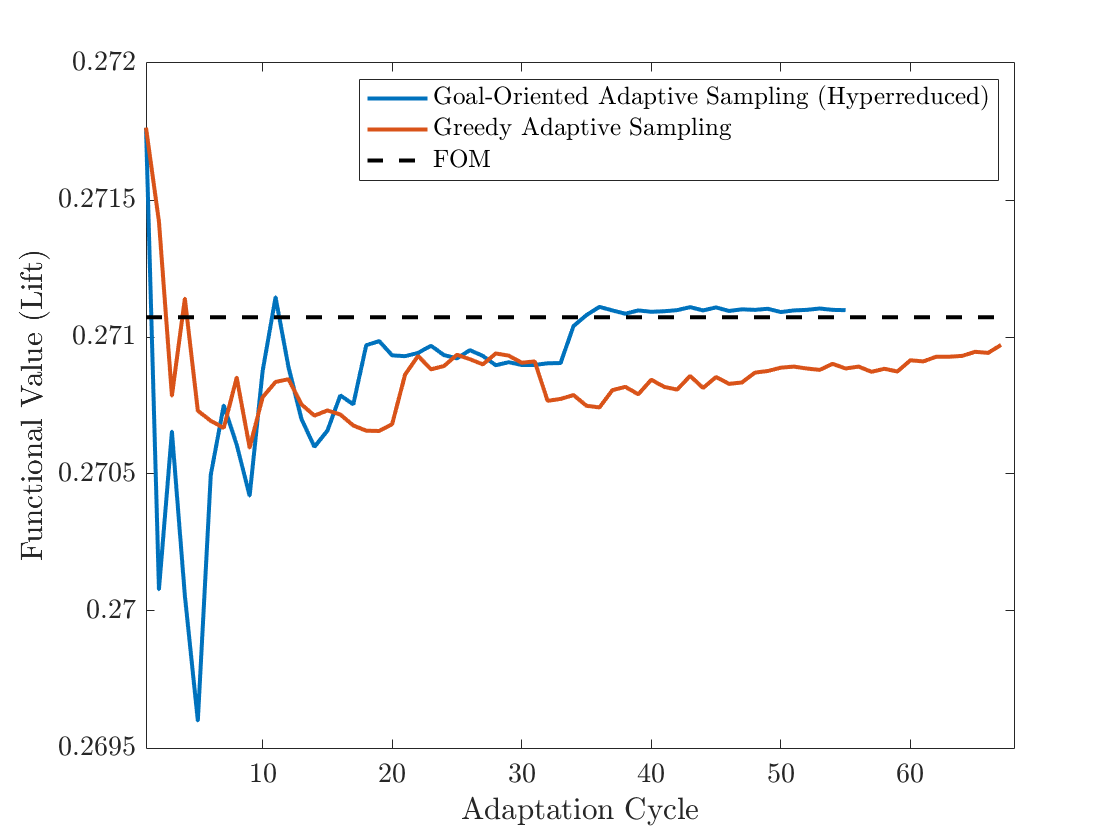}
\includegraphics[trim = 30 10 60 10, clip, width=0.49\textwidth]{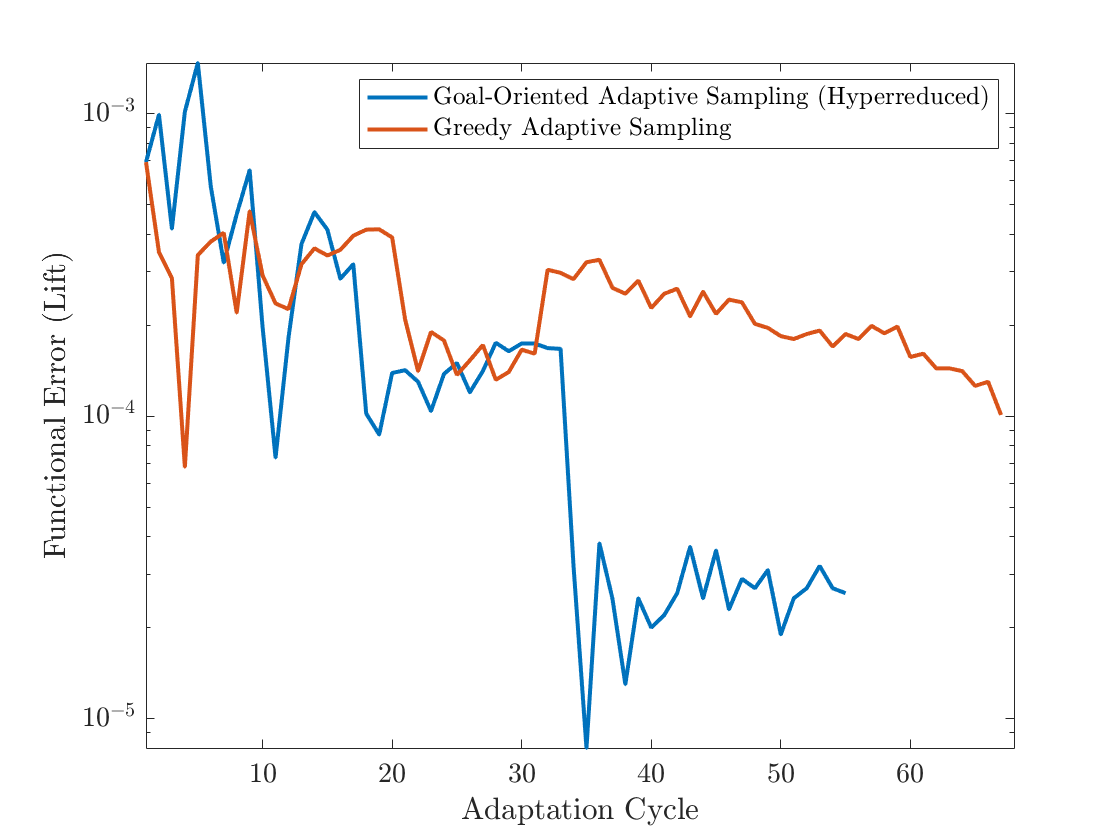}
\includegraphics[trim = 20 10 60 10, clip, width=0.49\textwidth]{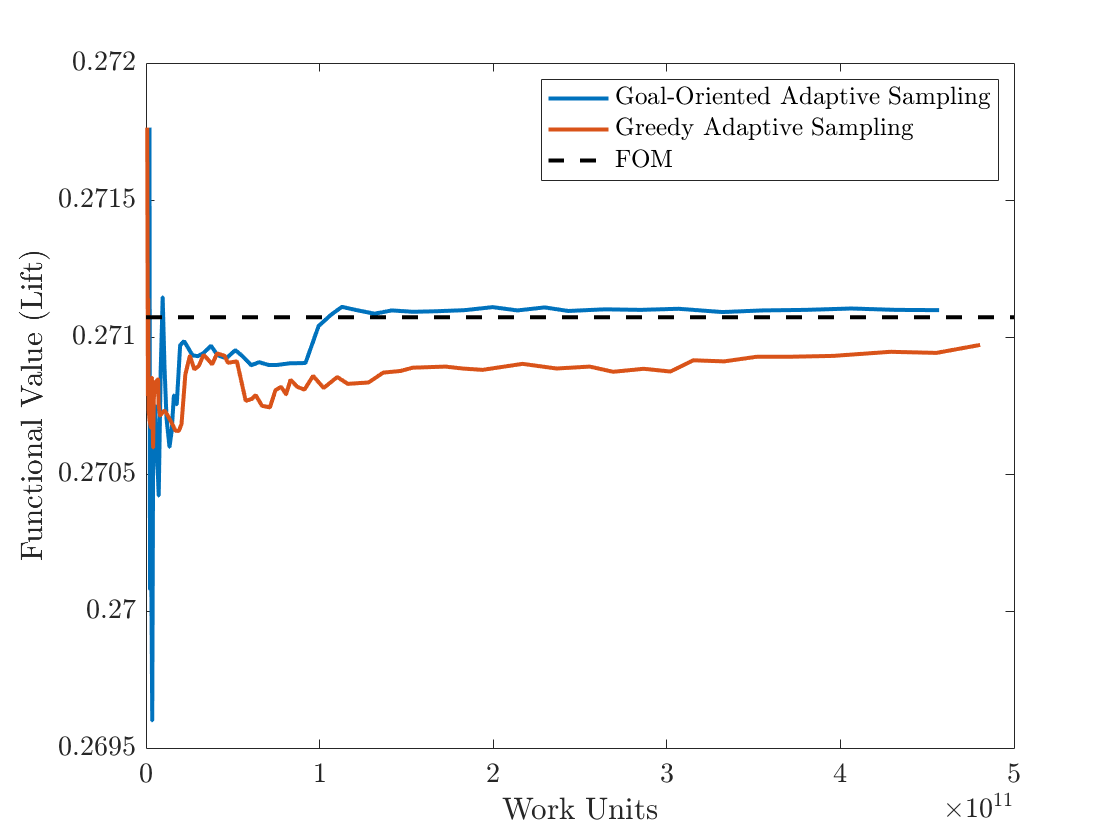}
\includegraphics[trim = 30 10 60 10, clip, width=0.49\textwidth]{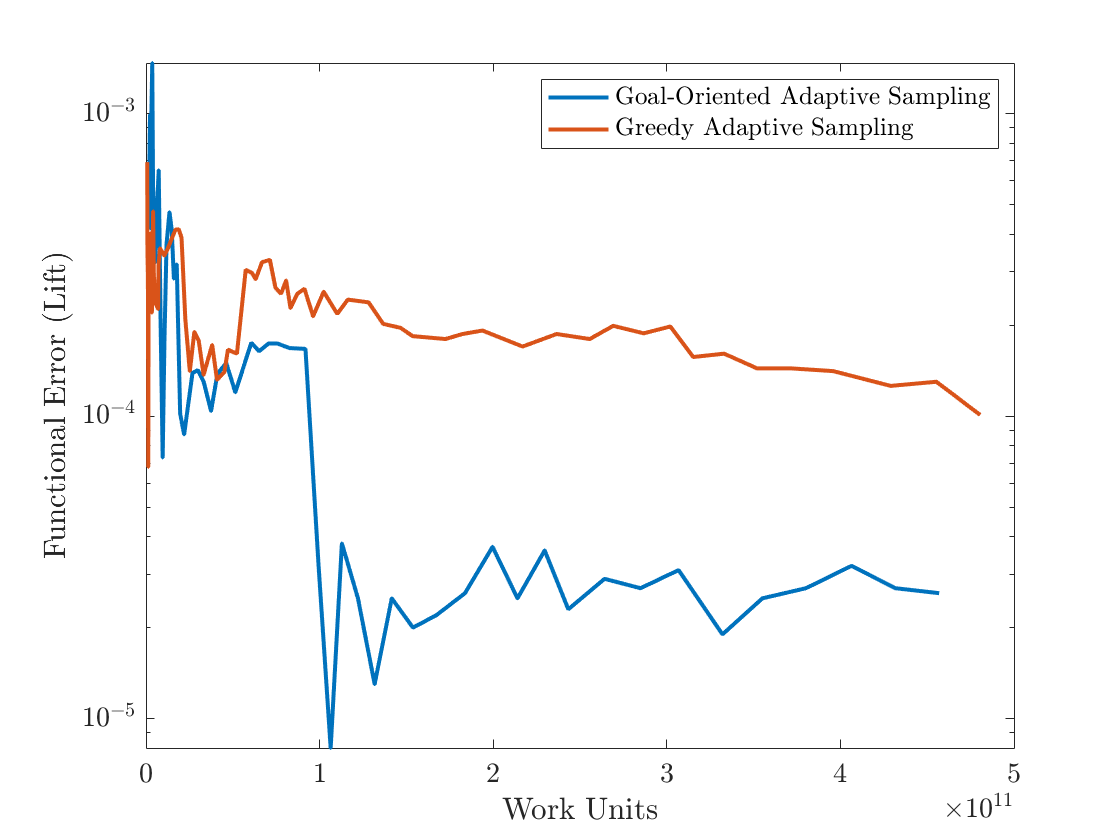}
\captionsetup{justification=centering}
\caption{Comparison of the Functional Value and Error in the Functional between Goal-Oriented and Greedy Adaptive Sampling when Plotted Against Adaptation Cycle and Works Units at Mach = 0.54 and Angle of Attack = 0.04885 radians}
\label{fig:sub_func_tracking}
\end{figure}

Figure \ref{fig:sub_func_tracking} shows the changes in the lift approximation and error in the functional from the two sampling procedures over the adaptive cycles and the work units at Mach 0.54 and an angle of attack of 2.8 degrees (0.04885 radians). The sampling procedure developed in this work converges to the FOM lift value in both fewer adaptation cycles and fewer work units compared to the greedy sampling approach. This suggests that the goal-oriented method is more effective at constructing a reduced-order model, as it leverages the DWR error indicators to guide the reduced basis construction by directly minimizing functional error. In contrast, the greedy sampling approach does not identify this region as having high errors, likely because it focuses on residual-based error rather than functional error. As a result, it fails to converge to the correct lift value at this location, even though the dimension of the reduced-order basis, $n$, is larger. While a higher-dimensional reduced-order basis should theoretically improve accuracy, its effectiveness depends on the strategic placement of snapshots, which the goal-oriented approach achieves more efficiently. Because the greedy approach is not informed by the DWR error indicators, which are indications of the expected online error in the functional, it is less likely that is meets the same error tolerance as the goal-oriented. This is shown in this sampling location.

\begin{figure}[H]
\centering
\includegraphics[trim = 30 10 60 10, clip, width=0.49\textwidth]{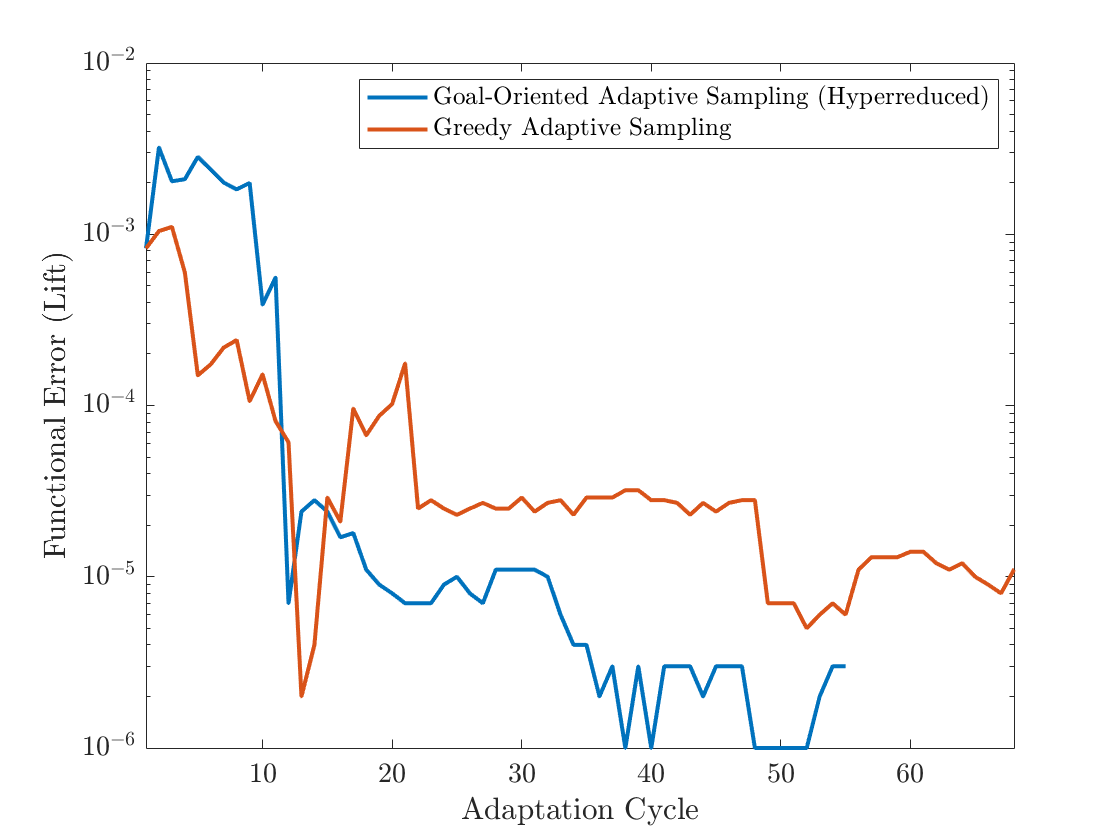}
\includegraphics[trim = 30 10 60 10, clip, width=0.49\textwidth]{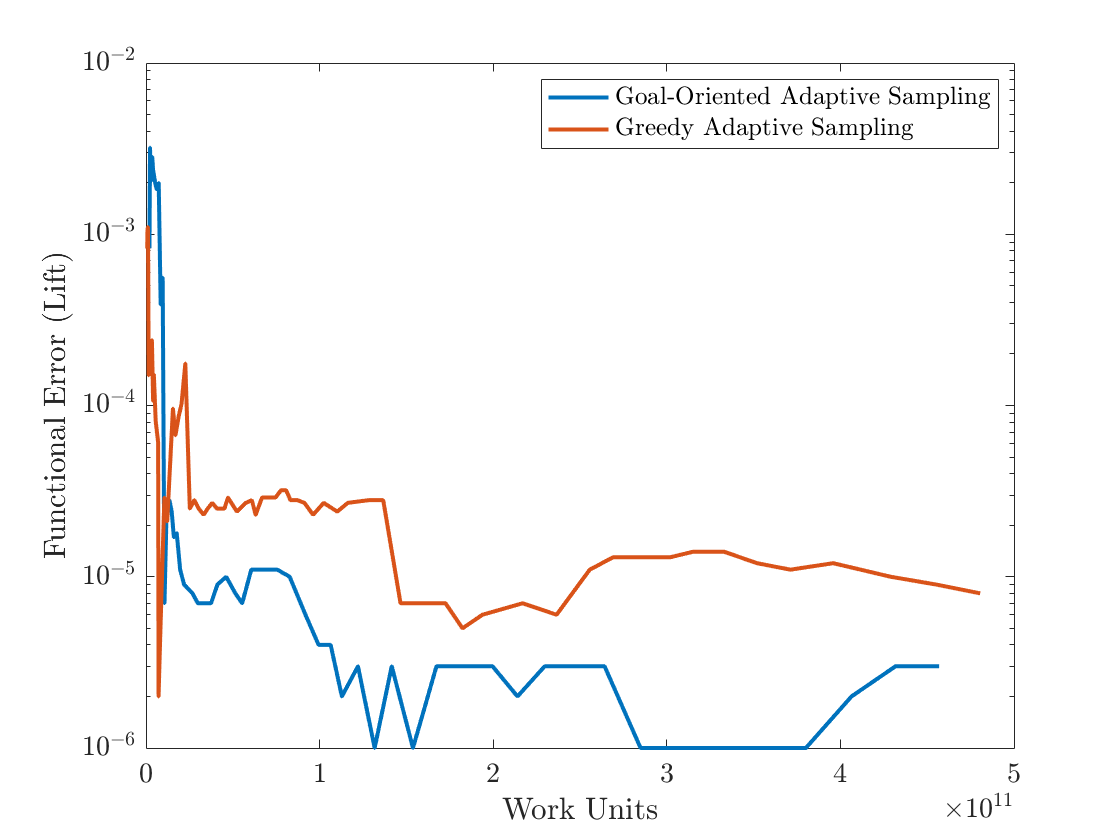}
\captionsetup{justification=centering}
\caption{Comparison of the Error in the Functional between Goal-Oriented and Greedy Adaptive Sampling when Plotted Against Adaptation Cycle and Works Units at Mach = 0.5 and Angle of Attack = 0.04014 radians}
\label{fig:sub_func_tracking_2}
\end{figure}

Figure \ref{fig:sub_func_tracking_2} presents the results at Mach 0.5 and an angle of attack of 2.3 degrees (0.04014 radians). Compared to the previous parameter location, both approaches converge more quickly to the correct lift value, indicating that the HROM can accurately predict the functional in this region without requiring many additional snapshots. Both sampling approaches can predict the lift to the requested accuracy of $1 \times 10^{-4}$, however the goal-oriented does this in fewer adaptation cycles and work units and also achieves a lower functional error at this parameter location. This suggests that the goal-oriented method is effective not only in identifying high-error regions but also in avoiding unnecessary snapshot additions where the functional error tolerance has already been met.

\begin{figure}[H]
\centering
\includegraphics[trim = 30 10 60 10, clip, width=0.49\textwidth]{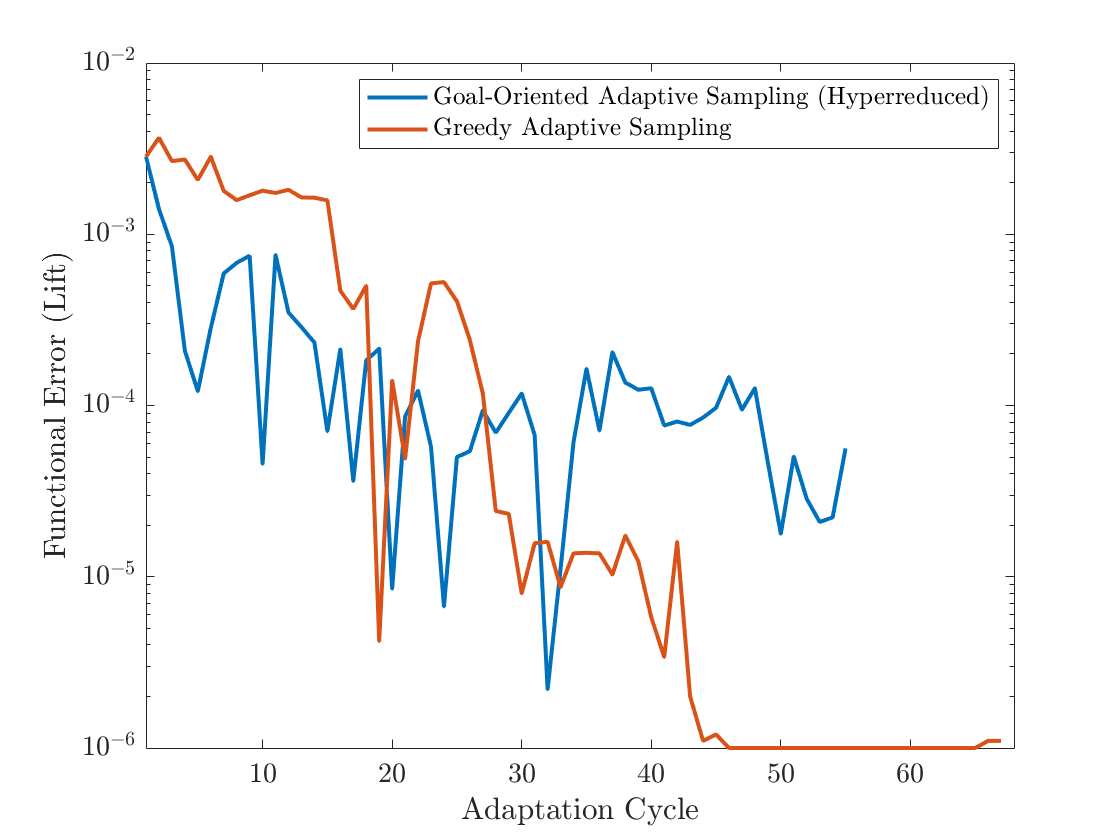}
\includegraphics[trim = 30 10 60 10, clip, width=0.49\textwidth]{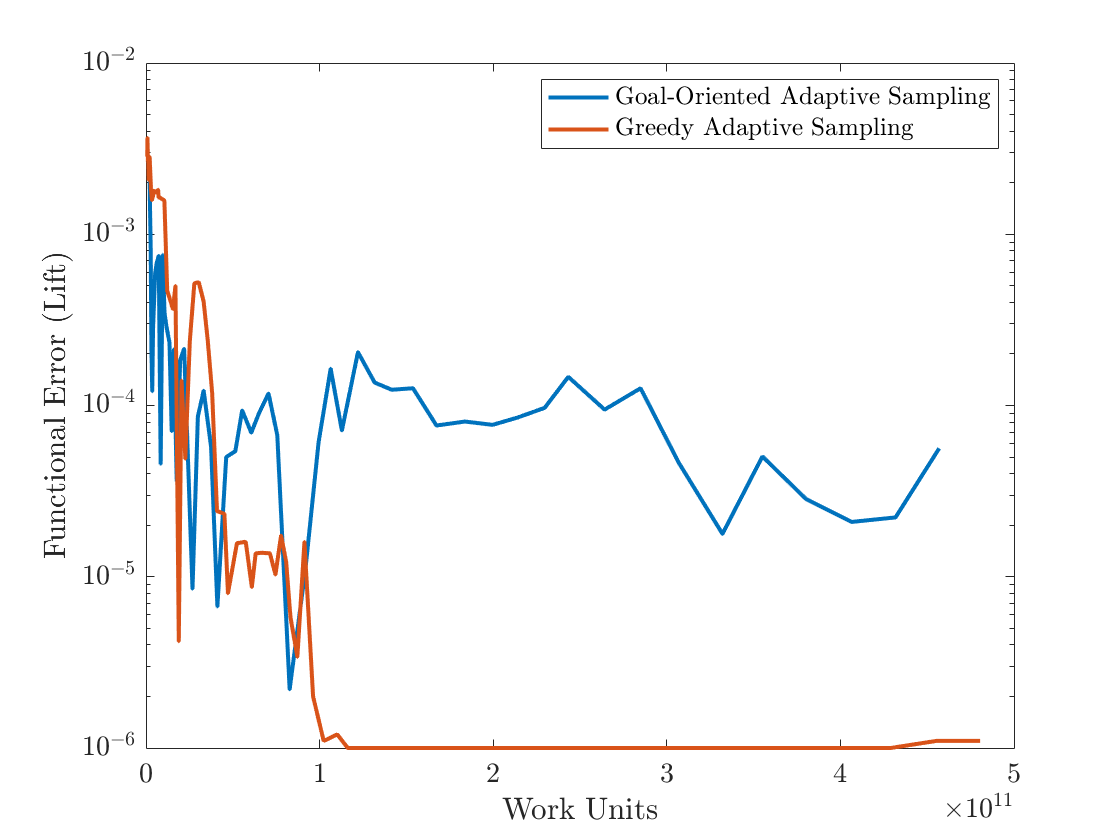}
\captionsetup{justification=centering}
\caption{Comparison of the Error in the Functional between Goal-Oriented and Greedy Adaptive Sampling when Plotted Against Adaptation Cycle and Works Units at Mach = 0.36 and Angle of Attack = 0.015 radians}
\label{fig:sub_func_tracking_3}
\end{figure}

Figure \ref{fig:sub_func_tracking_2} shows the results at the final sampling point of Mach 0.36 and an angle of attack of 0.86 degrees (0.015 radians). In contrast to the last two points, the greedy sampling approach outperforms the goal-oriented. It converges to the lift in fewer iterations and work units. However, it is not a requirement that the goal-oriented adaptive sampling approach will converge in the fewest iterations at all points in the grid. It's main purpose is to ensure that the tolerance bound set on the functional error is met at all points in the parameter space. We can see here that even without adding additional snapshots in this area, the goal-oriented approach still meets the intended bound of $1 \times 10^{-4}$ by the end of the sampling cycle. Here the greedy sampling adds more snapshots close to this location at the beginning of the sampling cycle, which is why the error in this area drops more quickly. However, as mentioned previously, we can see from the end of the goal-oriented sampling that these may be unnecessary or inefficient additions to the reduced-order basis. All three sampling points show how the goal-oriented sampling approach developed in this work strikes a balance at constructing an accurate reduced-order basis while limiting the total number of FOM snapshots computed to keep the online costs low.

%% file: 07_conclusion.tex
\section{Conclusion}\label{sec:conclusion}

This work explores the use of hyperreduction as a way to address the computational bottlenecks caused by the re-evaluation of high-dimensional quantities in a projection-based reduced-order modelling framework, specifically for parametric, highly nonlinear computational models. The energy-conserving sampling and weighting technique was selected out of the many available techniques for hyperreduction. It had been previously used in least-squares Petrov-Galerkin projection models and has been shown to be robust, accurate, and stable \cite{GRIM2021, GRIMAIAA, FARHAT2014, FARHAT2015}. 

First, it was verified that the ECSW method could successfully identify a reduced mesh for a simple 1D Burgers' Equation problem with one design parameter. As expected, an additional source of error is introduced once this set is used to approximate the residual and Jacobian in a ROM. The accuracy of the hyperreduction is impacted both by the source of the training data and the NNLS tolerance used. The Jacobian-based training data proved to be more robust and reliable, confirming the results found in \cite{TEZAUR2022}. The NNLS tolerance serves as a measure of the tradeoff between accuracy and computational savings. The lower the tolerance, the more elements are included in the reduced mesh resulting in a less efficient HROM but likely a more accurate approximation of the projected high-dimensional quantities.

Hyperreduction is not only important for achieving computational efficiency in online numerical predictions but also for accelerating the offline construction of a ROM through a greedy or adaptive sampling procedure. While the ECSW hyperreduction technique has been used in other greedy sampling frameworks \cite{TEZAUR2022}, these use residual or solution-based error indicators. In this work, we incorporate hyperreduction into a goal-oriented adaptive sampling method which uses dual-weighted residual errors that tailor a reduced-order model to a specific output of interest \cite{blaisTHESIS}.

One of the key contributions of \cite{blaisTHESIS} is the addition of a second dual-weighted residual error between a coarse and a fine reduced-order model, which makes it possible to recycle an existing reduced-order model solution computation. In this work, the error indicator is updated to measure the error between two hyperreduced reduced-order models. This allows the sampling procedure to capture the additional error introduced by the hyperreduction, and consider this when selecting new snapshot locations.

The updated sampling framework was tested on two-dimensional test cases, specifically the NACA0012 with various design parameter combinations and flow conditions. It was shown that the ECSW hyperreduction method can significantly reduce the number of elements required to re-evaluate the high-dimensional quantities and that with the updated DWR errors, the framework can still construct a ROM which provides a certain level of confidence in the expected error in the output functional. In all test cases, the average online error is within the vicinity of the tolerance specified when the model was trained, suggesting that the DWR error indicators accurately predict the online performance of the ROMs and HROMs. There are only small regions in the parameter space where the tolerance bound set are not met online. Furthermore, the work units required to assemble the model are significantly reduced through hyperreduction in the reduced-order solutions and in the error indicator. This suggests savings both in the offline and online stages of the reduced-order modeling procedure. Finally, the adaptive sampling procedure is compared to a greedy sampling approach which instead uses the residual to guide new snapshot location selection. These results showed that in terms of work units, the offline adaptation cycles of the greedy sampling are less expensive due to the lack of DWR calculations and adjoint solutions. Since the greedy sampling procedure could not be terminated at the same functional error tolerance as it was no longer guided by the DWR error measures, the procedure was instead concluded at a similar number of total work units so we could compare two ROMs that were equally expensive to construct offline. This meant that the greedy sampling produced a larger reduced-order subspace since more snapshots were added within those work units. If snapshot placement did not impact the ROM quality, the greedy sampling procedure should hypothetically be more accurate than the goal-oriented, as it has a larger reduced-order basis. However, it was found that the goal-oriented sampling tends to outperform the greedy sampling in terms of the speed and accuracy of the functional (lift) estimations at various points in the parameter space. Even when greedy sampling converges more quickly at a certain location, the goal-oriented will still reach the desired error tolerance at this location by the end of the procedure while greedy sampling may not at all locations. These results highlight the importance of intelligent snapshot selection in basis construction. Not only is the basis produced by the goal-oriented sampling more accurate, but it is also smaller, which means that its online calculations are less expensive. This shows the strength of using the DWR error indicators to guide the construction of both an efficient and accurate reduced-order basis.

\subsection{Future Work}

There are multiple avenues for future work in this research. The scope of the research can be extended to unsteady CFD problems, 3D problems, and problems with more complex parameter spaces. The cost and necessity of re-training the hyperreduction at every sampling cycle could also be reviewed; in particular, possible additional error measures, like those in \cite{TEZAUR2022}, could be used to determine whether the hyperreduction must be recomputed at a particular iteration. Alternatives to the non-negative least squares problem for finding the reduced set of important mesh elements could also be considered. Finally, full implementation of the hyperreduction approach, which only computes the residual and Jacobian on the elements in the reduced mesh should be tested in order to study the true computational cost and storage savings.

%% file: appendix.tex
\appendix

\section{ECSW Hyperreduction Verification}\label{sec:verf}

Before integrating hyperreduction into the adaptive sampling procedure, we will verify that the ECSW method has been implemented correctly. First, we will use the method presented in \cite{blaisTHESIS} to generate a POD basis. With the FOM snapshots used to build the basis, a reduced mesh will be determined for multiple NNLS tolerance values, employing both residual- and Jacobian-based training data. Next, we will test a HROM that utilizes the same POD basis, the identified reduced set of elements, and the corresponding weights, at various parameter locations using Algorithm \ref{alg:online}. The test case used is the one-dimensional Burgers' equation with one design parameter. The problem is described by the following differential equation \cite{rewienski2003}:
\begin{equation}
\begin{split}
            & \mathbf{R(w)} = \frac{\partial \mathbf{w}(x,t)}{\partial t} + 0.5\frac{\partial \mathbf{w}^2(x,t)}{\partial x} - \mathbf{S}(\mu),\\
            &\text{with initial and boundary conditions:} \\
            &\mathbf{w}(x,0) = 1, \forall x \in [0,100] ,\\
            &\mathbf{w}(0, t) = 1, \forall t > 0 .\\
\end{split}
\end{equation}

\begin{figure}[H]
\centering
 \includegraphics[trim=30 15 50 10,clip,width=0.7\textwidth]{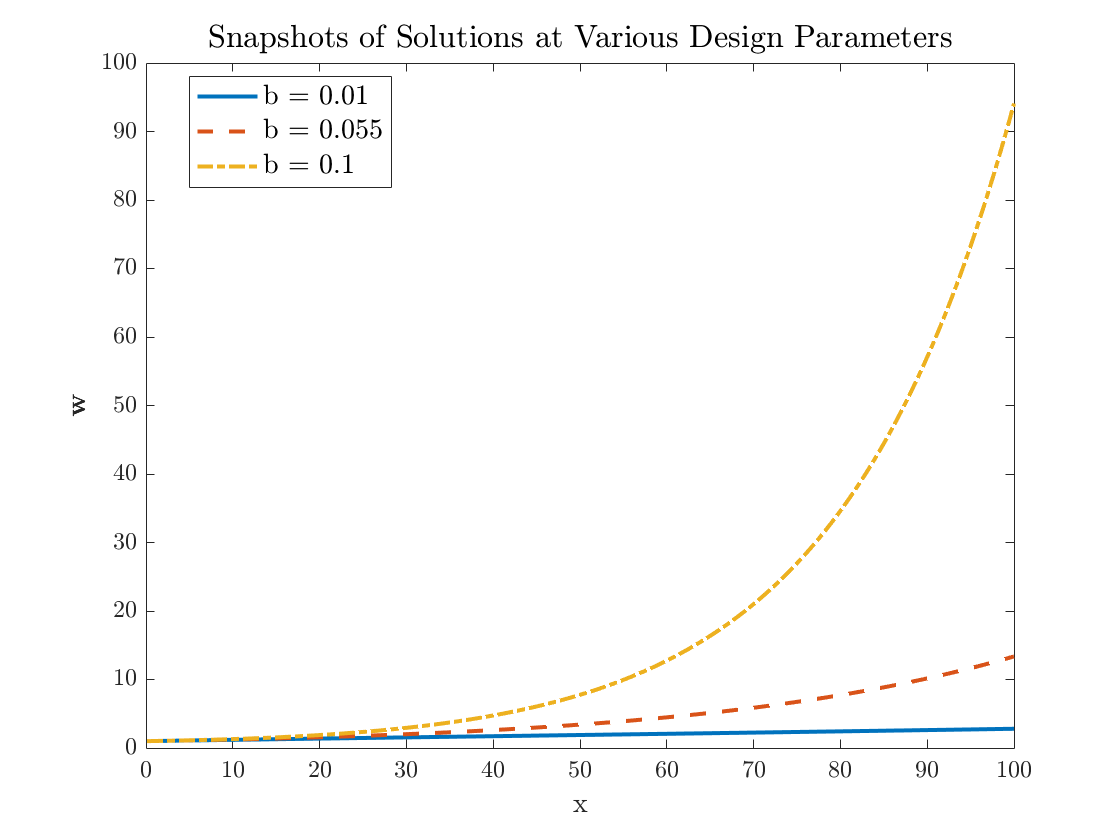}
\captionsetup{justification=centering}
\caption{Solutions of 1D Burgers' Equation with Exponential Source Term from Upwind Scheme}
\label{fig:simple_snaps}
\end{figure}

Only steady-state solutions are considered in this work. The functional to optimize is the integral of the steady-state solution $\mathbf{w}$ over the domain:
\begin{equation}\label{eq:func}
    \mathcal{J}(\mathbf{w}) = \int_{0}^{100} \mathbf{w} \;dx .
\end{equation}
The spatial domain is discretized using 1024 nodes (i.e. the FOM dimension $N$ is 1024), the FOM solution is approximated using the upwind scheme:
\begin{equation}
    \frac{w_j^{n+1} - w_j^{n}}{\Delta t} = \frac{F(w_{j+1}^{n}) - F(w_j^{n})}{\Delta x} + S(\mu, w_j^{n}),
 \end{equation}
where $F$ is the flux function $F(u) = \frac{u^2}{2}$, and $S$ is the source term. An exponential source term will be used: $S(b, x) = e^{b x}$, where $b$ is the one design parameter that can be varied to values between 0.01 and 0.1. Example solutions from the upwind scheme can be seen in Figure \ref{fig:simple_snaps} at three evenly spaced parameter locations.

The adaptive sampling procedure was initialized with three solution snapshots, shown in Figure \ref{fig:simple_snaps}, which are uniformly distributed over the parameter space. After running to an error tolerance of $1E-4$, the procedure solved the FOM at 4 additional parameter locations, resulting in a ROM dimension $n$ of 7. All 7 snapshots were then used as training data for both the residual and Jacobian-based ECSW approaches. Three values of $\epsilon$ in the NNLS problem were tested to evaluate the trade-off between computational/storage savings and accuracy in the ECSW method. For the residual-based approach, $\epsilon$ values of $1E-4$, $1E-6$ and $1E-8$ were tested. The original intention was to test the same values for the Jacobian-based approach; however, the algorithm used to solve the NNLS problem was only able to achieve a minimum tolerance of $1E-7$. This limitation highlights a potential shortcoming of the NNLS approach, as it does not guarantee convergence to a solution for all tolerance values. Ideally,  as the tolerance approaches very low values, the solution vector should converge to a vector of ones. Even in the residual-based approach, there are practical limits on acceptable $\epsilon$ values. These limitations may stem from the conditioning of the $\mathbf{C}$ matrix, which could be affected by factors such as linear dependence among rows in the training data or subtractive cancellation errors. Furthermore, since one of the main objectives of the goal-oriented adaptive sampling procedure is to reduce the number of FOM solutions computed, the amount of training data available is inherently limited, which could further impact the performance of the ECSW method.

\begin{figure}[H]
\centering
 \includegraphics[trim=50 10 90 20,clip,width=0.48\textwidth]{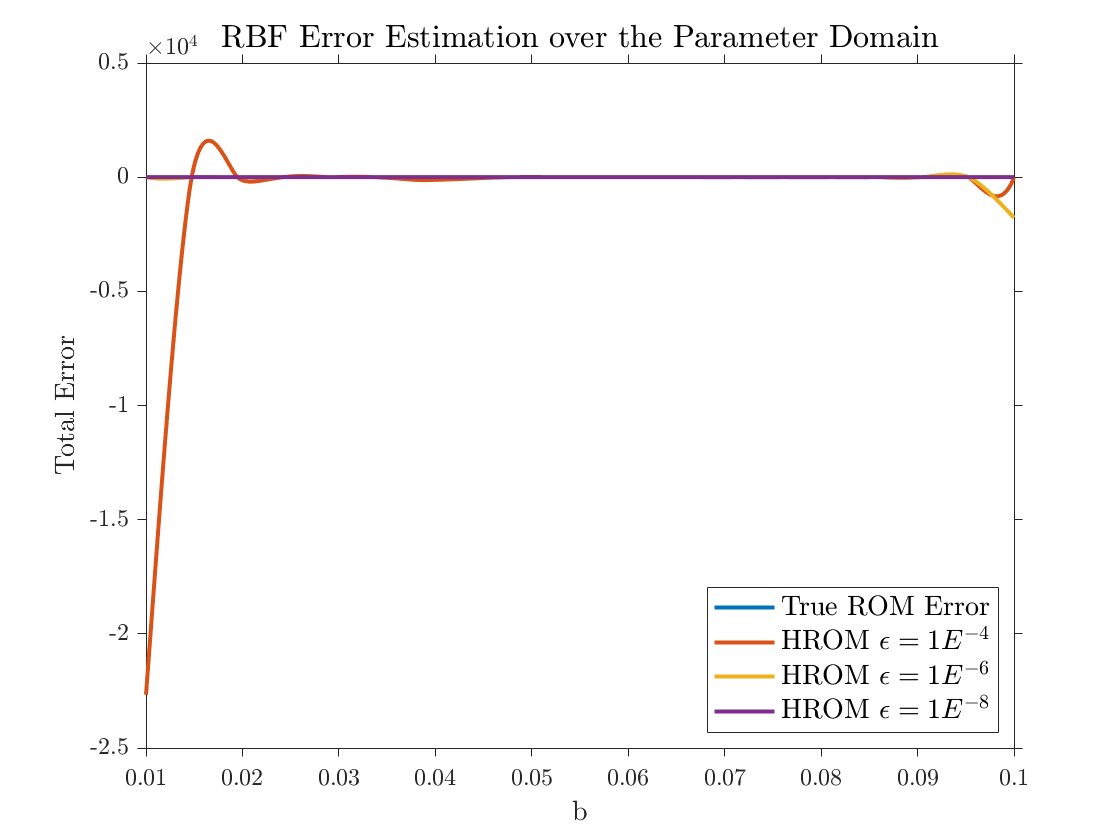}
  \includegraphics[trim=50 10 90 20,clip,width=0.48\textwidth]{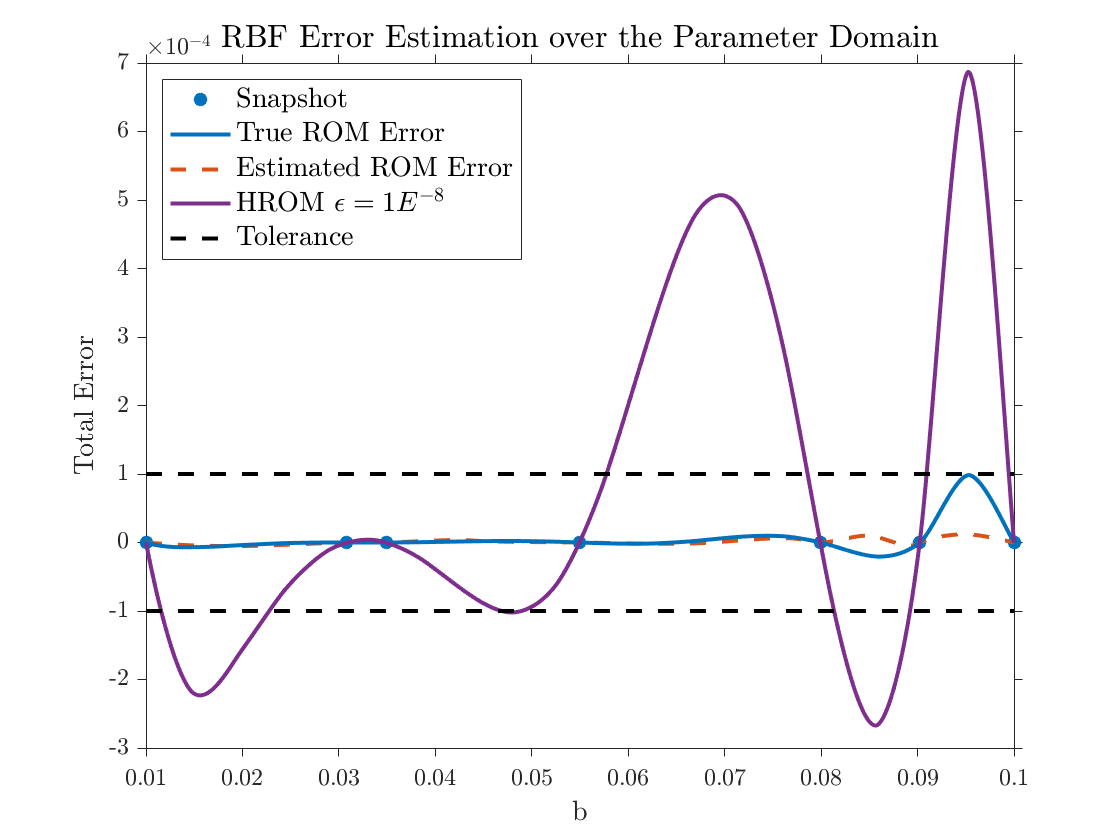}
\captionsetup{justification=centering}
\caption{Error in the Functional $\mathcal{J}$ over the Parameter Domain for the Final ROM and Three Residual-based ECSW HROM}
\label{fig:err_res_burg}
\end{figure}

\begin{figure}[H]
\centering
 \includegraphics[trim=40 10 90 20,clip,clip,width=0.48\textwidth]{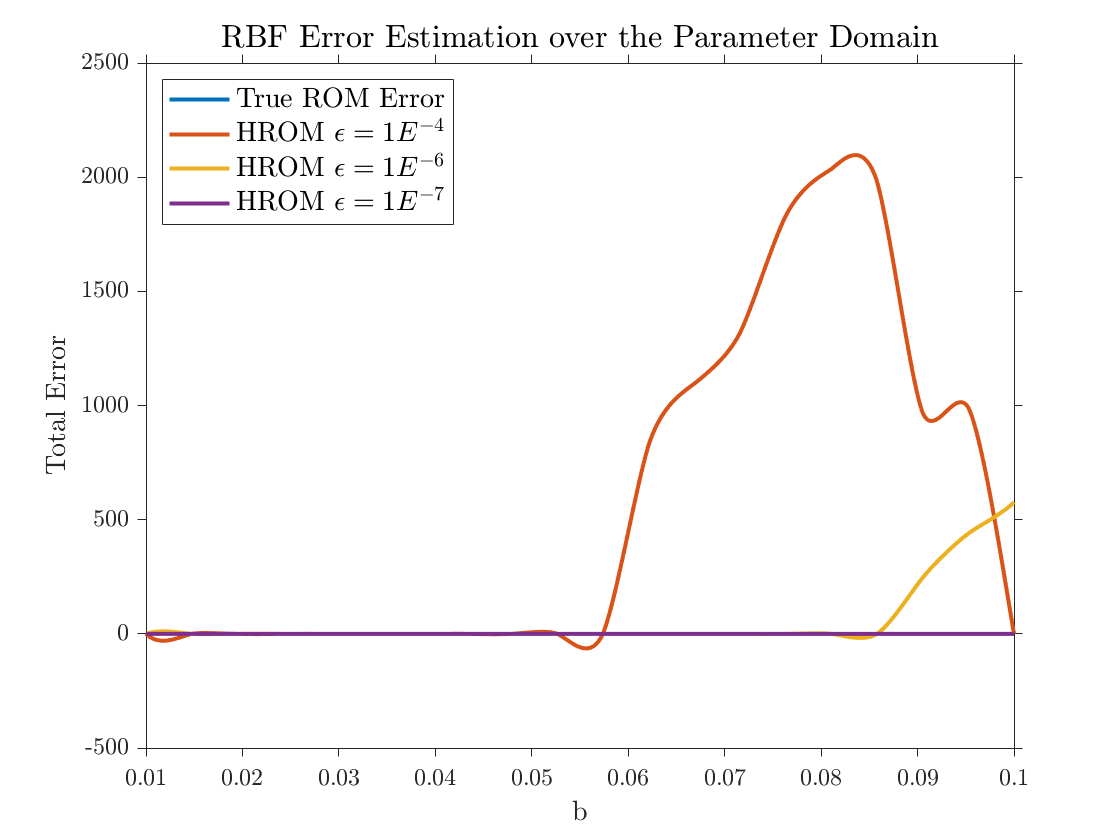}
  \includegraphics[trim=50 10 90 20,clip,width=0.48\textwidth]{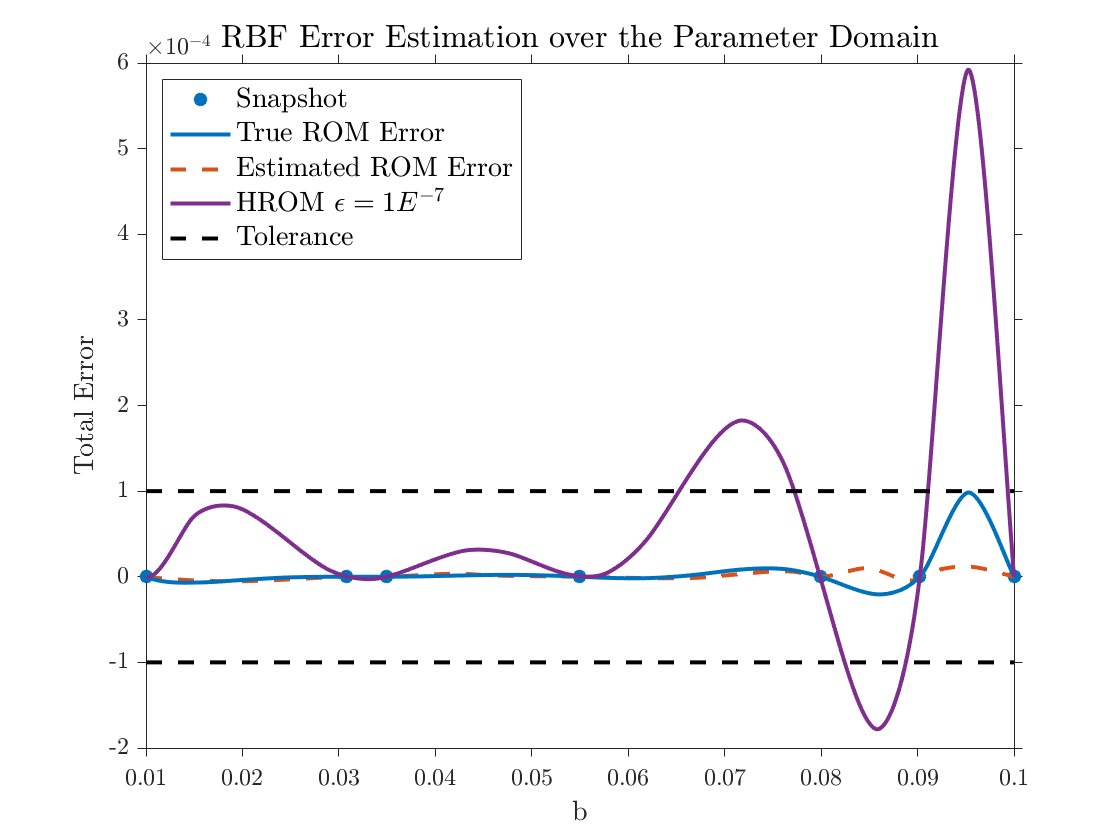}
\captionsetup{justification=centering}
\caption{Error in the Functional $\mathcal{J}  $ over the Parameter Domain for the Final ROM and Three Jacobian-based ECSW HROM}
\label{fig:err_jac_burg}
\end{figure}

Using the reduced mesh identified from each combination of training data type and NNLS tolerance value, six different hyperreduced ROMs (HROMs) were constructed. The error between the FOM and the HROMs was sampled at 20 evenly spaced points, and radial basis function (RBF) interpolation was employed to estimate the expected online or ``true" error distribution across the parameter domain for each model. The same process was applied to the ROM without hyperreduction from the sampling procedure. Figures \ref{fig:err_res_burg} and \ref{fig:err_jac_burg} show the functional error results for the HROMs built using residual and Jacobian-based training data, respectively. In both figures, the left image displays the ``true" error estimate for the ROM and three HROMs with different $\epsilon$ values. The right image shows only the HROM with the lowest NNLS tolerance value: $\epsilon = 1E-8$ for the residual-based ECSW and $\epsilon = 1E-7$ for the Jacobian-based ECSW. These plots also include the estimated error distribution from the final iteration of the sampling procedure, the ``true" error for the ROM, and the functional error tolerance bound. The results reveal that for both approaches, higher NNLS tolerance values lead to significant spikes in the error distribution, indicating that the reduced mesh cannot accurately approximate the residual and Jacobian online. When the tolerance is reduced to its lowest value, the performance becomes comparable to the ROM without hyperreduction. Among the tested models, the Jacobian-based HROM with $\epsilon = 1E-7$ performs best, though it only marginally outperforms the residual-based HROM with $\epsilon = 1E-8$. At its worst, the Jacobian-based HROM violates the original tolerance bound by no more than an order of magnitude at the ROM points.

\begin{figure}[H]
\centering
 \includegraphics[trim = 40 20 60 20, clip, width=0.7\textwidth]{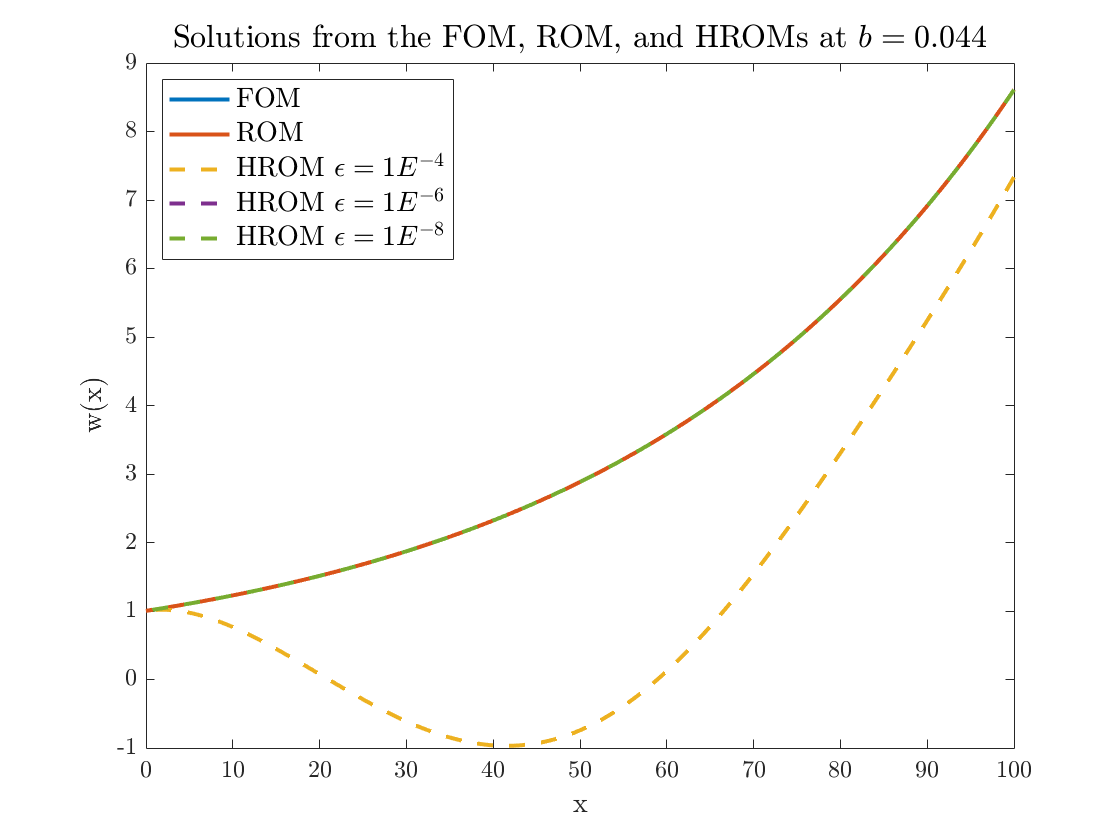}
\captionsetup{justification=centering}
\caption{Solutions of 1D Burgers' Equation from the FOM and Projection-Based ROM with and without Residual-based Hyperreduction}
\label{fig:res_sol}
\end{figure}

\begin{figure}[H]
\centering
 \includegraphics[trim = 40 20 60 20, clip, width=0.7\textwidth]{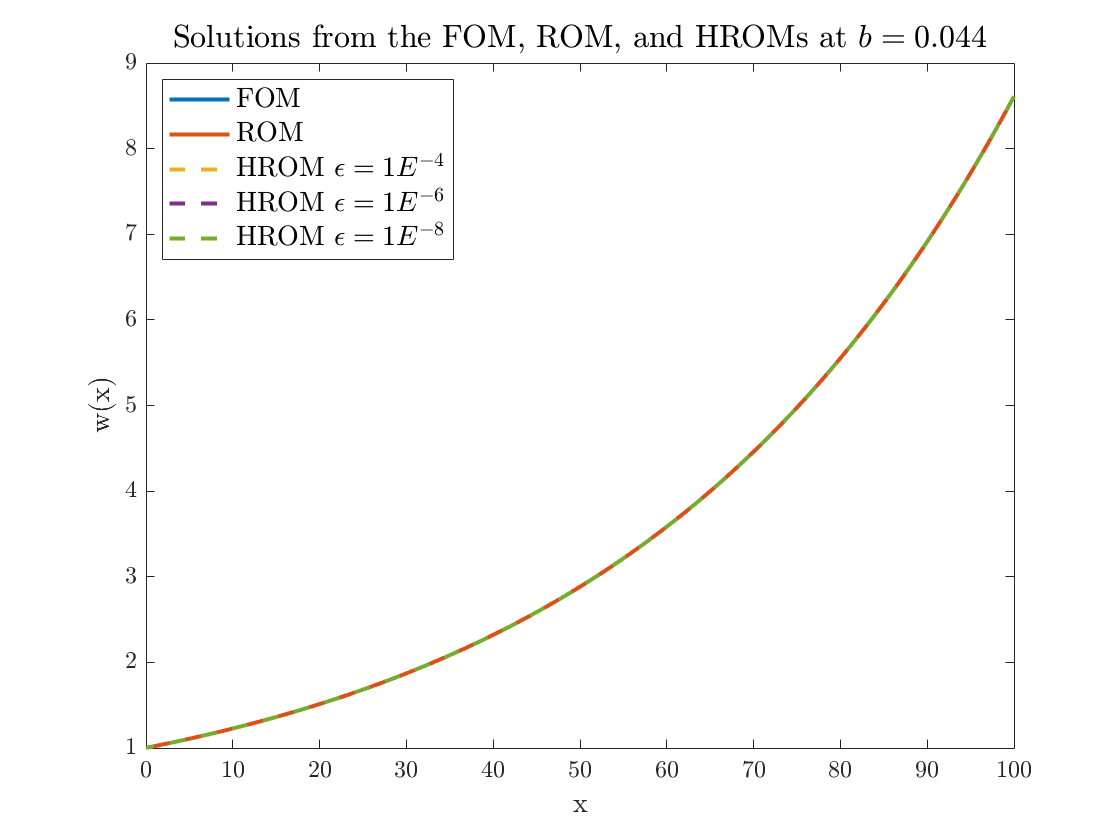}
\captionsetup{justification=centering}
\caption{Solutions of 1D Burgers' Equation from the FOM and Projection-Based ROM with and without Jacobian-based Hyperreduction}
\label{fig:jac_sol}
\end{figure}

The solution vectors from the FOM, ROM, and HROMs at a parameter location of $b = 0.044$ are plotted together, Figure \ref{fig:res_sol} displays the residual-based HROMs and Figure \ref{fig:jac_sol} displays the Jacobian-based HROMs.  With the exception of the residual-based HROM with $\epsilon = 1E-4$, the HROM and ROM solutions are nearly identical to the FOM solution, showing minimal error and appearing visually indistinguishable.
    
The solution and functional error results at this parameter location, along with the average ROM point functional error for each model, are summarized in Table \ref{tab:burg_results}. Hyperreduction has significantly decreased the number of elements required to compute the residual and Jacobian at each iteration, reducing the original 1024 elements to as few as 11 in some cases. The Jacobian-based HROMs outperform the residual-based HROMs in terms of errors introduced in both the solution vector $\mathbf{w}$ and the functional $\mathcal{J}$, as seen in the third and fourth column of Table \ref{tab:burg_results}. However, even the best-performing HROMs introduce additional errors relative to the original ROM without hyperreduction, as expected due to the added layer of approximation.

\begin{table}[H]
\caption{1D Burgers' ROM Results at b = 0.044 and Average Absolute ROM Point Error over the Parameter Domain} \label{tab:burg_results}
\centering
\begin{tabular}{|p{2.8cm}||p{0.8cm}|p{2.6cm}|p{2.6cm}|p{2.6cm}|}
        \hline
         ROM/HROM & $\lVert \boldsymbol{\xi} \rVert_0$  &$\lVert \mathbf{w}_{*} - \mathbf{w}_{FOM} \rVert_2$ & $\lVert \mathcal{J}_{*} - \mathcal{J}_{FOM} \rVert_2$ & ROM Error\\
        \hline \hline
        ROM & - & $2.0431E-8$     & $-3.7066E-6$ & $8.3270E-6$\\ \hline
        Res. $\epsilon = 1E-4$ & 18 & $6.061E-1$    & $-2.224E2$ & $1.1433E3$  \\ \hline
        Res. $\epsilon = 1E-6$ & 35 & $1.2079E-6$     & $4.5142E-4$ & $8.8868E1$ \\ \hline
        Res. $\epsilon = 1E-8$ & 46 & $2.1958E-7$     & $-8.0496E-5$ & $1.7082E-4$ \\ \hline
        Jac. $\epsilon = 1E-4$ & 11 & $3.8299E-6$     & $1.4331E-3$ & $5.5215E2$ \\ \hline
        Jac. $\epsilon = 1E-6$ & 25 & $1.1953E-6$     & $-4.45416-4$ & $6.3034E1$ \\ \hline
        Jac. $\epsilon = 1E-7$ & 28 & $8.8126E-8$     & $3.1820E-5$& $8.0756E-5$ \\ \hline
    \end{tabular}
\end{table}

\begin{figure}[H]
\centering
\includegraphics[trim = 50 20 70 20, clip, width=0.49\textwidth]{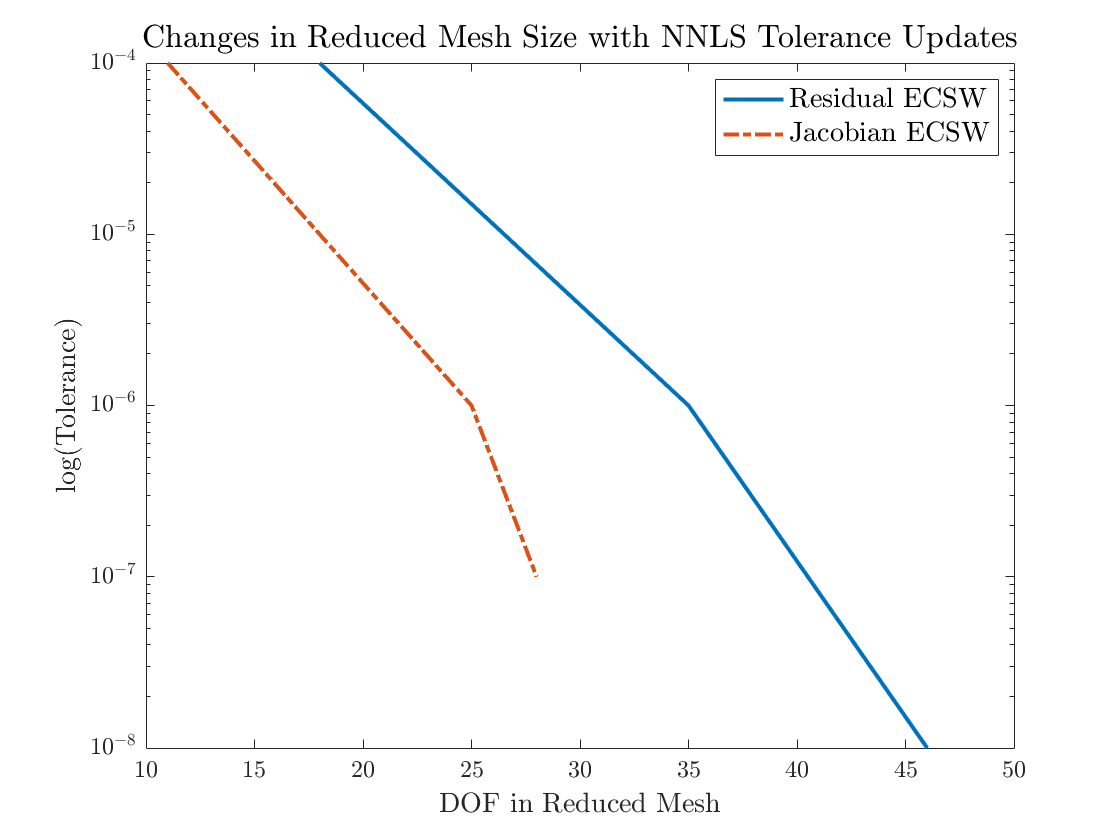}
\includegraphics[trim = 50 20 70 20, clip, width=0.49\textwidth]{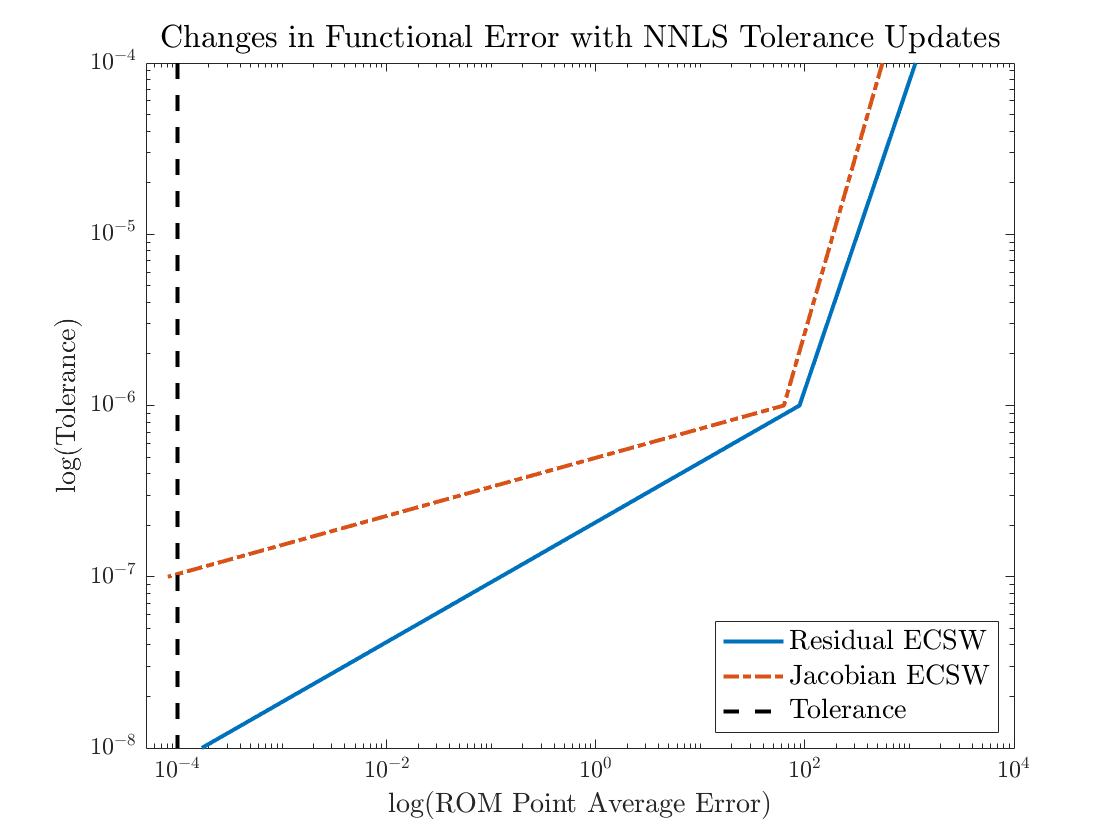}
\captionsetup{justification=centering}
\caption{Impact of NNLS Tolerance on the Size of the Reduced Mesh and the Accuracy of the HROM for the Residual and Jacobian-based Training Approaches}
\label{fig:comp_burg}
\end{figure}

The size of the reduced mesh and the average ROM point error for the HROMs are plotted against the NNLS tolerance in Figure \ref{fig:comp_burg}. From the plot on the left, it is evident that the Jacobian-based training data leads to a smaller reduced mesh for equivalent NNLS tolerance values. This suggests that the Jacobian ECSW approach more efficiently approximates the projected FOM quantities online, requiring only 28 of the total 1024 elements when $\epsilon = 1E-7$ to maintain the solution and functional error at $b = 0.044$, as well as the average ROM point error, within an order of magnitude of the ROM without hyperreduction. From the plot on the right, Jacobian-based ECSW exhibits lower average ROM point errors for the same or even higher NNLS tolerance values, indicating that it is both more accurate and computationally more efficient than residual-based ECSW.

\section{Work Units Derivations}\label{sec:appendix}

\subsection{ROM Solution Non-linear Iteration Work Units}
Equation \ref{eq:ROM_w_non_lin} is an evaluation of the work units required to solve:
\begin{equation}
\begin{split}
& \left[\mathbf{V}^T \frac{\partial \mathbf{R}}{\partial \mathbf{w}}^{(k)T} \frac{\partial \mathbf{R}}{\partial \mathbf{w}}^{(k)} \mathbf{V} \right] \mathbf{p}^{(k)} = - \mathbf{V}^T \frac{\partial \mathbf{R}}{\partial \mathbf{w}}^{(k) T} \mathbf{R}^{(k)}.
\end{split}
\end{equation}

We assume it requires one work unit to evaluate every entry in the residual and Jacobian. Since $\mathbf{R}^{(k)} \in \mathbb{R}^N$ and $\frac{\partial \mathbf{R}}{\partial \mathbf{w}}^{(k)} \in \mathbb{R}^{N \times N}$, they require $N$ and $N^2$ work units, respectively. To evaluate the LHS of the above equation, three matrix multiplications must be performed. By simply counting the number of operations required, the FLOPs necessary to perform matrix multiplication between a matrix of dimension $n \times p$ and another with dimension $p \times m$ is $nm(2p-1)$. 

Therefore, the total number of FLOPs to assemble the LHS (which are assumed to be one work unit) is:
\begin{equation}\label{eq:ROM_LHS}
\begin{split}
    W_{\text{LHS}} = & (\text{FLOPs to multiply } \mathbf{V}  \text{ and } \frac{\partial \mathbf{R}}{\partial \mathbf{w}}^{(k)}) + \\
    & (\text{FLOPs to multiply the matrix resulting from the previous step and } \frac{\partial \mathbf{R}}{\partial \mathbf{w}}^{(k)T}) + \\
    & (\text{FLOPs to multiply the matrix resulting from the previous step and } \mathbf{V}^T) \\
    = & (\text{FLOPs to multiply matrices of dimensions } (N \times N) \text{ and } (N \times n_i)) + \\
    & (\text{FLOPs to multiply matrices of dimensions } (N \times N) \text{ and } (N \times n_i)) + \\
    & (\text{FLOPs to multiply matrices of dimensions } (n_i \times N) \text{ and } (N \times n_i)) \\
    = & [Nn_i(2N-1)] + [Nn_i(2N-1)] + [n_i^2(2N-1)] \\
    = & (2Nn_i + n_i^2)(2N-1).
\end{split}
\end{equation}
Similarly, to assemble the RHS:
\begin{equation}\label{eq:ROM_RHS}
\begin{split}
    W_{\text{RHS}} = & (\text{FLOPs to multiply } \mathbf{R}^{(k)}  \text{ and } \frac{\partial \mathbf{R}}{\partial \mathbf{w}}^{(k)}) + \\
    & (\text{FLOPs to multiply the matrix resulting from the previous step and } \mathbf{V}^T) \\
    = & (\text{FLOPs to multiply matrices of dimensions } (N \times N) \text{ and } (N \times 1)) + \\
    & (\text{FLOPs to multiply matrices of dimensions } (n_i \times N) \text{ and } (N \times 1)) \\
    = & [N(2N-1)] + [n_i(2N-1)] \\
    = & (N + n_i)(2N-1).
\end{split}
\end{equation}
Therefore, the total cost to assemble the linear system is:
\begin{equation}\label{eq:ROM_assem}
\begin{split}
    W_{\text{assemble}} = & (2Nn_i + n_i^2)(2N-1)+ (N + n_i)(2N-1) \\
    = & (2Nn_i + n_i^2 + N + n_i)(2N-1).
\end{split}
\end{equation}

The work units, as shown in section \ref{sec:adapt}, for one non-linear iteration for the ROM solution is:
\begin{equation}
\begin{split}
    W_{\text{non lin}} = & (\text{units to evaluate the residual}) + \\
    & (\text{units to evaluate the Jacobian}) + \\
    &   W_{\text{assemble}} + \\
    & (\text{units to solve for solution update}) \\ 
    =  & (N) + (N^2) + [(2Nn_i + n_i^2 + N + n_i)(2N-1)] + (n_i^3).
\end{split}
\end{equation}

\subsection{HROM Solution Non-linear Iteration Work Units}
Once the ECSW hyperreduction is introduced, the cost to assemble and solve the Newton iterations changes. The following equation is the one to be evaluated:
\begin{equation}
\begin{split}\label{eq:HROM_non}
\left[\tilde{\mathbf{W}}^{(k)^{T}}\tilde{\mathbf{W}}^{(k)}\right]\textbf{p}^{(k)} = -\tilde{\mathbf{R}}^{(k)},
\end{split}
\end{equation}
where $\tilde{\mathbf{W}}^{(k)}$ is the hyperreduced test basis and $\tilde{\mathbf{R}}^{(k)}$ is the hyperreduced residual. Recall, the hyperreduced test basis can be written as:
\begin{equation}
\begin{split}
\tilde{\mathbf{W}}^{(k)} =  \bar{\mathbf{J}}^{(k)}\mathbf{V} = (\sum_{e \in \tilde{\mathcal{E}}} \xi_e \mathbf{L}_e^{T} \mathbf{J}_e^{(k)} \mathbf{L}_{e^+} )\mathbf{V}.
\end{split}
\end{equation}
It can be shown that for a sparse matrix, the number of FLOPs for matrix multiplication with a dense matrix (i.e. $S*D$ where $S$ is sparse) is two times the number of non-zeros in $S$, times the number of columns in $D$ \cite{ryan2015}. Again,assuming the cost of evaluating one entry in the matrix is one work unit, the cost to assemble the elemental Jacobian contribution $\mathbf{J}_e^{(k)}$ is $d_ed_{e^+}$. This must be then multiplied by the number of elements in the reduced mesh $n_{e_{i}}$. Assuming we can find the elemental Jacobian contribution in the global dimension (i.e. ignoring the matrix multiplication with $\mathbf{L}_{e}$ and $\mathbf{L}_{e^+}$), the cost of evaluating the hyperreduced test basis is:
\begin{equation}\label{eq:HROM_test}
\begin{split}
    W_{\text{test}} = & (\text{units to evaluate the $n_{e_{i}}$ Jacobian contributions from each element in Equation \ref{hyp_jac}}) + \\ 
    & (\text{units to sum each elemental contribution}) +\\
    & (\text{FLOPs to multiply the matrix resulting from the previous step and } \mathbf{V}) \\
    = & (n_{e_{i}}d_ed_{e^+}) + (n_{e_{i}}d_ed_{e^+}) + \\
    & (\text{FLOPs to multiply matrices of dimensions } (N \times N) \text{ and } (N \times n_i) \\ & \text{ where the first is sparse}) \\
    = & (2n_{e_{i}}d_ed_{e^+}) + (2n_{e_{i}}d_ed_{e^+}n_i) .
\end{split}
\end{equation}
where $n_{e_{i}}$ is the number of elements in the reduced mesh at a particular sampling cycle $i$. The hyperreduced residual can be written as:
\begin{equation}
    \mathbf{\hat{R}}^{(k)} \approx \tilde{\mathbf{R}}^{(k)} =  \sum_{e \in \tilde{\mathcal{E}}} \xi_e  \tilde{\mathbf{W}}^{(k)^{T}} \mathbf{L}_e^{T} \mathbf{R}_e^{(k)}.
 \end{equation}
Again neglecting the projection of the elemental residual into the FOM dimension, the work units to find the hyperreduced is:
\begin{equation}\label{eq:HROM_res}
\begin{split}
    W_{\text{residual}} = & (\text{units to evaluate the $n_e$ residual contributions from each element in Equation \ref{eq:hyp_res}}) +\\
    & (\text{FLOPs to multiply the matrix resulting from the previous step and } \tilde{\mathbf{W}}^{(k)}) + \\
    & (\text{units to sum each elemental contribution}) \\
    = & (n_{e_{i}}d_e) + (2n_{e_{i}}n_id_e) + (n_{e_{i}}n_i) \\
    = & n_{e_{i}}  (d_e + 2n_i d_e + n_i).
\end{split}
\end{equation}

The cost to assemble the LHS and RHS of Equation \ref{eq:HROM_non} is then:
\begin{equation}\label{eq:HROM_assem}
\begin{split}
    W_{\text{assemble}} = & (\text{FLOPs to multiply }\tilde{\mathbf{W}}^{(k)} \text{ and } \tilde{\mathbf{W}}^{(k)})  \\
    = &[n_i^2 (2N -1)].
\end{split}
\end{equation}

The work units, as shown in section \ref{sec:adapt}, for one non-linear iteration for the HROM solution is:
\begin{equation}
\begin{split}
    W_{\text{non lin}} = & W_{\text{residual}}  +  W_{\text{test}} + W_{\text{assemble}} + (\text{units to solve for solution update}) \\ 
    = & [n_{e_{i}}  (d_e + 2n_i d_e + n_i)] + [2n_{e_{i}} d_e d_e^+  +
        2 d_e d_e^+ n_{e_{i}} n_i] +  [n_i^2 (2N -1)] + (n_i^3).
\end{split}
 \end{equation}

\subsection{DWR Error $\epsilon_r$ Work Units}
Prior to incorporating hyperreduction, the adjoint problem for the second DWR error indicator between the coarse and fine ROM is:
\begin{equation}
\begin{split}
& \left[\mathbf{V}_h^T \frac{\partial \mathbf{R}}{\partial \tilde{\mathbf{w}}_h} \bigg|_{\tilde{\mathbf{w}}_H}^{T}  \frac{\partial \mathbf{R}}{\partial \tilde{\mathbf{w}}_h} \bigg|_{\tilde{\mathbf{w}}_H} \mathbf{V}_h \right]^T \hat{\psi}_h = -\left[ \frac{\partial \mathcal{J}}{\partial \tilde{\mathbf{w}}_h} \bigg|_{\tilde{\mathbf{w}}_H} \mathbf{V}_h\right]^T.
\end{split}
\end{equation}
The cost to assemble the LHS of the adjoint problem is then:
\begin{equation}\label{eq:ROM_LHS_adj}
\begin{split}
    W_{\text{LHS adj}} = & (\text{FLOPs to multiply } \mathbf{V}_h  \text{ and } \frac{\partial \mathbf{R}}{\partial \tilde{\mathbf{w}}_h} \bigg|_{\tilde{\mathbf{w}}_H}) + \\
    & (\text{FLOPs to multiply the matrix resulting from the previous step and } \frac{\partial \mathbf{R}}{\partial \tilde{\mathbf{w}}_h} \bigg|_{\tilde{\mathbf{w}}_H}^T + \\
    & (\text{FLOPs to multiply the matrix resulting from the previous step and } \mathbf{V}_h^T) \\
    = & (\text{FLOPs to multiply matrices of dimensions } (N \times N) \text{ and } (N \times n_i)) + \\
    & (\text{FLOPs to multiply matrices of dimensions } (N \times N) \text{ and } (N \times n_i)) + \\
    & (\text{FLOPs to multiply matrices of dimensions } (n_i \times N) \text{ and } (N \times n_i)) \\
    = & [Nn_i(2N-1)] + [Nn_i(2N-1)] + [n_i^2(2N-1)] \\
    = & (2Nn_i + n_i^2)(2N-1).
\end{split}
\end{equation}
The cost to assemble the RHS of the adjoint problem is:
\begin{equation}\label{eq:ROM_RHS_adj}
\begin{split}
    W_{\text{RHS adj}} = & (\text{FLOPs to multiply } \mathbf{V}_h  \text{ and } \frac{\partial \mathcal{J}}{\partial \tilde{\mathbf{w}}_h} \bigg|_{\tilde{\mathbf{w}}_H}) \\ 
    = & (\text{FLOPs to multiply matrices of dimensions } (1 \times N) \text{ and } (N \times n_i)) \\
    = & [n_i(2N-1)].
\end{split}
\end{equation}

The error estimate found using the adjoint from above is then:
\begin{equation}
\begin{split}
\epsilon_r = \hat{\psi}_h^T \hat{\mathbf{R}}_h(\tilde{\mathbf{w}}_H) =   \hat{\psi}_h^T \mathbf{V}_h^T \frac{\partial \mathbf{R}}{\partial \tilde{\mathbf{w}}_h} \bigg|_{\tilde{\mathbf{w}}_H}^{T} \mathbf{R}_h(\tilde{\mathbf{w}}_H) .
\end{split}
\end{equation}

The work units required to solve for $\epsilon_r$ can be broken down into:
\begin{equation}
\begin{split}
    W_{\epsilon_r} = & (\text{FLOPs to multiply } \mathbf{R}_h(\tilde{\mathbf{w}}_H)   \text{ and } \frac{\partial \mathbf{R}}{\partial \tilde{\mathbf{w}}_h} \bigg|_{\tilde{\mathbf{w}}_H}^{T}) + \\
    & (\text{FLOPs to multiply the matrix resulting from the previous step and } \mathbf{V}_h^T) + \\
    & (\text{FLOPs to multiply the matrix resulting from the previous step and } \hat{\psi}_h^T) \\
    = & (\text{FLOPs to multiply matrices of dimensions } (N \times N) \text{ and } (N \times 1)) + \\
    & (\text{FLOPs to multiply matrices of dimensions } (n_i \times N) \text{ and } (N \times 1)) + \\
    & (\text{FLOPs to multiply matrices of dimensions } (1 \times n_i) \text{ and } (n_i  \times 1)) \\
    = & [N(2N-1)] + [n_i(2N-1)] + (2n_i-1) \\
    = & (N + n_i)(2N-1) + (2n_i-1).
\end{split}
\end{equation}

Recalling the result in section \ref{sec:adapt}, the work units associated with solving the DWR at one ROM point is:
\begin{equation}
\begin{split}
    W_{\text{DWR error}} = & (\text{units to evaluate the residual}) + \\
    & (\text{units to evaluate the Jacobian}) + \\
    & (\text{units to evaluate the derivative of the functional w.r.t. $\tilde{\mathbf{w}}$}) + \\
    & W_{\text{LHS adj}} + W_{\text{RHS adj}} + \\
    & (\text{units to solve for the adjoint}) +  W_{\epsilon_r} \\
    =  & (N) + (N^2) + (N) + [(2Nn_i + n_i^2 + n_i)(2N-1)] + \\
    & (n_i^3) + [(N + n_i)(2N -1) + (2n_i -1)].
\end{split}
\end{equation}

\subsection{Hyperreduced DWR Error $\epsilon_r$ Work Units}
As discussed in section \ref{sec:adapt}, the updated error indicator with hyperreduction can be found using the following adjoint:
\begin{equation}
\left[ \left\{ \sum_{e\in\widetilde{\mathcal{E}}}\xi_e\mathbf{W}_h^T\mathbf{L}_e^T \frac{\partial \mathbf{R}_e}{\partial \tilde{\mathbf{w}}_h} \bigg|_{\tilde{\mathbf{w}}_H}  \right\}\mathbf{V}_h \right]^T  \tilde{\psi}_h = -\left[ \frac{\partial \mathcal{J}}{\partial \tilde{\mathbf{w}}_h} \bigg|_{\tilde{\mathbf{w}}_H} \mathbf{V}_h\right]^T
\end{equation}

The cost to assemble the LHS of the adjoint problem is then:
\begin{equation}
\begin{split}
    W_{\text{LHS adj}} = & (\text{FLOPs to evaluate quantity in the curly brackets}) + \\
    & (\text{FLOPs to multiply the matrix resulting from the previous step and }\mathbf{V}_h) \\
    = & (\text{FLOPs to evaluate quantity in the curly brackets}) + \\
    & (\text{FLOPs to multiply matrices of dimensions } (n_i \times N) \text{ and } (N \times n_i)) \\
    = & [n_{e_{i}}(d_e n_i + n_i)] + [n_i^2(2N-1)]. \\
\end{split}
\end{equation}
The cost to assemble the RHS of the adjoint problem is:
\begin{equation}
\begin{split}
    W_{\text{RHS adj}} = & (\text{FLOPs to multiply } \mathbf{V}_h  \text{ and } \frac{\partial \mathcal{J}}{\partial \tilde{\mathbf{w}}_h} \bigg|_{\tilde{\mathbf{w}}_H}) \\ 
    = & (\text{FLOPs to multiply matrices of dimensions } (1 \times N) \text{ and } (N \times n_i)) \\
    = & [n_i(2N-1)].
\end{split}
\end{equation}

The error indicator can then be found using:
\begin{equation}
\begin{split}
\epsilon_r &=  - \tilde{\psi}_h^T \tilde{\mathbf{R}}_h(\tilde{\mathbf{w}}_H) .
\end{split}
\end{equation}

The work units required to solve for $\epsilon_r$ can be broken down into:
\begin{equation}
\begin{split}
W_{\epsilon_r} = & (\text{FLOPs to multiply } \tilde{\mathbf{R}}_h(\tilde{\mathbf{w}}_H)  \text{ and } \tilde{\psi}_h^T) \\ 
    = & (\text{FLOPs to multiply matrices of dimensions } (1 \times n_i) \text{ and } (n_i \times 1)) \\
    = & (2n_i-1).
\end{split}
\end{equation}

Recalling the result in section \ref{sec:adapt}, the work units associated with solving the hyperreduced DWR at one ROM point is:
\begin{equation}
\begin{split}
    W_{\text{DWR error}} = & (\text{units to evaluate the hyperreduced residual}) + \\
    & (\text{units to evaluate the hyperreduced test basis}) + \\
    & (\text{units to evaluate the derivative of the functional w.r.t. $\tilde{\mathbf{w}}$}) + \\
    & W_{\text{LHS adj}} + W_{\text{RHS adj}} +\\
    & (\text{units to solve for the adjoint}) + W_{\epsilon_r} \\
    =  & [n_{e_{i}}  (d_e + 2n_i d_e + n_i)] + [n_{e_{i}} d_e d_e^+  +
        2 d_e d_e^+ n_{e_{i}} n_i] + (N) \\
        & +  [n_{e_{i}}(d_e n_i + n_i) + 
        (n_i^2)(2N-1) + n_i(2N-1)] + (n_i^3) + (2n_i -1).
\end{split}
\end{equation}